\documentclass[12]{amsart}
\usepackage{graphicx}
\usepackage{amssymb}
\usepackage{epstopdf}
\usepackage{nicefrac}
\usepackage{float}
\usepackage{color}
\usepackage{tabularx}
\usepackage{hyperref}

\newtheorem{thm}{THEOREM}[section]

\newtheorem{cor}[thm]{COROLLARY}
\newtheorem{defn}[thm]{DEFINITION}

\newtheorem{lemma}[thm]{LEMMA}
\newtheorem{prob}[thm]{PROBLEM}
\newtheorem{prop}[thm]{PROPOSITION}
\newtheorem{quest}[thm]{QUESTION}

\newcommand{\ds}{\displaystyle}

\newcommand{\EF}{{\mathbf{E}_{\F}}}
\newcommand{\PF}{{\mathbf{P}_{\F}}}
\newcommand{\HF}{{\mathbf{H}_{\F}}}

\newcommand{\ET}{{\mathbf{E}_{\cT}}}
\newcommand{\PT}{{\mathbf{P}_{\cT}}}
\newcommand{\HT}{{\mathbf{H}_{\cT}}}

\newcommand{\RF}{{\mathcal{R}_{\F}}}

\newcommand{\F}{{\mathcal F}} 
\newcommand{\G}{\Gamma}

\newcommand{\e}{{\epsilon}} 

\newcommand{\cGF}{\cG_{\F}} 
\newcommand{\GF}{\Gamma_{\F}} 

\newcommand{\wtU}{{\widetilde U}}

\newcommand{\wtvp}{{\widetilde{\varphi}}}

\newcommand{\whw}{{\widehat w}}
\newcommand{\whx}{{\widehat x}}
\newcommand{\why}{{\widehat y}}
\newcommand{\whz}{{\widehat z}}

 \newcommand{\whF}{{\widehat{\F}}}
\newcommand{\whM}{{\widehat{M}}}
\newcommand{\whQ}{{\widehat{Q}}}

\newcommand{\oU}{{\overline{U}}}

\newcommand{\bZ}{{\bf Z}}
\newcommand{\bC}{{\bf C}}

\newcommand{\mC}{{\mathbb C}}
\newcommand{\mD}{{\mathbb D}}

\newcommand{\mG}{{\mathbb G}}

\newcommand{\mN}{{\mathbb N}}

\newcommand{\mR}{{\mathbb R}}
\newcommand{\mS}{{\mathbb S}}
\newcommand{\mT}{{\mathbb T}}
\newcommand{\mZ}{{\mathbb Z}}

\newcommand{\cE}{{\mathcal E}}
\newcommand{\cG}{{\mathcal G}}
\newcommand{\cH}{{\mathcal H}}

\newcommand{\cO}{{\mathcal O}}
\newcommand{\cP}{{\mathcal P}}
\newcommand{\cQ}{{\mathcal Q}}
\newcommand{\cR}{{\mathcal R}}
\newcommand{\cS}{{\mathcal S}}
\newcommand{\cT}{{\mathcal T}}
\newcommand{\cU}{{\mathcal U}}

\newcommand{\cX}{{\mathcal X}}
\newcommand{\cZ}{{\mathcal Z}}

\newcommand{\fH}{\mathfrak{H}}

\newcommand{\fM}{{\mathfrak{M}}}

\newcommand{\fT}{{\mathfrak{T}}}

\newcommand{\vp}{{\varphi}}

\begin{document}

\title{Lectures on Foliation Dynamics: Barcelona 2010}

\thanks{2000 {\it Mathematics Subject Classification}. Primary 57R30, 37C55, 37B45; Secondary 53C12 }

\author{Steven Hurder}
\address{Steven Hurder, Department of Mathematics, University of Illinois at Chicago, 322 SEO (m/c 249), 851 S. Morgan Street, Chicago, IL 60607-7045}
\email{hurder@uic.edu}
\thanks{Preprint date: October 24, 2011; final edits: August 25, 2014}

\date{}


\keywords{}

\maketitle
 
\tableofcontents

 \section{Introduction}

 The study of foliation dynamics seeks to understand the asymptotic properties of   leaves of foliated manifolds, their statistical properties  such as orbit growth rates and geometric entropy, and to classify geometric and topological ``structures'' which are associated to the dynamics, such as the minimal sets of the foliation. The study is inspired by the seminal work of 
  Smale \cite{Smale1967} (see also the comments by Anosov \cite{Anosov2006}) outlining a program of study for the differentiable dynamics for a $C^r$-diffeomorphism $f \colon N \to N$ of a closed manifold $N$, $r \geq 1$. The themes of this approach included:  
\begin{enumerate}
\item  Classify dynamics as hyperbolic, or otherwise;
\item Describe the minimal/transitive closed invariant sets and attractors;
\item Understand when the system is structurally stable under $C^r$-perturbations, for $r \geq 1$; 
\item Find   invariants (such as cohomology, entropy or zeta functions) characterizing the system.
\end{enumerate}
   Smale also suggested to study these topics for large group actions, which leads directly to the topics of these notes. 
       The study of the dynamics of   foliations began in ernest   in the 1970's with the research programs of    Georges Reeb,  Stephen Smale, Itiro Tamura, and their   students.

A strict analogy between foliation dynamics and the theory for diffeomorphisms cannot be exact.
Perhaps the most fundamental problem is the role played by invariant probability measures in the analysis of dynamics of diffeomorphisms. A diffeomorphism of a compact manifold generates an action of the group of integers $\mZ$, which is amenable, so every minimal set carries at least one invariant measure. Many of the techniques of smooth dynamics use such invariant measures to analyze and  approximate the ``typical dynamics'' of the map. In contrast, a foliation need not have any transverse, holonomy-invariant measures. Moreover, the dynamics of foliations which do not admit such invariant measures provide some of the most important examples in the subject.

Even when such invariant measures exist, there is the additional problem with ``time''. In foliation dynamics,  the concept of linearly or time-ordered trajectories is replaced with the vague notion of multi-dimensional futures for points, as defined by the leaves through the points. The geometry of the leaves thus plays a fundamental role in the study of foliation dynamics, which is a fundamentally new aspect of the subject, in contrast to the study of diffeomorphisms, or $\mZ$-actions. 

Issues with other basic concepts also arise, such as the existence of periodic orbits, which for foliations corresponds   most precisely to compact leaves. However, analogs of hyperbolicity almost never imply the existence of compact leaves, while this is a fundamental tool for the study of diffeomorphisms. In spite of these obstacles, there is a robust theory   of foliation dynamics. 

Another aspect of foliation dynamics, is that the collection of examples illustrating ``typical behavior'' is woefully incomplete. There is a vast richness of dynamical behaviors  for foliations, much greater than for flows and diffeomorphisms, yet the constructions of examples to illustrate these behaviors is still very incomplete. We will highlight in these notes some examples of a more novel nature, with the caveat that those presented are far from being close to a complete set of representatives. There is much work to be done!  The following are some of the  topics  we  discuss in these notes:
\begin{enumerate}
\item The asymptotic properties of leaves of $\F$ 
\begin{itemize}
\item How do the leaves accumulate onto the minimal sets?
\item What are the topological types of minimal sets? Are they ``manifold-like''?
\item Invariant measures: can you quantify the rates of recurrence of leaves?
\end{itemize}

\item Directions of ``stability'' and ``instability'' of leaves 
\begin{itemize}
\item Exponents: are there directions of exponential divergence?
\item Stable manifolds:  show the existence of  dynamically defined transverse invariant manifolds, and how do they influence the global behavior of leaves?
\end{itemize}

\item Quantifying chaos 
\begin{itemize}
\item Define a measure of transverse chaos -- foliation entropy
\item Estimate the entropy using linear approximations
\end{itemize}

\item Dynamics of minimal sets 
\begin{itemize}
\item Hyperbolic   minimal sets
\item Parabolic    minimal sets
\end{itemize}

\item Shape of minimal sets 
\begin{itemize}
\item Matchbox manifolds
\item Approximating minimal sets
\item Algebraic invariants
\end{itemize}

\end{enumerate}

The subject of foliation dynamics is very broad, and includes many other topics to study beyond what is discussed in these notes, such as rigidity of the dynamical system defined by the leaves, the behavior of random walks on leaves and properties of harmonic measures, and the Hausdorff dimension of minimal sets, to name a few additional important ones.

  This survey is based on a series of  five lectures, given  May 3--7, 2010,   at the {\it Centre de Recerca Matem\`{a}tica, Barcelona}. The goal of the lectures was to present  aspects of the theory of foliation  dynamics  which have particular importance  for the   classification of foliations of compact manifolds.  The lectures emphasized intuitive concepts and informal discussion, as can be seen from the slides    \cite{Hurder2010a}. Due to its origins, these notes will eschew formal definitions when convenient, and the reader is referred to the   sources 
  \cite{CN1985,CandelConlon2000,Hurder2011b,Lawson1977,MS2006,Tamura1992,Walczak2004} for further details. 
  
  Many of the illustrations in the following text were drawn by   Lawrence Conlon, circa 1994. Our thanks for his permission to use them.

  The author would like to sincerely thank the organizers of this workshop, Jes\'{u}s \'{A}lvarez L\'{o}pez (Universidade de Santiago de Compostela) and Marcel Nicolau (Universitat Aut\`{o}noma de Barcelona) for their efforts to make this month long event happen, and the C.R.M. for the excellent hospitality offered to the participants.

  \section{Foliation Basics}\label{sec-basics}

\begin{figure}[H]
\setlength{\unitlength}{1cm}
\begin{picture}(20,4.7)
\put(5,.2){\includegraphics[width=50mm]{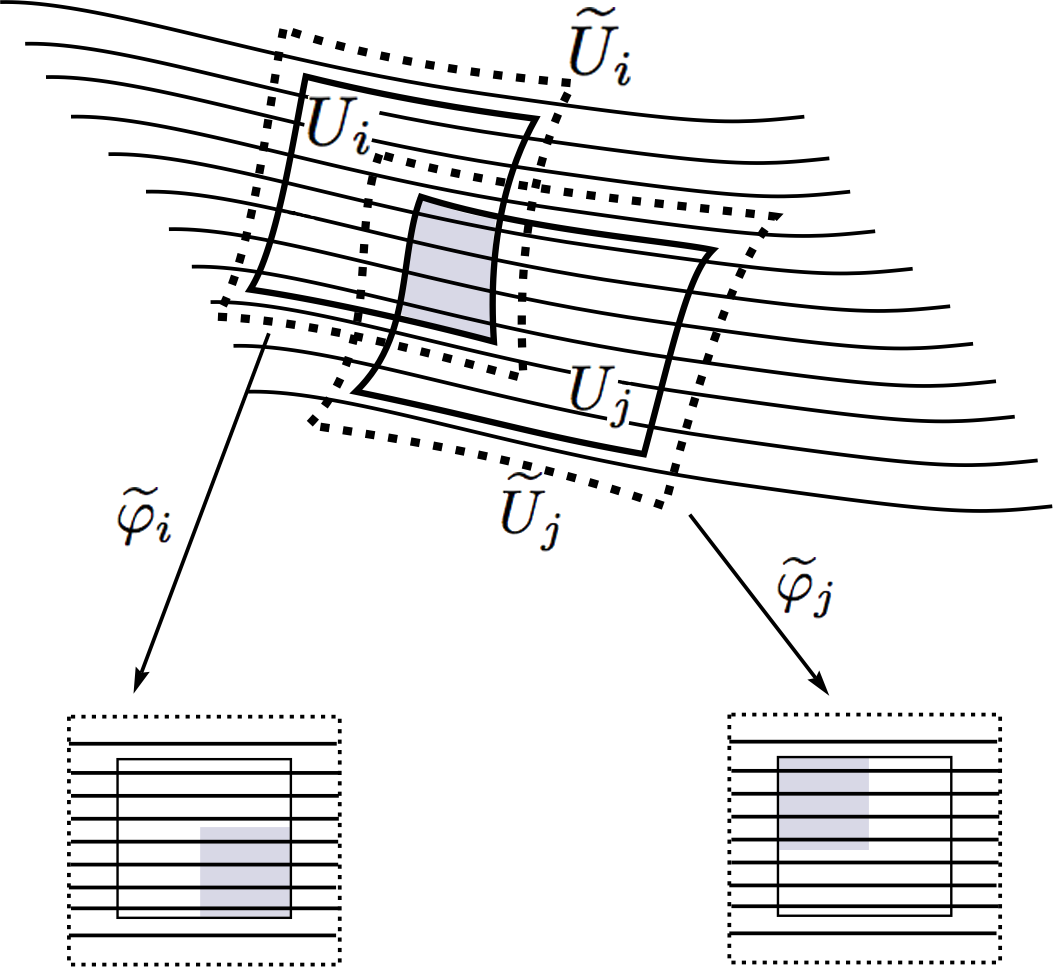}} 
\end{picture}
\caption{Foliation charts}
\label{Figure1}
\vspace{-20pt}
\end{figure}

 A foliation $\mathcal F$ of dimension $p$ on a manifold $M^m$ is a 
 decomposition  into ``uniform layers'' -- the leaves -- which are immersed submanifolds of codimension $q = m -p$:  there is an open covering of $M$ by   coordinate charts 
so that the leaves are mapped into linear planes of dimension $p$, and the transition functions preserve these planes.
See Figure~\ref{Figure1}.

More precisely,  we   require that around each point $x$ there is an open neighborhood $x \in U_x$ and a  ``foliation chart'' 
$\vp_x \colon U_x \to (-1,1)^m$ for which each inverse image $\cP_x(y) = \vp_x^{-1}( (-1,1)^p \times \{y\}) \subset U_x$, $y \in (-1,1)^q$, is a connected component of $L \cap U_x$ for some leaf $L$.
The foliation $\F$ is said to have differentiability class $C^r$, $0 \leq r \leq \omega$, if the charts $\vp_x$ can be chosen to be $C^r$-coordinate charts for the manifold $M$. For a compact manifold $M$, we can always choose a finite covering of $M$ by foliation charts 
$\ds \cU \equiv \{ (\vp_i , U_i) \mid i = 1, 2, \ldots, k\}$ with the additional   property that each chart $\vp_i$ admits an extension to a foliation chart $\wtvp_i \colon \wtU_i \to (-2,2)^m$ where the closure $\oU_i \subset \wtU_i$.

The subject of foliations tends to be quite abstract, as it is difficult to illustrate in full the implications of the above definition    in   dimensions greater than two. 
 One is typically presented with    a few ``standard examples'' in dimensions two and three, that   hopefully   yield  intuitive insight from which to gain a deeper understanding of the more general cases.   
For example, many talks  with ``foliations'' in the title start with the following example,
 the 2-torus  $\mT^2$ foliated by lines of irrational slope:

\begin{figure}[H]
\setlength{\unitlength}{1cm}
\begin{picture}(20,2.2)
\put(5.6,0){\includegraphics[width=46mm]{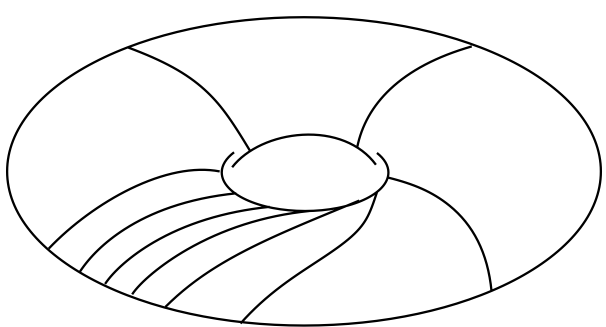}}
\end{picture}
\caption{Linear foliation with all leaves dense}
\label{Figure2}
\vspace{-20pt}
\end{figure} 

Never trust a talk which starts with this example! It is just too simple, in that the leaves are parallel and contractible, hence the foliation  has no germinal holonomy. Also, every leaf of $\F$ is uniformly dense in $\mT^2$ so the topological nature of the minimal sets for $\F$ is trivial to determine. The key dynamical information about this example is given by the rates of returns to open subsets, which is more analytical than topological information.

At the other extreme of examples of  foliations defined by flows on a surface are those with a compact leaf as the unique minimal set, such  as in Figure~\ref{Figure3}:

\begin{figure}[H]
\setlength{\unitlength}{1cm}
\begin{picture}(20,2.2)
\put(5.6,.2){\includegraphics[width=46mm]{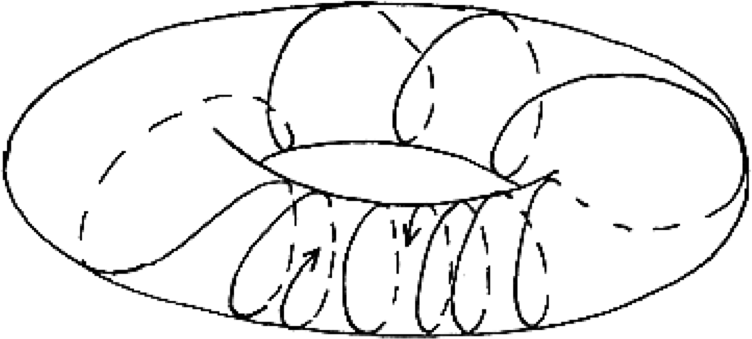}}
\end{picture}
\caption{Flow with one attracting leaf}
\label{Figure3}
\vspace{-20pt}
\end{figure}

Every orbit limits to the circle, which is the  forward (and backward) limit set for all leaves. 
 
One other canonical example is that of the Reeb foliation of the solid $3$-torus as pictured in Figure~\ref{Figure4}, which has a similar  dynamical description:

\begin{figure}[H]
\setlength{\unitlength}{1cm}
\begin{picture}(20,3.2)
\put(5.6,.2){\includegraphics[width=50mm]{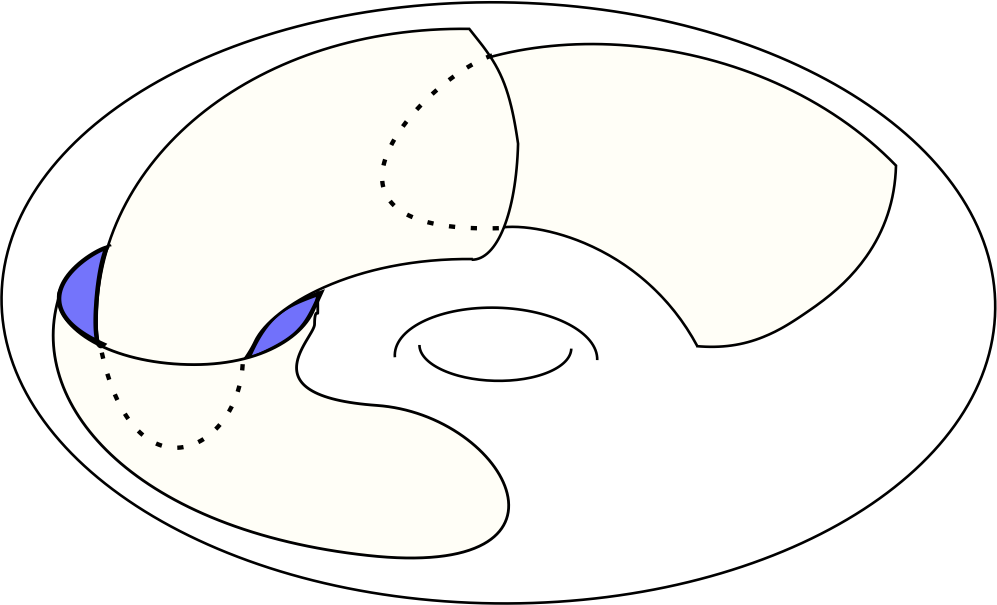}}
\end{picture}
\caption{Reeb foliation of solid torus}
\label{Figure4}
\vspace{-20pt}
\end{figure}

This example  illustrates several concepts: the limit sets of leaves, the existence of attracting holonomy for  the compact toral leaf, and also the (possible) existence of multiple hyperbolic measures for the foliation geodesic flow, as in Definition~\ref{def-thm}.  

 We will introduce further examples in the text, that illustrate more advanced dynamical properties of foliations. Although, as mentioned above, it becomes more difficult to illustrate concepts that only arise for foliations of manifolds of more complicated 3-manifolds, or in higher dimensions. The interested reader should view   the illustrations in the beautiful article by \'{E}tienne Ghys and Jos Ley for flows on 3-manifolds \cite{GhysLey2006} to get some intuitive insights  of the complexity that is ``normal'' for foliation dynamics  in higher dimensions.

 \section{Topological Dynamics}\label{sec-topdyn}

The study of the topological dynamics   for continuous   actions of non-compact groups on compact spaces is a venerable topic, as in 
\cite{Auslander1963,Auslander1988, Furstenberg1963,Ellis1969, Veech1970}, or the more   recent works       \cite{AAG2008,AA2003,Lindenstrauss1999}. 
The holonomy along leafwise paths of $\F$ defines local homeomorphisms between open subsets of $\mR^q$, and many    concepts of topological dynamics adapt   to this pseudogroup context.

Recall   the concept of holonomy pseudogroup for a foliation. The point of view we adopt is best illustrated by starting with the classical case of flows.   Recall for a non-singular  flow $\varphi_t \colon M \to M$ the orbits define a $1$-dimensional foliation $\F$, whose leaves are the orbits of points.
 
 Choose a cross-section $\cT \subset M$ which is transversal to the orbits, and intersects each orbit (so $\cT$ need not be connected.) Then for each $x \in \cT$ there is some least $\tau_x > 0$ so that $\varphi_{\tau_x} (x) \in \cT$. 
 The positive constant $\tau_x$ is called the   \emph{return time} for $x$. See the illustration Figure~\ref{Figure5} below.

\begin{figure}[H]
\setlength{\unitlength}{1cm}
\begin{picture}(20,3.8)
\put(5,.2){\includegraphics[width=56mm]{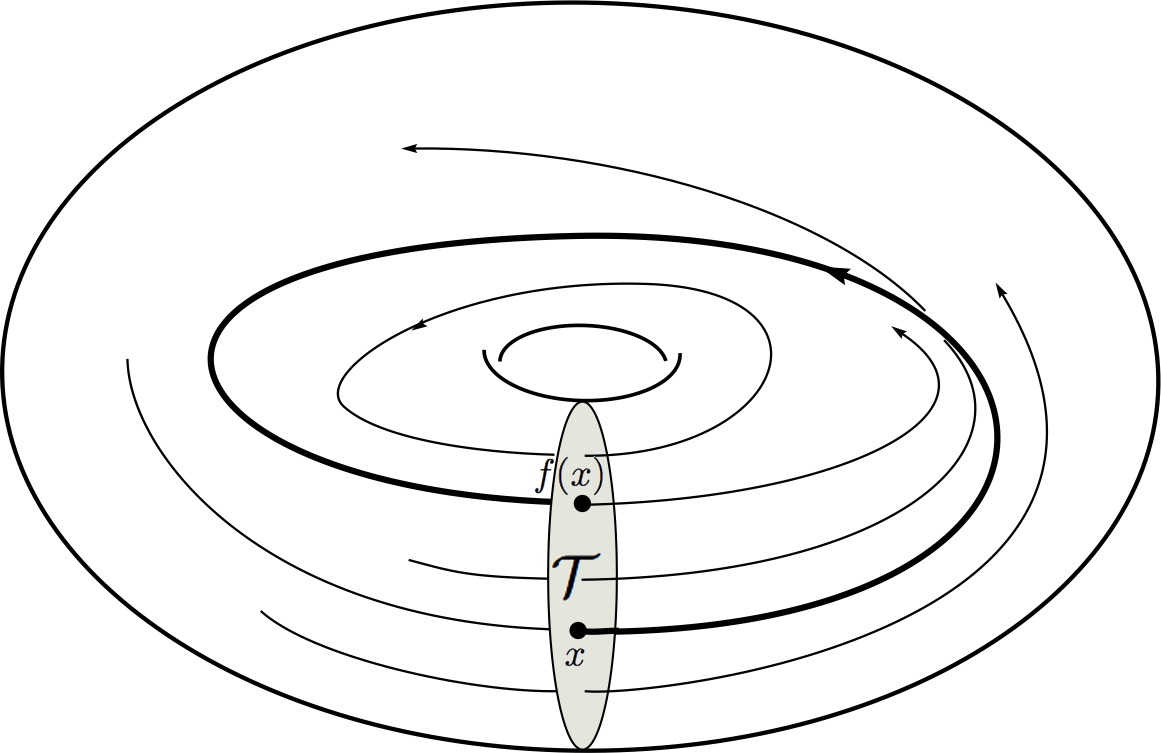}}
\end{picture}
\caption{Cross-section to a flow}
\label{Figure5}
\vspace{-20pt}
\end{figure}

 The induced map $f(x) = \varphi_{\tau_x} (x)$ is a \emph{Borel map} $f \colon \cT \to \cT$,  
  called the \emph{holonomy} of the flow. The choice of a cross-section for a flow reduces the study of its dynamical properties to that of the discrete dynamical system  $f \colon \cT \to \cT$.

The holonomy for foliations is defined similarly to the case for flows, as local $C^r$-diffeomorphisms  associated to paths along leaves, starting and ending at a fixed transversal, except that there is a fundamental difference. For the orbit of a flow $L_w$ through a point $w$, there exists two choices of   trajectory along a unit speed path, either forward and backward. 
 However, for a leaf  $L_w$ of a foliation $\F$ of dimension at least two, there is  no such concept as ``forward'' or ``backward'', and all directions   yield paths along which one may discover dynamical properties of the foliation. 
The correct analog   is thus the holonomy pseudogroup $\cGF$ construction,    introduced by Haefliger \cite{Haefliger1958, Haefliger1970, Haefliger1984}.

 We fix the following conventions. 
$M$ is a compact Riemannian manifold without boundary, and 
$\F$ is a codimension $q$-foliation, transversally $C^r$ for $1 \leq r \leq \infty$. Fix also a finite covering by foliation charts, 
$\ds \cU \equiv \{ (\vp_i , U_i) \mid i = 1, 2, \ldots, k\}$. 
  The projections along plaques in each chart defines submersions $\phi_i \colon U_i \to (-1,1)^q$. When $U_i \cap U_j \ne \emptyset$ we say that the pair $(i,j)$ is \emph{admissible}, and can define the transition function
  $\gamma_{i,j} \colon T_{i,j} \to T_{j,i}$ where $T_{i,j} = \phi_i(U_i \cap U_j) \subset T_i = (-1,1)^q$.

The finite collection of local diffeomorphisms $\cGF^{(1)} \equiv \{\gamma_{i,j} \mid (i,j) ~ {\rm admissible}\}$ generates a pseudogroup $\cGF$ of local $C^r$-diffeomorphisms modeled on the transverse space $\mR^q$. The choice of $\cU$ yielding the collection $\cGF^{(1)}$ is analogous to the notion of a generating set for a   group $\G$.

Now assume, without loss of generality, that the submanifolds $\cT_i = \vp_i^{-1}(\{0\} \times T_i) \subset U_i$ have pairwise disjoint closures, so for each $x \in \cT = \cT_1 \cup \cdots \cup \cT_k$ there is a unique $i$ with $x \in \cT_i$. Also, we assume that there is given a metric on $M$, which restricts via the embedding $\cT \subset M$ to a metric $d_{\cT}$ on $\cT$.
 
\begin{figure}[H]
\setlength{\unitlength}{1cm}
\begin{picture}(20,3.8)
\put(4,.2){\includegraphics[width=76mm]{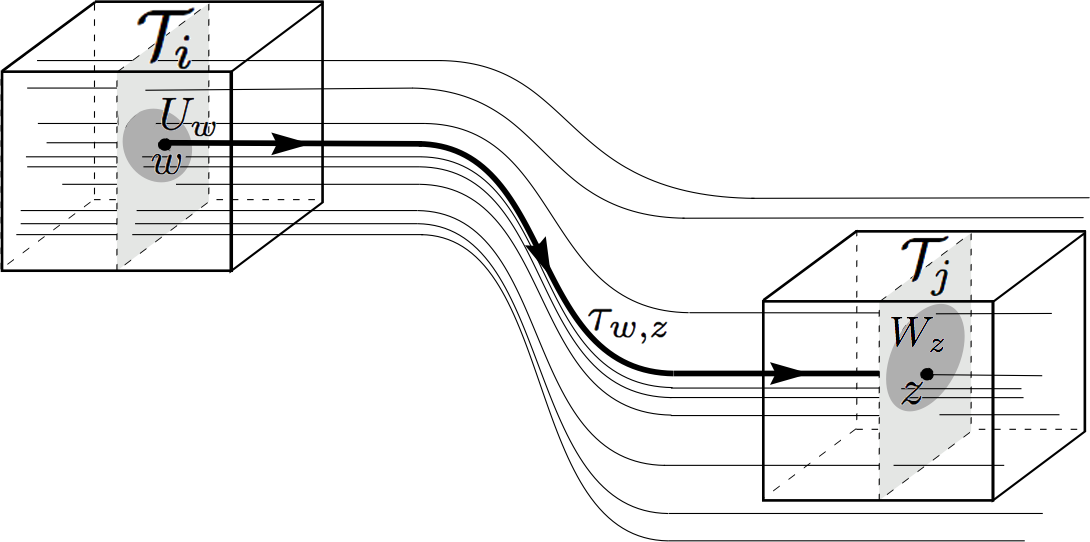}}
\end{picture}
\caption{Holonomy along a leafwise path}
\label{Figure6}
\vspace{-20pt}
\end{figure} 

 Fix $w \in \cT$, then  
choose $z \in L_w \cap \cT$ and a smooth path $\tau_{w,z} \colon [0,1] \to L_w$. 
 Cover the path $\tau_{w,z}$ by foliation charts, which determines a plaque chain from $w$ to $z$ which contains the path  $\tau_{w,z}$.  Then there exists an open subset $w \in U_w \subset \cT_i$ such that every $w' \in U_w$ admits a plaque chain that shadows the one along $\tau_{w,z}$ and so defines an image point $h_{\tau_{w,z}}(w') \in  W_z \subset \cT_j$ for some $j$.   This defines a local homeomorphism  $h_{\tau_{w,z}} \colon U_w \to W_z$ of open subsets of $\mR^q$, and hence  the holonomy pseudogroup $\cGF$ for $\F$ modeled on $\cT$, which is compactly generated in sense of Haefliger \cite{Haefliger2002a}. 
 This most basic concept of foliation theory is developed in detail in all standard texts \cite{CN1985, CandelConlon2000, CandelConlon2003, HecHir1981, Tamura1992}.

We introduce a notational convention that is quite convenient. For a   leafwise path  $\tau_{w,z} \colon [0,1] \to L_w$ with $w,z \in \cT$, let $Dom(h_{\tau_{w,z}}) \subset \cT$ denote the largest domain of definition for $h_{\tau_{w,z}}$ obtained from some covering of the path $\tau_{w,z}$ by foliation charts. Note that this is an abuse of notation, as the domain is well defined, once the covering is chosen, but different coverings may yield distinct maximal domains, although all such consist  of open sets in $\cT$ containing the initial point $w$.

To eliminate the issue of domains, one  introduces the germ at $w$ of the local homeomorphism $h_{\tau_{w,z}}$,  denoted by $[h_{\tau_{w,z}}]_w$. 
The collection of all such germs $\{[h_{\tau_{w,z}}]_w   \mid  w \in \cT,  z \in L_w \cap \cT \}$ 
generates the \emph{holonomy groupoid}, denoted by  $\GF$.
We summarize these properties as follows:

\begin{prop} Let $\F$ be a foliation of a manifold $M$. Then 
\begin{enumerate}
\item $[h_{\tau_{w,z}}]_w$ depends only on the leafwise homotopy class of the path, relative endpoints.
\item  The maximal sizes of the domain $U_w$ and range $W_z$ representing an equivalence class $[h_{\tau_{w,z}}]_w$ depends on the path $\tau_{w,z}$.
\item The collection of all such maps $\{ h_{\tau_{w,z}} \colon  U_w \to  W_z\ \mid  w \in \cT,  z \in L_w \cap \cT \}$ generates the \emph{holonomy pseudogroup} $\cGF$.
\end{enumerate}
\end{prop}
Assume that $\F$ is a $C^1$-foliation of a compact Riemannian manifold, with smoothly immersed leaves. Then  for each leaf $L_w$ of $\F$ the induced Riemannian metric on $L_w$ is complete, so there exists a length minimizing geodesic in each homotopy class, modulo endpoints, of a path in $L_w$.
\begin{cor} 
Given a leafwise path  $\tau_{w,z} \colon [0,1] \to L_w$ let  $\sigma_{w,z} \colon [0,1] \to L_w$ be a leafwise geodesic segment which is   homotopic relative the  endpoints to $\tau_{w,z}$. Then $[h_{\tau_{w,z}}]_w = [h_{\sigma_{w,z}}]_w$.
\end{cor}

 While the germ $\gamma$ of the holonomy along a leafwise path $\tau_{w,z}$ is well defined, up to conjugation, the size of the domain of a representative map $h_{\tau_{w,z}} \in \cGF$ need not be. It is a strong restriction on the dynamics of $\cGF$ or the topology of $M$  that a uniform estimate on the sizes of the domains exists. This is a very delicate technical point that arises in many proofs about the dynamics of a foliation. Our choice of notation for the domains of the holonomy maps  suppresses this technical issue, for the purpose of simplicity of exposition.

 In the study of the topological dynamics of group actions, the   domains of definition  for the transformations are always well defined. On the other hand, in the following formulations for pseudogroups, we are careful to specify that the behavior is with respect to domains of the holonomy transformations associated to  leafwise paths.  Thus, while we will say that the groupoid $\cGF$   has a particular  dynamical property, more precisely this is with respect to the subcollection of transformations in $\cGF$ defined  geometrically by the holonomy parallel transport. 

First,    recall that    a  \emph{minimal set} for $\F$ is a closed, saturated subset   $\cZ \subset M$ for which every leaf $L \subset \cZ$ is dense.
In terms of pseudogroup $\cGF$, a subset $\cX \subset \cT$ is minimal if it is invariant under the action of $\cGF$ and every orbit is dense.

  A related notion is that of a   \emph{transitive set} for $\F$, which is a closed saturated subset   $\cZ \subset M$  such that there exists at least one dense leaf $L_0 \subset \cZ$. In other words, these are the subsets of a foliated manifold which are the closure of a single leaf.

The   concept  of an \emph{equicontinuous action}  is   classical for the dynamics of group actions. 

\begin{defn} \label{def-equic}
The dynamics of $\F$ restricted to a saturated subset $\cZ \subset M$ is \emph{equicontinuous} if for all $\epsilon > 0$, there exists $\delta > 0$ such that for all $w \ne w' \in \cZ \cap \cT$, and for any leafwise path $\tau_{w,z} \colon [0,1] \to L_w$ starting at $w$ and ending at some $z \in \cZ \cap \cT$ with $w,w' \in Dom(h_{\tau_{w,z}})$, then 
$\ds d_{\cT}(w,w') < \delta$ implies  $ d_{\cT}( h_{\tau_{w,z}}(w) , h_{\tau_{w,z}}(w')) < \epsilon$.
\end{defn}
The typical example is provided by the foliation defined by a flow with dense orbits on the $2$-torus given at the start of these notes. This is a special case of a Riemannian foliation, which admits a transverse metric so that all the holonomy maps $h_{\tau_{w,z}}$ are isometries.
We note that equicontinuity is a strong hypothesis on a pseudogroup. In particular, we have:

\begin{thm}[Sacksteder \cite{Sacksteder1965}] \label{thm-sacksteder} If $\cGF$ is an equicontinuous pseudogroup acting on a compact Polish space $\cX$, then there exists a Borel probability measure $\mu$ on $\cX$ which is $\cGF$-invariant. \hfill $\Box$
\end{thm}

The   concept of a  \emph{distal action} is  closely related to the above. 

\begin{defn} \label{def-distal}
The dynamics of $\F$ restricted to a  saturated subset $\cZ \subset M$ is \emph{distal} if for all $w \ne w'  \in \cX = \cZ \cap \cT$, there exists $\delta_{w,w'} > 0$  so that for any leafwise path $\tau_{w,z} \colon [0,1] \to L_w$ starting at $w$ and ending at some $z \in \cX$, with $w,w' \in Dom(h_{\tau_{w,z}})$,  then $\ds  d_{\cT}( h_{\tau_{w,z}}(w) , h_{\tau_{w,z}}(w'))  \geq \delta_{w,w'}$. 
\end{defn}
In other words, the metric distortion of the distance between two points in $\cX \subset \cT$  under the action of leafwise holonomy transformations is bounded from  below. The typical examples are provided by the foliations defined by an action of a parabolic subgroup on a compact quotient of a nilpotent Lie group by a lattice subgroup.
Distal and equicontinuous pseudogroups are closely related \cite{AAG2008,EEN2001,Furstenberg1963,Lindenstrauss1999,Veech1970}.

The   concept is of a   \emph{proximal action}  is   opposite   to that of a  distal action.   

\begin{defn} \label{def-proximal}
The dynamics of $\F$ restricted to a  saturated subset $\cZ \subset M$ is \emph{proximal}  if there exists $\epsilon > 0$ so that for all $w \ne w' \in \cZ \cap \cT_i$ for some $1 \leq i \leq k$ with $d_{\cT}(w,w') < \epsilon$, then for all $0 < \delta \leq \epsilon$ there exists a leafwise path $\tau_{w,z} \colon [0,1] \to L_w$ starting at $w$ and ending at some $z \in \cZ \cap \cT$ with $w,w' \in Dom(h_{\tau_{w,z}})$ such that  
$\ds   d_{\cT}( h_{\tau_{w,z}}(w) , h_{\tau_{w,z}}(w')) \leq \delta$.
\end{defn}
Proximality asserts for any   pair  of       points that are sufficiently close, there is a   holonomy map  for which  the distance between their images can be made arbitrarily close.   The typical examples are provided by the foliations defined by an action of a Borel subgroup on a compact quotient of a simple Lie group by a lattice subgroup. Typical examples of this special case are the weak-stable foliations  associated to the geodesic flows of  compact hyperbolic manifolds.

Finally, there is the fundamental   concept   of an  \emph{expansive action}.   

\begin{defn} \label{def-expansive}
The dynamics of $\F$ restricted to a  saturated subset $\cZ \subset M$ is \emph{expansive}, or more precisely $\epsilon$-expansive,  if there exists $\epsilon > 0$ so that for all $w \ne w' \in \cZ \cap \cT_i$ for some $1 \leq i \leq k$ with $d_{\cT}(w,w') < \epsilon$, then   there exists a leafwise path $\tau_{w,z} \colon [0,1] \to L_w$ starting at $w$ and ending at some $z \in \cZ \cap \cT$ with $w,w' \in Dom(h_{\tau_{w,z}})$ such that  
$\ds   d_{\cT}( h_{\tau_{w,z}}(w) , h_{\tau_{w,z}}(w')) \geq \epsilon$.
\end{defn}

 The simplest approach to classifying foliation topological dynamics, is to ask if a given closed invariant set $\cZ \subset M$ is either equicontinuous, distal, proximal, or expansive. There are many interesting examples of foliations with invariant sets exhibiting each of these dynamics.

 \section{Derivatives} \label{sec-derivatives}

The properties of foliation dynamics introduced above have been topological in nature. However, it has been known at least since the discovery of the Denjoy-type examples \cite{Denjoy1932} that the   topological dynamics of   flows, and more generally foliations,   are strongly influenced and restricted by the degree of differentiability of its holonomy maps.   A deeper understanding of foliation dynamics necessarily proceeds with a more detailed study of the differential properties of the holonomy pseudogroups.

  To begin, 
  we introduce the {\it transverse differentials} for the holonomy groupoid. 
 Consider first the case of a foliation $\F$ defined by a   smooth flow 
  $\varphi \colon \mR \times M   \to M$ generated by a non-vanishing  vector field $\vec{X}$.
Then   $T\F =  \langle \vec{X} \rangle \subset TM$.

For  $z = \varphi_t(w)$, consider the Jacobian matrix $D\varphi_t \colon  T_wM  \to T_{z}M$. The flow satisfies the 
group law ~ $\varphi_s \circ \varphi_t = \varphi_{s+t}$, which  implies the identity 
$D\varphi_s (\vec{X}_w) = \vec{X}_z$    by the chain rule for derivatives.
Introduce the normal bundle to the flow $Q  = TM / T\F$. For  each $w \in M$,  we identify $Q_w =  T_w\F^{\perp}$.
Thus,   $Q$ can be considered as  a subbundle of $TM$, and thereby  the   Riemannian metric on $TM$ induces metrics on each fiber $Q_w  \subset T_wM$. 
The derivative transformation preserves the normal bundle $Q \to M$, so  defines the \emph{normal derivative cocycle},
$$D\varphi_t \colon Q_w \longrightarrow Q_{z} \quad , \quad t \in \mR $$
We can then define the  norms of the normal derivative maps, 
$$\| D\varphi_t   \| = \| D\varphi_t \colon Q_w \longrightarrow Q_{z}\|$$
It is also useful to introduce the symmetric norm
$$\| D\varphi_t |_w  \|^{\pm} = \max \left\{ \| D\varphi_t \colon Q_w \longrightarrow Q_{z}\| , \| D\varphi_t ^{-1} \colon Q_z \longrightarrow Q_{w}\| \right\}$$
For $M$ compact and $t$ fixed, the norms $\| D\varphi_t  \mid_w  \|^{\pm}$ are uniformly bounded  for $w \in M$.

The maps $D\varphi_t \colon Q_w \longrightarrow Q_{z}$ can be thought of as ``non-autonomous local approximations'' to the transverse behavior of the flow $\varphi_t$. The actual values of these derivatives is only well defined up to a global choice of framing of the normal bundle $Q$, so extracting useful dynamical information from transverse  derivatives presents a challenge. One solution to this problem was solved by seminal work of Pesin in the 1970's.
 ``Pesin Theory'' is a collection of results about the dynamical properties of flows, based on defining non-autonomous linear approximations of the normal behavior to the flow. Excellent discussions and references for this theory are in these references \cite{Oseledets1968,Mane1987,BarreiraPesin2002}.  We use only a small amount of the full Pesin theory in the discussion in these notes. 
 
 First, let us recall a basic fact for the dynamics induced by a linear map. Given a matrix $A \in GL(q,\mR)$, let $L_A \colon \mR^n \to \mR^n$ be the linear map defined by multiplication by $A$. We say the action $L_A$ is \emph{partially hyperbolic} if $A$ has an eigenvalue of norm not equal to $1$. In this case, there is an eigenspace for $A$ which is defined dynamically as the direction of maximum rate of expansion (or minimum contraction) for the action   $L_A$. If $A$ is conjugate to an orthogonal matrix, then we say that $A$ is \emph{elliptic}. In this case, the action  $L_A$ preserves ellipses in $\mR^n$, and all orbits of $L_A$ and its inverse are bounded. Finally, if all eigenvalues of $A$ have norm $1$, but $A$ is not elliptic, then we say that $A$ is \emph{parabolic}. In this case, $A$ is conjugate (over $\mC$) to an upper triangular matrix with all diagonal entries of norm $1$, and so the norm $\|A^{\ell}\|$ grows as a polynomial function of the power $\ell$. The dynamics of $L_A$ in this case is distal, which is also  a dynamically defined property. 
 
 One key idea of Pesin theory is that the hyperbolicity property is well defined also for  non-autonomous linear approximations to smooth dynamical systems, so we look for this behavior on the level of derivative cocycles. This is the provided by the following concept.

\begin{defn}\label{def-hyperflow}
$w \in M$ is a \emph{hyperbolic point} of the flow if 
$$e_{\varphi}(w) \equiv \limsup_{T \to \infty} ~  \left\{ \frac{1}{T} \cdot \max_{s \leq T} ~ \left\{ \ln \left\{ \| (D\varphi_s  \colon Q_w \to Q_{z})\|^{\pm} \right\}  \right\}\right\} ~ > ~  0$$
\end{defn}

\begin{lemma}
The set of hyperbolic points $\cH(\varphi) = \{ w \in M \mid  e_{\varphi}(w) > 0 \}$ is flow-invariant.
\end{lemma}

\medskip

One of the first basic results if that if the set of hyperbolic points is non-empty, then the flow itself has   hyperbolic behavior on special subsets where the ``lim sup'' is replaced by a limit:

\begin{prop} \label{prop-hypmeas} 
Let $\varphi$ be a $C^1$-flow. Then 
the closure  $\overline{\cH(\varphi)} \subset M$ supports an invariant ergodic probability measure $\mu_*$ for $\varphi$, for which there exists $\lambda > 0$ such that  for  $\mu_*$-a.e. $w$,
$$e_{\varphi}(w) = \lim_{s \to \infty}   \left\{ \frac{1}{s} \cdot \ln \{ \| D\varphi_s   \colon Q_w \to Q_{z}\| \right\} = \lambda  $$ 
\end{prop}
 \proof  This follows from the continuity of the derivative and its cocycle property, the definition of the asymptotic Schwartzman cycle associated to a flow \cite{Schwartzman1957}, plus the usual subadditive techniques  of Oseledets Theory \cite{Oseledets1968, Pesin1977, BarreiraPesin2002}. 
 \endproof

We want to apply the ideas behind Proposition~\ref{prop-hypmeas}   to the derivatives of the maps in the holonomy pseudogroup $\cGF$.  The difficulty is that the orbits of the pseudogroup are not necessarily  ordered into a single direction along which the leaf hyperbolicity is to be found, and hence along which the integrals are defined in obtaining  the Schwartzman asymptotic cycle as in the above. One  approach is  to associate a   flow to a foliation $\F$, such that this flow      captures the dynamical information for $\F$. Such a flow exists,  and was introduced in the papers \cite{Hurder1988, Walczak1988}.

Let $w \in M$ and consider $L_w$ as a complete Riemannian manifold. 
For $\vec{v} \in T_w\F = T_wL_w$ with $\|\vec{v}\|_w = 1$, there is unique geodesic $\tau_{w, \vec{v}}(t)$ starting at $w$ with $\tau_{w,\vec{v}}'(0) = \vec{v}$.

 Define the map 
 $\varphi_{w, \vec{v}} \colon \mR \to M$ by $\varphi_{w,\vec{v}}(w) = \tau_{w, \vec{v}}(t)$. 
Let  $\whM = T^1\F$ denote the unit tangent bundle to the leaves, then the maps $\varphi_{w, \vec{v}} $  define the \emph{foliation geodesic flow}
$$\varphi^{\F}_{t} \colon \mR \times \whM \to \whM$$
Let $\whF$ denote the foliation on $\whM$ whose leaves are the unit tangent bundles to leaves of $\F$. Then the following is immediate from the definitions:

 \begin{lemma} \label{lem-normal}
 $\varphi^{\F}_{t}$ preserves the leaves of the foliation $\whF$ on $\whM$,  and hence $D\varphi^{\F}_{t}$  preserves the normal bundle $\whQ  \to \whM$ for $\whF$. 
 \end{lemma}

 Lemma~\ref{lem-normal}  makes is possible to give  an  extension of  Definition~\ref{def-hyperflow},  to the case of the normal derivative cocycle for  the foliation geodesic flow.   Consider the following  three possible cases for the asymptotic behavior of this cocycle.

\begin{defn}\label{def-hyperfol} 
Let $\varphi^{\F}_{t}$ be the foliation geodesic flow for a $C^1$-foliation $\F$. Then  $\whw \in \whM$  is:
\begin{description}
\item[H]    \emph{hyperbolic} if 
$$e_{\F}(\whw) \equiv \limsup_{T \to \infty} ~   \left\{ \frac{1}{T} \cdot \max_{s \leq T} ~       \left\{ \ln \| (D\varphi^{\F}_{s}   \colon \whQ_{\whw} \to \whQ_{\whz})\|^{\pm} \right\}  \right\} ~ > ~  0$$
\item[E]   \emph{elliptic}  if $e_{\F}(\whw) =0$, and there exists $\kappa(\whw)$ such that 
$$   \| (D\varphi^{\F}_{t}   \colon \whQ_{\whw} \to Q_{\whz}) \|^{\pm}  \leq \kappa(\whw) ~ \text {for all} ~  t \in \mR  $$
\item[P] \emph{parabolic}  if $e_{\F}(\whw) =0$, and $\whw$ is not elliptic.
\end{description}
\end{defn}

There is  a   variation on this definition which is also very    useful, which   takes into account the fact that for foliation dynamics, one does not necessarily have a preferred direction for the foliation geodesic flow: one considers all possible directions simultaneously in  Definition~\ref{def-hyperflow}.

Let $\|\gamma\|$ denote the minimum length of a geodesic $\sigma$ whose holonomy $h_{\sigma_{w,z}}$ defines the germ $\gamma = [h_{\sigma_{w,z}}]_w \in \GF$. Let  $D_w\gamma = D_w h_{\sigma_{w,z}}$ denote the derivative at $w$.

\begin{defn}  \label{def-expansionfunction}
The   transverse expansion rate  function for $\cGF$ at $w$ is
 \begin{equation}\label{ep-atgn}
 e(\cGF , T, w) =    \max_{\|\gamma\|  \leq T} ~ \left\{   \ln   \left\{ \| D_w \gamma  \|^{\pm}  \right\}  \right\}
\end{equation}
  \end{defn}
 Note that $ e(\cGF , d, w)$ is a Borel function of $w \in \cT$, as each norm function $\| D_{w'} h_{\sigma_{w,z}} \|$ is continuous for $w' \in D(h_{\sigma_{w,z}})$ and the maximum of Borel functions is Borel.
  
\begin{defn}  
The asymptotic transverse expansion rate at $w \in \cT$ is
 \begin{equation}\label{eq-atg}
e_{\F}(w)  = e(\cGF, w) = \limsup_{T \to \infty} ~  \left\{ \frac{1}{T} \cdot e(\cGF, T, w) \right\} ~ \geq ~ 0
  \end{equation}
\end{defn}  
The limit of Borel functions is Borel, and each $ e(\cGF, d, w)$ is a Borel function of $w$, hence $e(\cGF, w)$ is Borel. 
 The value $e_{\F}(w)$ can be thought of as the ``maximal Lyapunov exponent'' for the holonomy groupoid   at $w$.
  Analogous to the flow case, the chain  rule and the definition of  $e_{\F}(w)$ imply:
\begin{lemma}\label{lem-orbitinv}
  $e_{\F}(z) = e_{\F}(w)$ for all $z \in L_w \cap \cT$. Moreover, the value of $e_{\F}(w)$ is independent of the choice of Riemannian metric on $TM$. 
  Hence,  the expansion function $e(w)$ is constant along leaves of   $\F$, and is a dynamical invariant of $\F$. 
\end{lemma}

There is a  trichotomy for   the expansion rate function $ e(\cGF , d, w) $ analogous to that in Definition~\ref{def-expansionfunction}. Thus, there is   a decomposition of the manifold $M$ into those leaves which satisfy one of the three types of asymptotic behavior for the normal derivative cocycle: 

\begin{thm} [Dynamical decomposition of foliations] \label{thm-decomp} 
  Let $\F$ be a $C^1$-foliation on a compact manifold $M$. Then $M$ has a decomposition into disjoint saturated Borel subsets, 
  \begin{equation}\label{eq-trichotomy}
M = \EF \cup \PF \cup \HF
\end{equation}
which are the leaf saturations of the sets defined by:
\begin{enumerate}
\item Elliptic:  $\ET   = \{ w \in \cT \mid \forall ~ T \geq 0, ~  e(\cGF , T, w)  \leq \kappa(w) \}$
\item Parabolic:  $\PT    = \{ w \in \cT \setminus  \EF  \mid  e(\cGF, w) = 0 \}$
\item Hyperbolic:  $\HT   = \{ w \in \cT \mid  e(\cGF, w) > 0 \}$
\end{enumerate}
\end{thm}
Note that $w \in \ET$ means that  the holonomy homomorphism $D_w  \gamma$   has bounded image in  $GL(q, \mR)$, contained in a ball of radius  
 $\exp\{ \kappa(w) \} = \sup \{ \|D_w \gamma \| \mid  \gamma \in \cGF^w \}$, where $\cGF^w$ denotes the germs of holonomy transport along paths starting at $w$.

The nomenclature in Theorem~\ref{thm-decomp} reflects the trichotomy for the dynamics of a matrix $A \in GL(q,\mR)$ acting via the associated  linear transformation $L_A \colon \mR^q \to \mR^q$: The elliptic points are the regions where the infinitesimal holonomy transport ``preserves ellipses up to bounded distortion''.  
The parabolic points are where the    infinitesimal holonomy acts similarly to that of a parabolic subgroup of $GL(q,\mR)$; for example, the action is ``infinitesimally distal''. The   hyperbolic points are where the  the infinitesimal holonomy has some degree of exponential expansion. Perhaps more properly, the set $\HF$   should be called ``non-uniform, partially hyperbolic leaves''.  The study of the dynamical properties of the set of hyperbolic leaves $\HF$  has close analogs with the study of non-uniformly hyperbolic dynamics for flows, as in  \cite{BDV2005}.

\medskip
 
 The decomposition in Theorem~\ref{thm-decomp}  has many applications to the study of foliation dynamics and classification results, as  discussed  for example in   \cite{Hurder2009}, and also  \cite{HurderKatok1987,Hurder1988,Hurder2000b}.    
We illustrate some of these applications with examples and selected results. Here is one important concept:
 
\begin{defn}\label{def-thm}
 An invariant probability measure $\mu_*$ for the foliation geodesic flow on $\whM$ is said to be \emph{transversally hyperbolic} if 
$e_{\F}(\whw)   > 0$ for $\mu_*$-a.e. $\whw$. 
\end{defn}
The function $e_{\F}(\whw) $ is constant on orbits, so is constant on the ergodic components of  $\mu_*$. 
Thus,  if $\mu_*$ is an ergodic invariant measure for the foliation geodesic flow, then $\mu_*$   transversally hyperbolic means there is some constant $\lambda(\mu_*) > 0$ with $\lambda(\mu_*) = e_{\F}(\whw)$ for for $\mu_*$-a.e. $\whw$.

Also,   note that the support of a transversally hyperbolic  measure $\mu_*$ is contained in the unit tangent bundle $\whM$, and not $M$ itself. A generic point $\whw$ in the support of $\mu_*$ specifies both a point in a   leaf,  and the direction along which to follow a geodesic to find   infinitesimal normal hyperbolic behavior.

\begin{thm}\label{thm-thm}
 Let $\F$ be a $C^1$-foliation of a compact manifold.  If $\HF \ne \emptyset$, then the foliation geodesic flow has at least one transversally hyperbolic ergodic measure, which is contained in the closure of unit tangent bundle over $\HF$. 
\end{thm}
\proof  The proof is technical, but   basically follows from  calculus techniques applied to the foliation pseudogroup, as in Oseledets Theory. The key point is that if $L_w \subset \HF$, then there is a sequence of geodesic segments of lengths going to infinity on the leaf $L_w$, along which the transverse infinitesimal expansion grows at an exponential rate. Hence, by continuity of the normal derivative cocycle and the cocycle law, 
these geodesic segments converge to a transversally hyperbolic  invariant probability measure $\mu$  for the foliation geodesic flow. The existence of an ergodic component $\mu_*$ for this measure with positive exponent then follows from the properties of the ergodic decomposition of $\mu$. 
\endproof

\begin{cor}\label{cor-thm}
Let $\F$ be a $C^1$-foliation of a compact manifold with $\HF \ne \emptyset$. Then there exist $w \in \overline{\HF}$ and a unit vector $\vec{v} \in T_w\F$ such that the forward orbit of the geodesic flow through $(w, \vec{v})$ has a transverse direction which is \emph{uniformly exponentially contracting}. 
\end{cor}

Let us return to  the examples introduced earlier, and consider what the trichotomy decomposition means in each case.
 For the linear foliation of the $2$-torus in Figure~1, every point is elliptic, as the foliation  is Riemannian. 
However, if $\F$ is a $C^1$-foliation which is topologically semi-conjugate to a linear foliation, so is a generalized Denjoy example, then $M_{\cP}$ is not empty!  Shigenori Matsumoto has given  a new  construction of  Denjoy-type $C^1$-foliations on the $2$-torus  for which the exceptional minimal set consists of elliptic points, and the points in the wandering set are all parabolic \cite{Matsumoto2010}.
  
Consider next the   Reeb foliation of the solid torus, as in Figure~3.  
Pick $w \in M$ on an interior parabolic leaf, and a direction $\vec{v} \in T_wL_w$.   Follow the geodesic $\sigma_{w,\vec{v}}(t)$ starting from $w$. 
It is asymptotic to the boundary torus, so defines a limiting Schwartzman cycle on the boundary torus for some   flow. Thus, it limits on either a circle, or a lamination. This will be a hyperbolic measure if the holonomy of the compact leaf is hyperbolic. Note that the exponent of the invariant measure for the foliation geodesic  flow   depends on the direction of the geodesic used to define it. 

 One of the basic problems about the foliation geodesic flow is to understand   the support of its transversally hyperbolic invariant measures whose generic starting points lie in $\HF$, and if the leaves intersecting the supports of these measures have ``chaotic'' behavior.

\section{Counting} \label{sec-counting}

The decomposition of the foliated manifold $\ds M = \EF \cup \PF \cup \HF$
uses  the asymptotic properties of the normal  ``derivative cocycle''  $D \colon \cGF \to GL(n, \mR)$, where the   transverse expansion is allowed to  ``develop in any direction'' when the leaves are higher dimensional. 

 A basic question is then, how do you tell whether one of the Borel, $\F$-saturated components, such as the hyperbolic set $\HF$, is non-empty? Moreover, it  is natural to speculate whether  the ``geometry of the leaves'' influences the structure of the sets in the trichotomy (\ref{eq-trichotomy}). To this end, we consider in this section  the notion of the  \emph{growth rates of leaves}. This leads to a variety of ``counting type'' invariants for foliation dynamics, and various insights into  the behavior of the derivative cocycle.

\begin{figure}[H]
\setlength{\unitlength}{1cm}
\begin{picture}(30,4)
\put(5.6,.2){\includegraphics[width=46mm]{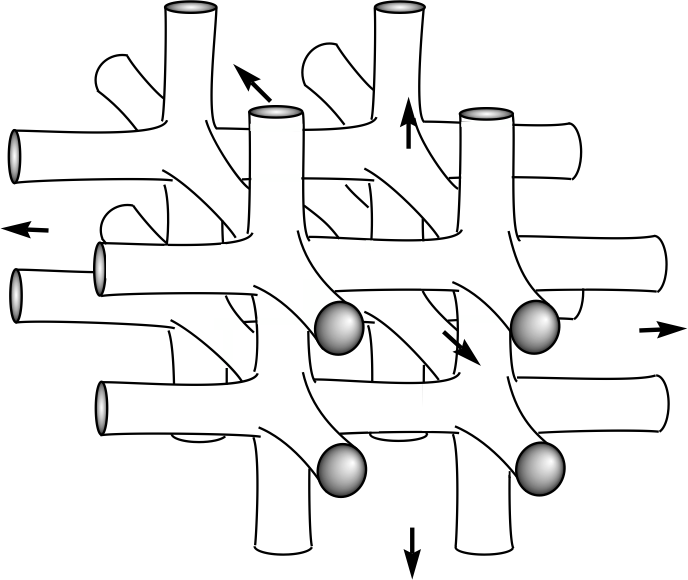}}
\end{picture}
\caption{The infinite ``Jungle Gym''}
\label{Figure7} 
\vspace{-20pt}
\end{figure}

 Let us first  consider some        examples with more complicated leaf geometry than seen above. 
   Figure~\ref{Figure7} depicts what  is called the    ``Infinite Jungle Gym'' in the foliation literature \cite{PS1981, CandelConlon2000}.

The surface in Figure~\ref{Figure7} can be realized as a leaf of a circle bundle over a compact surface, where the holonomy consists of three commuting linearly independent  rotations of the circle. Thus, even though this is a surface of infinite genus, the transverse holonomy is just a generalization of that for the Denjoy example, in that it consists of a group of isometries with dense orbits for the circle $\mS^1$.

\begin{figure}[H]
\setlength{\unitlength}{1cm}
\begin{picture}(30,3)
\put(4.6,.1){\includegraphics[width=52mm]{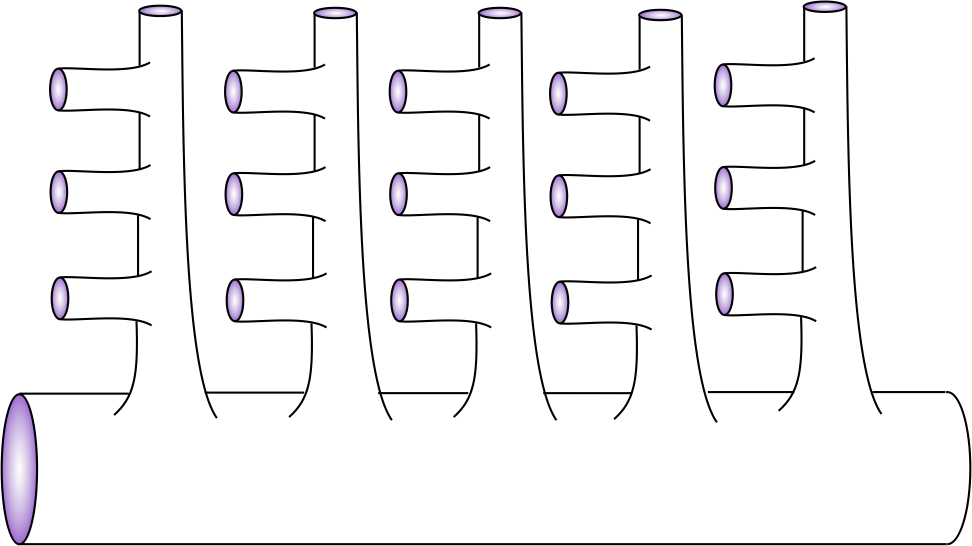}}
\end{picture}
\caption{A leaf of ``Level 2''}
\label{Figure8} 
\vspace{-20pt}
\end{figure}

 The next manifold $L$ in Figure~\ref{Figure8} doesn't have a cute name, but has the interesting property that its space of ends $\cE(L_1)$  has non-empty derived set, but the second derived set is empty.

 This manifold  can be realized as a leaf in a smooth foliation which is asymptotic to a compact surface of genus two.  
The construction of the foliation in which this occurs is given in \cite{CandelConlon2000}. It is just one example of a large class of foliations with a proper leaf of finite depth \cite{CantwellConlon1978, CantwellConlon1981b, Hector1983, Tsuchiya1979a, Tsuchiya1980b}. As with the Reeb foliation, the hyperbolic invariant measures for the flow are concentrated on the limiting compact leaf.  The dynamics is not chaotic.

\begin{figure}[H]
\setlength{\unitlength}{1cm}
\begin{picture}(20,3.4)
\put(5,.2){\includegraphics[width=46mm]{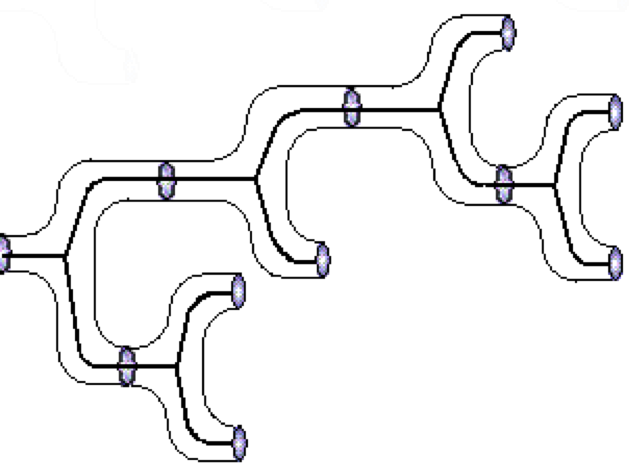}}
\end{picture}
\caption{A leaf with Cantor endset}
\label{Figure9} 
\vspace{-20pt}
\end{figure}

 The manifold in Figure~\ref{Figure9}  has  endset $\cE(L_2)$ which is a Cantor set,    equal to its own derived set. 
 
  We present in more detail the construction of a foliated manifold containing  this   as a leaf,  called  the ``Hirsch foliation'', introduced in  \cite{Hirsch1975}, as it illustrates a basic theme of the lectures and the elementary construction yields sophisticated dynamics.  See \cite{BHS2006} for generalizations of this construction. 
    
 \underline{Step 1}: Choose an analytic embedding of  $\mS^1$ in the  solid torus $\mD^2 \times \mS^1$ so that its image is twice a generator 
of  the fundamental group of the solid torus. Remove an open tubular 
neighborhood of the embedded $\mS^1$,  resulting in the manifold with boundary in Figure~\ref{Figure10}:

\begin{figure}[H]
\setlength{\unitlength}{1cm}
\begin{picture}(20,4)
\put(5,.2){\includegraphics[width=50mm]{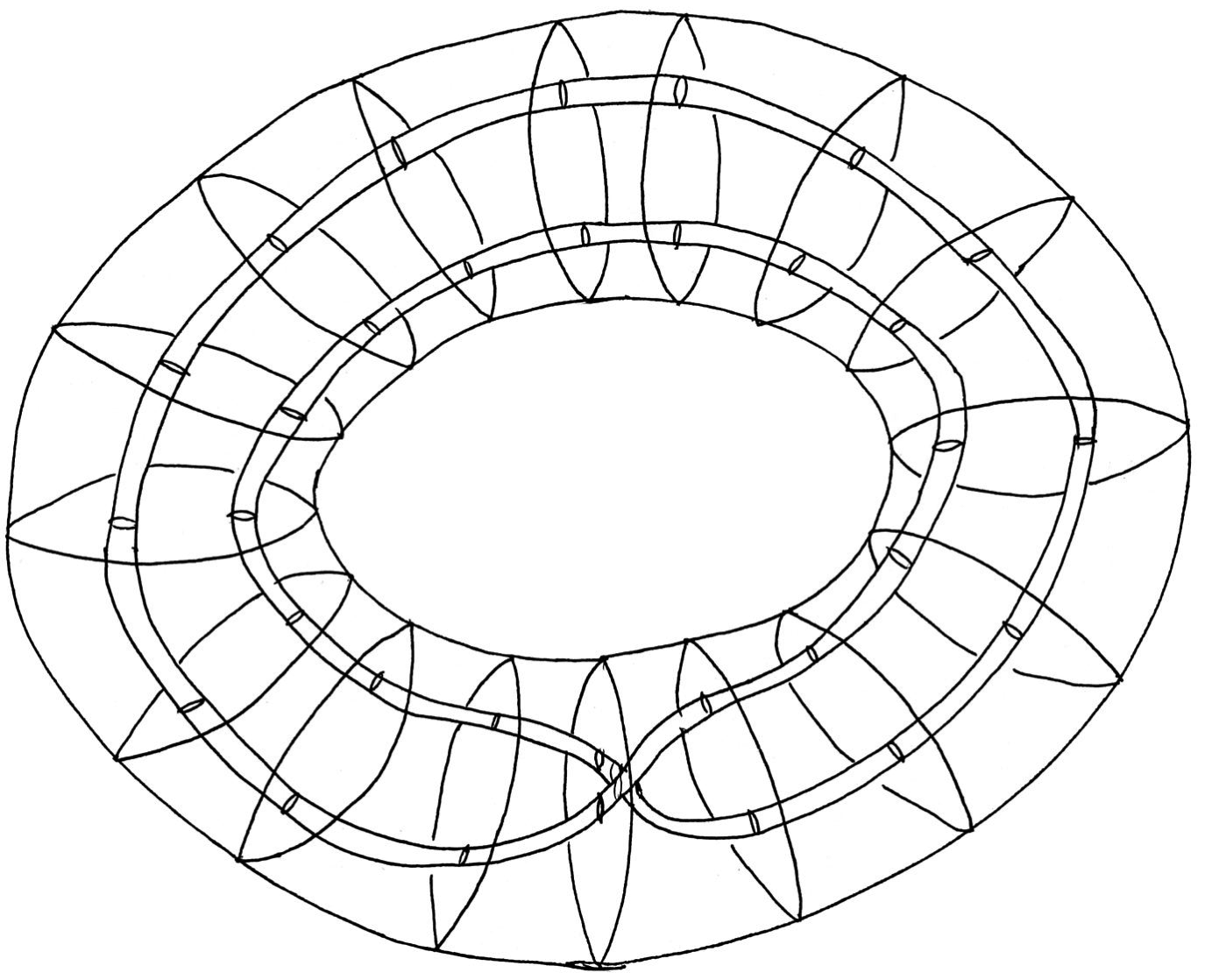}}
\end{picture}
\caption{Solid torus with tube drilled in it}
\label{Figure10} 
\vspace{-20pt}
\end{figure}

\underline{Step 2}:   What remains is a three  dimensional  manifold $N_1$ whose boundary is two disjoint copies of $\mT^2$.  
$\mD^2 \times \mS^1$ fibers over $\mS^1$ with fibers the 2-disc.  
This  fibration --  restricted to $N_1$ -- foliates $N_1$ with leaves consisting of 2-disks with two open subdisks removed.

\begin{figure}[H] 
\setlength{\unitlength}{1cm}
\begin{picture}(20,2)
\put(6.4,.2){\includegraphics[width=25mm]{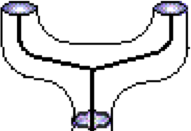}}
\end{picture}
\caption{A ``pair of pants''}
\label{Figure11} 
\vspace{-20pt}
\end{figure}

  Identify the two components of the boundary of $N_1$ by a 
diffeomorphism which covers the map $h(z) = z^2$ of $S^1$  to 
obtain the  manifold $N$.  Endow $N$ with a Riemannian metric; then the punctured 2-disks foliating $N_1$ can now be viewed as    pairs of pants, as in  Figure~\ref{Figure11}.

 \vfill
 
\underline{Step 3}: The foliation of $N_1$  is transverse to the boundary, so  
the punctured 2-disks assemble to yield  a foliation of   foliation $\F$ on $N$, where the leaves  without holonomy   (corresponding to irrational points for the chosen doubling map of $S^1$) are infinitely branching surfaces, decomposable into pairs-of-pants which correspond to the punctured disks in $N_1$.

A basic  point is that this works for any covering map $f \colon \mT^2 \to \mT^2$ homotopic to the doubling map $h(z)$ along a meridian. 
In particular, as Hirsch remarked in his paper, the proper choice of such a ``bonding map''  results in a codimension-one, real analytic  foliation, such that all leaves accumulate on a unique exceptional minimal set.

\begin{figure}[H]
\setlength{\unitlength}{1cm}
\begin{picture}(20,3.6)
\put(6.4,.2){\includegraphics[width=36mm]{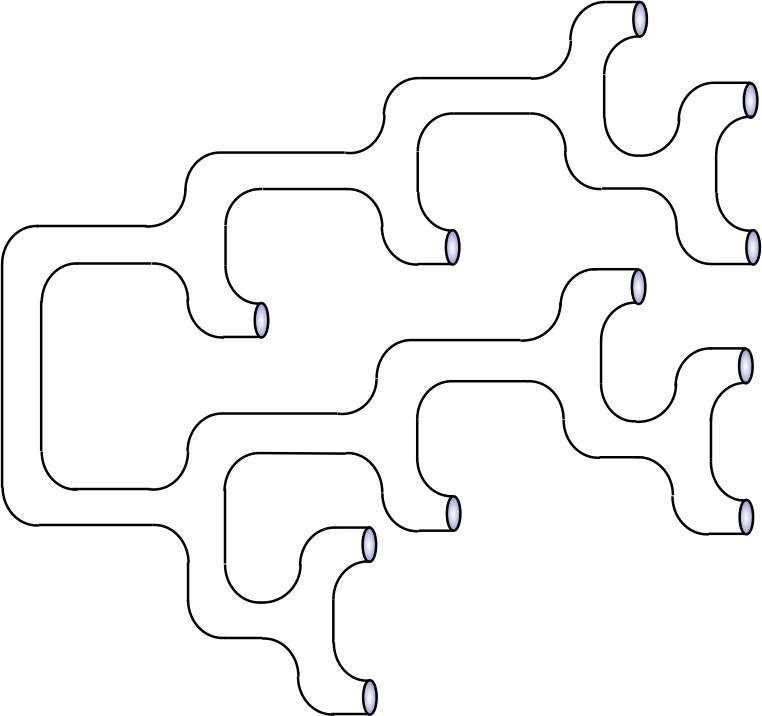}}
\end{picture}
\caption{Leaf for eventually periodic orbit}
\label{Figure12} 
\vspace{-20pt}
\end{figure}

  The Hirsch foliation always has a leaf $L_w$ pictured as in Figure~\ref{Figure12}, corresponding to a forward periodic orbit of the doubling map $g \colon \mS^1 \to \mS^1$. 
Consider the behavior of the geodesic flow, starting at the ``bottom point'' $w \in L_w$.
For a each radius $R \gg 0$, the terminal points of the geodesic rays of length at most $R$ will ``jump'' between the $\mu^R$ ends of this compact subset of the leaf, for some $\mu > 1$. Thus, for these examples, a small variation of the initial vector $\vec{v}$ will result in a large variation of the terminal end of the geodesic $\sigma_{w, \vec{v}}$.
  
 \medskip

 The constant $\mu$   appearing in the above example seems to be an ``interesting'' property of the foliation dynamics, and a key point is that it can be obtained by ``counting'' the complexity of the leaf at infinity, following a scheme introduced by Joseph Plante for leaves of  foliations \cite{Plante1975}. 
 
 Recall the holonomy pseudogroup $\cGF$ constructed in section~\ref{sec-basics}, modeled on a complete transversal $\cT = \cT_1 ~ \cup \cdots \cup ~ \cT_k$ associated to a finite covering of $M$ by foliations charts. Given $w \in \cT$ and $z \in L_w \cap \cT$ and a leafwise path $\tau_{w,z}$ joining them, we obtain an element $h_{\tau_{w,z}} \in \cGF$.
 
 The orbit of $w \in \cT$ under $\cGF$ is  
  $$\cO(w) = L_w \cap \cT = \{ z \in \cT \mid  g(w) = z, ~ g \in \cGF,  w \in Dom(g)\} $$

 Introduce the \emph{word norm} on elements of $\cGF$. 
Given  open sets $U_i \cap U_j \ne \emptyset$ in the fixed covering of $M$ by foliation charts,   they define an element $h_{i,j} \in \cGF$. By the definition of holonomy along a path, for each $\tau_{w,z} \colon [0,1] \to L_w$ there is a sequence of indices 
  $\{i_0 , i_1 , \ldots , i_{\ell}\}$ such that 
  $$[h_{\tau_{w,z}}]_w = [h_{i_{\ell -1} , i_{\ell}} \circ \cdots \circ h_{i_1 , i_0}]_w$$
  That is, the germ of the holonomy map $h_{\tau_{w,z}}$ at $w$ can be expressed as the composition of $\ell$ germs of the basic maps $h_{i,j}$.
  We then say that $\gamma = [h_{\tau_{w,z}}]_w$ has word length at most $\ell$. Let $\|\gamma\|$ denote the least such $\ell$ for which this is possible.
  The norm of the identity  germ is defined to be $0$. 
 
 Define the ``orbit of $w$ of radius $\ell$ in the groupoid  word norm'' to be:
   $$\cO_{\ell}(w) =   \{ z \in \cT \mid  g(w) = z, ~ g \in \cGF,  w \in Dom(g), \| [g]_w\| \leq \ell\} $$
   
\begin{defn}\label{def-growthfcn}
 The growth function of an orbit is defined as  $Gr(w, \ell) = \#  \cO_{\ell}(w)$.
\end{defn}
   
 Of course, the growth function for $w$ depends upon many choices.   However, its ``growth type function'' is   independent of   choices, as observed by Plante. This follows from  one of the basic facts of the theory, that the word norm on $\cGF$ is bounded above by a multiple of the length of geodesic paths.
  
  \begin{prop} \cite{Milnor1968, Plante1975}\label{prop-word}
  Let $\F$ be a $C^1$-foliation of a compact manifold $M$. Then there exists a constant $C_m > 0$ such that for all $w \in \cT$ and $z \in L_w \cap \cT$, if $\sigma_{w,z} \colon [0,1] \to L_w$ is a leafwise geodesic segment from $w$ to $z$ of length $\| \sigma_{w,z}\|$, then
  $$ \| [h_{\sigma_{w,z}}]_w \|  ~ \leq ~ C_m \cdot \| \sigma_{w,z}\|$$
  \end{prop}
  
 In order to obtain a well defined invariant of growth of an orbit, one introduces the notion the \emph{growth type} of a function. The one which  we use (there are many - see \cite{Hector1977, Egashira1993,Badura2005}) is essentially the weakest one.  Given given   functions $f_1, f_2 \colon [0,\infty) \to [0, \infty)$ say that $f_1 \lesssim f_2$ if there exists constants $A, B, C > 0$ such that for all $r \geq 0$, we have that $  f_2(r) ~ \leq ~ A \cdot f_1(B \cdot r) + C$.
 Say that    $f_1 \sim f_2$ if both $f_1 \lesssim f_2$ and $f_2 \lesssim f_1$ hold.    This defines   equivalence relation on functions, which defines their   \emph{growth class}.

  One can consider a variety of special classes of growth types. For example, note that if $f_1$ is the constant function and $f_2 \sim f_1$ then $f_2$ is constant also. 
  
We say that $f$ has \emph{exponential growth type} if $f(r) \sim \exp(r)$. Note that $\exp(\lambda \cdot r) \sim \exp(r)$ for any $\lambda > 0$, so there is only one growth class of ``exponential type''.
  
  A function has \emph{nonexponential growth type} if $f  \lesssim \exp(r)$,  but $\exp(r) \not\lesssim f$.

We also have a subclass of uniform nonexponential growth type, called in the author's papers by  \emph{subexponential growth type},  if for any $\lambda > 0$ there exists $A, C > 0$ so that $f(r) \leq A \cdot \exp(\lambda \cdot r) + C$.

  Finally, $f$ has  \emph{polynomial growth type} if there exists $d \geq 0$ such that $f  \lesssim r^d$. The growth type is exactly polynomial of degree $d$ if $f \sim r^d$. 
  
\begin{defn}\label{def-orbitgrowthtype}
The growth type of an orbit $\cO(w)$  is the growth type of  $Gr(w, \ell) = \#  \cO_{\ell}(w)$. 
\end{defn}  
A basic result of Plante  is that:
 \begin{prop}\label{prop-orbitgrowthtype} Let $M$ be a compact manifold.  Then for all $w \in \cT$, $Gr(z, \ell) \lesssim \exp(\ell)$. Moreover, 
 for  $z \in L_w \cap \cT$, then $Gr(z, \ell) \sim  Gr(w, \ell)$.
 Thus, the growth type of a leaf $L_w$ is well defined, and we say that  $L_w$ has the growth type of the function $Gr(w, \ell)$.
  \end{prop}
  We can thus speak of a leaf $L_w$ of $\F$ having exponential growth type, and so forth. For example, the Infinite Jungle Gym  manifold  in Figure~\ref{Figure7}  has growth type exactly polynomial of degree $3$, while the leaves of the Hirsch foliations (in Figures~\ref{Figure10} and \ref{Figure12}) have exponential growth type.  
 
 Before continuing with the discussion of the growth types of leaves, we note the correspondence between these ideas and a basic problem in geometric group theory. 
 Growth functions for finitely generated groups are a basic object of study in geometric group theory. 
  
  Let $\Gamma = \langle \gamma_0 = 1, \gamma_1, \ldots , \gamma_k\rangle$ be a finitely generated group. Then $\gamma \in \Gamma$ has word norm 
 $\|\gamma \| \leq \ell$ if we can express $\gamma$ as a product of at most $\ell$ generators, 
 $\gamma = \gamma_{i_1}^{\pm} \cdots \gamma_{i_d}^{\pm}$. Define the ball of radius $\ell$ about the identity of $\Gamma$ by 
   $$\Gamma_{\ell}  \equiv \{\gamma \in \Gamma \mid \|\gamma\| \leq \ell\}$$
The growth function  $Gr(\Gamma, \ell) = \# \Gamma_{\ell}$ depends upon the choice of generating set for $\Gamma$, but its growth type does not.
The following is a  celebrated  theorem  of Gromov:
\begin{thm} \cite{Gromov1981}
Suppose $\Gamma$ has polynomial growth type for some generating set. Then there exists  a subgroup of finite index $\Gamma' \subset \Gamma$ such that $\Gamma'$ is a nilpotent group.
  \end{thm}
In general, one asks to what extent does the growth type of a group determine its algebraic structure?  

Questions of a similar nature can be asked about leaves of foliations;  especially,  to what extent does the growth function of leaves determine how they are embedded in a compact manifold, and  the dynamical properties of the foliation? 
 
 Note that  there is a fundamental difference between the growth types of groups   and for leaves.   
  The homogeneity of groups implies that the growth rate is uniformly the same for balls in the word metric about any point $\gamma_0 \in \Gamma$.  That is, one can choose the constants $A, B, C > 0$ in the definition of growth type which are independent of the center $\gamma_0$.    
       For  foliation pseudogroups, there is a basic question about the uniformity of the growth function as the basepoint within an orbit varies:  
        
 \begin{quest} 
  How does the    function $d \mapsto Gr(w,d)$ behave, as a Borel function of $w \in \cT$?
\end{quest}

 Examples of Ana Rechtman \cite{Rechtman2009, AR2010} (see also \cite{Kaimanovich2001}) show that even for smooth foliations of compact manifolds, this function is not uniform as function of $w \in \cT$. If one requires uniformity of the growth function $\ell \mapsto Gr(w, \ell)$, as a function of $w \in \cT$, then one can ask if there is some form of analog  for foliation pseudogroups of the classification program for finitely generated groups.

\section{Exponential Complexity} \label{sec-complexity}

Section~\ref{sec-derivatives}   introduced exponential growth criteria for the normal derivative cocycle of the pseudogroup $\cGF$ acting on the transverse space $\cT$, and  section~\ref{sec-counting} discussed the growth types  for the orbits of the groupoid. In both cases, exponential behavior represents a type of exponential complexity for the dynamics of $\cGF$. These  examples are part of a larger theme, that when studying  classification problems, 
\emph{Exponential Complexity is Simplicity}.
In this section, we develop this theme further. First, we give an aside, presenting a well-known phenomenon for map germs.
 
Recall a simple example from advanced calculus. Let $f(x) = x/2$, and  let $g \colon (-\epsilon, \epsilon) \to \mR$ be a smooth map with $g(0) = 0$,   $g'(0) = 1/2$.
 Then $g \sim f$ near $x=0$.  That is, for $\delta >0$ sufficiently small, there is a smooth map $h \colon (-\epsilon, \epsilon) \to \mR$ such that 
 $h^{-1} \circ g \circ h = f (x)$ for all $|x| < \delta$.
 
  This illustrates the principle that   exponentially contracting maps,  or more generally hyperbolic maps in higher dimensions, the derivative  is a complete  invariant for  their germinal conjugacy class at the fixed point. 
 For maps which are ``completely flat'' at the origin,  where $g(0) = 0$, $g'(0) = 1$, $g^k(0) = 0$ for all $k > 1$,   their ``classification'' is much more difficult \cite{Mather1968,Wall1971}. So, in contrast we have  \emph{Subexponential Complexity is Most Nettlesome}. On the other, there are invariants for foliations which are only defined for amenable systems, as discussed later, so the real point is that this distinction between exponential and subexponential complexity pervades   classification problems.

Analogously, for foliation dynamics, and the related   problem of studying the dynamics    of a finitely generated group acting smoothly on a compact manifold, exponential complexity in the dynamics   often gives rise to hyperbolic behavior for the holonomy pseudogroup. Hyperbolic maps can be put into a standard form, and so one obtains a fundamental  tool for studying the dynamics of the pseudogroup.  The problem   is thus, how does one pass from exponential complexity to hyperbolicity?

One issue with the   ``counting argument'' for the growth of leaves  is that just  counting the growth rate  of an orbit  ignores fundamental information about expansivity of the dynamics. 
The orbit growth rate counts the number of times the leaf crosses a transversal $\cT$ in a fixed distance within the leaf, but does not take into account whether these crossings are ``nearby'' or ``far apart''. 
For example,  there are Riemannian foliations with all leaves of exponential growth type. See   \cite{Richardson2001}, for example. Thus, exponential orbit growth rate need not imply   transversally hyperbolic behavior. 

 On the other hand,  in  the Hirsch examples, the handles at the end of each  ball of radius $\ell$ in a leaf  appear to be widely separated transversally, so somehow this is different.  The holonomy pseudogroup $\cGF$ of the Hirsch example is topologically semi-conjugate  to the pseudogroup  generated by the doubling map $z \mapsto z^2$ on $\mS^1$.
  After $\ell$-iterations, the inverse map to $h(z) = z^{2^{\ell}}$ has derivative of norm $2^{\ell}$, and so  for a Hirsch foliation modeled on this map, every leaf is transversally hyperbolic.

The  \emph{geometric entropy}   for pseudogroup $C^1$-actions, introduced by Ghys, Langevin, and  Walczak \cite{GLW1988},  gives a measure    of their exponential complexity which   combines    the two types of complexity. It has   found many applications in the study of foliation dynamical systems. One example of this is the surprising role of these invariants in showing that  certain secondary classes of $C^2$-foliations are zero in cohomology if the entropy vanishes  \cite{CantwellConlon1984, Hurder1986,HurderKatok1987, HLa2000}.

 We begin with the basic notion of $\epsilon$-separated sets, due to Bowen \cite{Bowen1971} for diffeomorphisms, and extended to groupoids in \cite{GLW1988}.
  Let $\epsilon > 0$ and $\ell > 0$. 
   A subset $\cE \subset \cT$ is said to be $(\epsilon, \ell)$-separated if   for all $w,w' \in \cE \cap \cT_i$   there exists $g \in \cGF$ with $w,w' \in Dom(g) \subset \cT_i$, and  $\|g\|_w \leq \ell $  so that  $d_{\cT}(g(w), g(w')) \geq \epsilon$. 
  If $w  \in \cT_i$ and $w' \in \cT_j$ for $i \ne j$ then   they are   $(\epsilon, \ell)$-separated by default. 
 
 \medskip

The  ``expansion growth function'' counts the maximum  of this quantity:
  $$h(\cGF, \epsilon, \ell) =     \max \{ \# \cE \mid \cE \subset \cT ~ \text{is} ~ (\epsilon,\ell) \text{-separated} \}   $$
  
   If the pseudogroup $\cGF$ consists of isometries, for example, then applying elements of $\cGF$ does not increase the separation between points,  so the growth functions $h(\cGF, \epsilon, \ell)$ have polynomial growth of degree equal the dimension of $\cT$, as functions of $\ell$. Thus, if the functions $h(\cGF, \epsilon, \ell)$ have greater than polynomial growth type, then the action of the pseudogroup cannot be elliptic, for example.  
   
Introduce the measure of the exponential growth type of the expansion growth function:
 \begin{eqnarray}
h(\cGF, \epsilon) ~ & = & ~ \limsup_{d \to  \infty} ~ \frac{ \ln \left\{ \max \{ \# \cE \mid \cE ~ \text{is} ~ (\epsilon,d) \text{-separated} \} \right\}}{d} \label{eq-entropy1} \\
h(\cGF) ~ & = & ~ \lim_{\epsilon \to 0} ~  h(\cGF, \epsilon) \label{eq-entropy2}
\end{eqnarray}
   Then we have the fundamental result of Ghys, Langevin, and  Walczak \cite{GLW1988}:
\begin{thm} Let $\F$ be a $C^1$-foliation of a compact manifold $M$. Then 
 $h(\cGF) $ is finite.
 
Moreover, the property $h(\cGF) > 0$ is independent of all choices.
  \end{thm}
  
For example,  if $\F$ is defined by a $C^1$-flow $\phi_t \colon M \to M$,  then  $h(\cGF) > 0$ if and only if  $h_{top}(\phi_1) > 0$.  Note that $h(\cGF)$ is defined using the word growth function for orbits, while the topological entropy of the map $\phi_1$ is defined using the geodesic length function (the time parameter) along leaves. These two notions of ``distance along orbits'' are comparable, which can be used to give estimates, but not necessarily any more precise relations. This point  is discussed in detail in \cite{GLW1988}.  

In any case, the essential information contained in the invariant $h(\cGF)$ is simply whether the foliation $\F$ exhibits exponential complexity for its orbit dynamics, or not. Exploiting further the information contained in this basic invariant of $C^1$-foliations has been one of the fundamental problems in the study of foliation dynamics since the introduction of the concept of geometric entropy    in 1988. 
 
 One aspect of the geometric entropy $h(\cGF)$ is that it is a ``global invariant'', which does not indicate ``where'' the chaotic dynamics is happening.  
  The author introduced a  variant of $h(\cGF)$ in \cite{Hurder2009}, the  \emph{local geometric entropy} $h(\cGF, w)$ of $\cGF$ which is a refinement of  the global entropy. The local geometric entropy  is analogous to the local measure-theoretic entropy for maps introduced by Brin and Katok    \cite{BrinKatok1983, GlasnerYe2009}. 
    
Given a subset $X \subset \cT$, in the definition  of $(\epsilon, \ell)$-separated sets above, we can demand that  the separated points   be contained in $X$,   yielding  the relative  expansion growth function:
  $$h(\cGF, X, \epsilon, \ell) =     \max \{ \# \cE \mid \cE \subset X ~ \text{is} ~ (\epsilon,\ell) \text{-separated} \}   $$
Form the corresponding limits as in (\ref{eq-entropy1}) and  (\ref{eq-entropy2}), to obtain 
 the \emph{relative geometric entropy} $h(\cGF, X)$. 
 
 Now, fix $w \in \cT$ and let    $X =  B(w,\delta) \subset \cT$ be the open $\delta$--ball about $w \in \cT$.
Perform the same double limit process as used to define $h(\cGF)$ for the sets $B(w,\delta)$, but then also let the radius of the balls tend to zero, to obtain:
\begin{defn}\label{def-localentropy}
The \emph{local geometric entropy} of $\cGF$ at $w$ is
\begin{equation}\label{eq-localentropy}
h_{loc}(\cGF, w) = \lim_{\delta \to 0} \Big\{  \lim_{\epsilon \to 0} \Big\{  \limsup_{\ell \to \infty} \frac{\ln \{h(\cGF, B(w,\delta),\ell,\epsilon)\}}{\ell}\Big\} \Big\}
\end{equation}
\end{defn}
The quantity $h_{loc}(\cGF, w)$   measures of the amount of ``expansion'' by the pseudogroup in an open neighborhood of $w$. 
The following estimate is elementary to show.
\begin{prop}[Hurder, \cite{Hurder2009}] \label{prop-locentprop}
Let $\cGF$   a $C^1$-pseudogroup. Then $h_{loc}(\cGF,w)$ is a Borel function of $w \in \cT$, and $h_{loc}(\cGF,w) = h_{loc}(\cGF,z)$ for $z \in L_w \cap \cT$. Moreover,
\begin{equation}\label{eq-locentest}
 h(\cGF)  ~  = ~  \sup_{w \in \cT} ~ h_{loc}(\cGF,w)
\end{equation}
\end{prop}
It follows that there is a disjoint Borel decomposition of $\cT$ into $\cGF$-saturated subsets
$\ds \cT = \bZ_{\F} \cap \bC_{\F}$, 
where $\bC_{\F} = \{w \in \cT \mid h(\cGF, w) > 0\}$ consists of the ``chaotic'' points for the groupoid action,  and $\bZ_{\F} = \{w \in \cT \mid h(\cGF, w) = 0\}$ are the ``tame'' points.  

\begin{cor}  $h(\cGF) > 0$   if and only if $\bC_{\F} \ne \emptyset$. 
\end{cor}

We   discuss in the next section  the    relationship     between local entropy   $h(\cGF, w) > 0$  and the transverse Lyapunov spectrum of ergodic  invariant measures for the leafwise geodesic flow on the closure $\overline{L_w}$.

\medskip

   Next, we  consider  some examples where $h(\cGF) > 0$.  
    \begin{prop}
   The Hirsch foliations always have positive geometric entropy.
   \end{prop}
   \proof
 The idea of the proof is as follows. 
    The holonomy pseudogroup $\cGF$ of the Hirsch examples is topologically semi-conjugate  to the pseudogroup  generated by the doubling map $z \mapsto z^2$ on $\mS^1$.
  
  After $\ell$-iterations, the inverse map to $z \mapsto z^{2^{\ell}}$ has derivative of norm $2^{\ell}$ so we have a rough estimate
  $ h(\cGF, \epsilon, \ell) \sim  (2\pi/\epsilon) \cdot 2^{\ell}$. 
  Thus, $h(\cGF) \sim \ln 2$.
  \endproof

For these examples, the relationship between ``orbit growth type''  and   expansion growth type is transparent. Observe that in 
   the Hirsch example,   as we wander out the tree-like leaf, the exponential growth of the ends of the typical leaf yield an exponential growth for the number of $\epsilon$-separated points along the ``core circle'' representing  the transversal space $\cT$. This is suggested by comparing the two illustrations below, where the ends of the left side wrap around the core, with a branching of a pair of pants corresponding to a double covering of the core:

\begin{figure}[H]
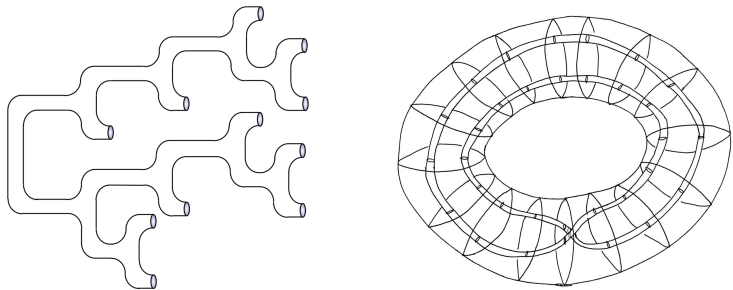

\setlength{\unitlength}{1cm}
\begin{picture}(30,3.8)
\put(3.4,.2){\includegraphics[width=40mm]{figure12.png} \quad \quad \quad\includegraphics[width=45mm]{figure10.png}}
\end{picture}
\caption{Comparing orbit with endset}
\label{Figure13}
\vspace{-20pt}
\end{figure}

\begin{figure}[H]
\setlength{\unitlength}{1cm}
\begin{picture}(30,4)
\put(6,.2){\includegraphics[width=40mm]{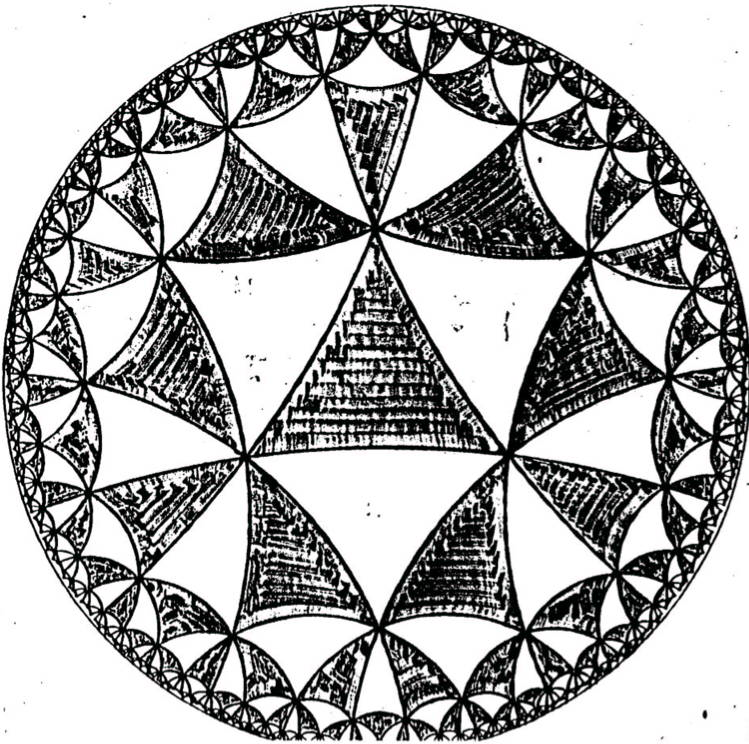}}
\end{picture}
\caption{Tiling by fundamental domains for hyperbolic manifold cover}
\label{Figure14}
\vspace{-20pt}
\end{figure}

It is natural to ask whether there are other classes of foliations for which this phenomenon occurs, that exponential growth type of the leaves is equivalent to positive foliation geometric entropy? It turns out that for the weak stable foliations of Anosov flows, this is also the case in general. First, let us recall a   result of Anthony Manning \cite{Manning1979}:
    
\begin{thm}\label{thm-manning}
Let $B$ be a compact manifold of negative curvature,  let 
  $M = T^1B$ denote the unit tangent bundle to $B$, and
let $\phi_t \colon M \to M$ denote the geodesic flow of $B$.
Then $h_{top}(\phi)= Gr(\pi_1(B, b_0))$. 
 That is, the entropy for the geodesic flow of $B$ equals the  growth rate of the volume of balls in the universal covering of $B$.
 \end{thm}
\proof The idea of proof for this result is conveyed by  the illustration Figure~\ref{Figure14}, representing the fundamental domains for the universal covering.
 The assumption that $B$ has non-positive curvature implies that its universal covering $\widetilde{B}$ is a disk, and we can ``tile'' it with fundamental domains.

From the center basepoint, there is a unique geodesic segment to the corresponding basepoint in each translate. The number of such segments in a given radius is precisely the growth function for the fundamental group $\pi_1(B, b_0)$. 
On the other hand, the negative curvature hypothesis implies that these geodesics separate points for the geodesic flow as well. 
\endproof
    
 We include this example, because it is actually a result about foliation entropy!
The assumption that  $B$ has uniformly negative sectional curvatures implies that the geodesic flow $\phi_t \colon M \to M$ defines a foliation on $M$, its weak-stable foliation. Then by a result of  Pugh and Shub,  the weak-stable manifolds $L_w$ form the leaves of a $C^1$-foliation  of $M$,  called the \emph{weak-stable foliation} for $\phi_t$. Moreover, the orbits of the geodesic flow $\phi_t(w)$ are contained in the leaves of $\F$.
  Then again one has  $h(\cGF) \sim   h_{top}(\phi_1)$ which equals the growth type of the leaves.

 \medskip
  
Besides special cases such as  for the Hirsch foliations and their generalizations in \cite{BHS2006} where one has uniformly expanding holonomy groups, and the weak stable foliations for Anosov flows, how does one determine when a foliation $\F$ has positive entropy? 

There is a third case where $h(\cGF) > 0$ can be concluded, as noted in \cite{GLW1988}, when the   dynamics of $\cGF$ admits a ``ping-pong game''. 
The term ``ping-pong game'' is adopted from the paper \cite{delaHarpe1983} which gives a more geometric proof of Tits Theorem \cite{Tits1972} for the dichotomy of the growth types of countable subgroups of linear groups.  
To say that the dynamics of $\cGF$ admits a ping-pong game, means that there are disjoint open sets $U_0 , U_1 \subset V \subset \cT$ and maps $g_0, g_1 \in \cGF$ such that  for $i = 0,1$:
\begin{itemize}
\item the closure $\overline{V} \subset Dom(g_i)$ for $i = 0,1$
\item $g_i(\overline{V}) \subset U_i$
\end{itemize}
It follows that for each $w \in V$ the forward orbit 
$$\cO^+_{g_0 , g_i}(w) = \{ g_I(w) \mid I = (i_1, \ldots , i_k) ~, i_{\ell} \in \{0,1\} ~ , ~g_I = g_{i_k} \circ \cdot \circ g_{i_1}\}$$
consists of distinct points, and so the full orbit $\cO(w)$ has exponential growth type. Moreover, if $\epsilon>0$ is less than the distance between the disjoint closed subsets $g_0(\overline{V})$ and $g_1(\overline{V})$, then the points in $\cO^+_{g_0 , g_i}(w)$ are all $(\epsilon, \ell)$-separated for appropriate $\ell > 0$, and hence $h(\cGF, K) > 0$.

  For   codimension-one foliations, the existence of ping-pong game dynamics for its pseudogroup  is equivalent to the existence of a ``resilient leaf'', which in turn is analogous to the existence of homoclinic orbits for a diffeomorphism. It is a well-known principle that the existence of homoclinic orbits for a diffeomorphism implies positive topological entropy.

 We conclude with a general question:
 \begin{quest}  
  Are there other canonical classes of $C^1$-foliations  where positive entropy is to be  ``expected''?
  For example, if $\F$ has leaves of exponential growth, when does there   exist a $C^1$-close perturbation of $\F$ with positive entropy?
  \end{quest}

\section{Entropy and Exponent} \label{sec-entropy}

 Three aspects of ``exponential complexity'' for foliation dynamics have been introduced: Lyapunov spectrum for the foliation geodesic flow, exponential growth of orbits, and the geometric entropy which measures  the transverse exponential expansion. In this section,  we discuss the   relationships between these invariants, as is currently understood. 
The theme is summarized by:
\begin{center}   
 {\it Positive Entropy $\leftrightarrow$ Chaotic Dynamics $\leftrightarrow$ ??}
 \end{center}

As always, we assume that   $\F$ is a $C^r$-foliation of a compact manifold $M$, for $r \geq 1$. We formulate three problems illustrating the themes of research:
\begin{prob}\label{prob-ged}
 If $h(\cGF) > 0$, what conclusions can we reach about the dynamics of $\F$?
\end{prob} 
\begin{prob}\label{prob-dge}
What hypotheses on the dynamics of $\F$ are sufficient to imply that $h(\cGF) > 0$?
\end{prob} 
\begin{prob}\label{prob-gecoho}
Are there cohomology hypotheses on $\F$ which would ``improve'' our understanding of its dynamics? How does leafwise cohomology $H^*(\F)$ influence dynamics? How are the secondary cohomology invariants for $\F$ related to entropy?
\end{prob} 

The solution to Problems~\ref{prob-ged} and \ref{prob-dge} are well-known for foliations defined by a $C^2$-flow, due to work of Margulis and Man\'{e} \cite{Mane1987}. The problem with extending these results to foliations of higher dimensions, is that a foliation rarely has any holonomy-invariant measures, and if such exists, there still do not exists methods for estimating recurrence of leaves to the support of the measure, so that the techniques in \cite{Mane1987} do not directly apply.  Thus, given asymptotic data about either the transverse derivative cocycle, or the transverse expansion growth function, one has to develop new techniques to extract from this data dynamical conclusions.

On the other hand, there are examples supporting the hope that such   relationships as suggested in  Problems~\ref{prob-ged}-\ref{prob-gecoho} should exist, and remain to be discovered.  We discuss below some ``deterministic'' techniques, based on the orbit behavior of the foliation geodesic flow which  relate transverse expansion growth with the transverse Lyapunov spectrum of the foliation geodesic flow, and in special cases to the foliation geometric entropy. 
  Another approach, an active area of current research, is to study the relation between exponent and recurrence for ``typical'' orbits of appropriately chosen leafwise harmonic measures \cite{Candel2003, DK2007, DKN2007, DKN2009,Garnett1983a, Garnett1983b}.

 We begin by recalling a   result of Ghys, Langevin, and Walczak \cite{GLW1988} which  gives a straightforward conclusion valid in all codimensions.
 
\begin{thm}  [G-L-W 1988] \label{thm-GLWmeas0}
Let $M$ be compact with a   $C^{1}$-foliation $\F$ of codimension $q \geq 1$, and $X \subset \cT$ a closed subset. 
 If $h(\cGF, X) = 0$,  then the restricted  action of $\cGF$ on $X$ admits an invariant probability measure. 
\end{thm}
The idea of the proof is to interpret the condition $h(\cGF) = 0$ as a type of \emph{equi-distribution} result, and   form averaging sequences over the orbits, which   yield $\cGF$-invariant probability measures on $X$.

\begin{cor}\label{cor-app1}
Let $M$ be compact with a   $C^{1}$-foliation $\F$ of codimension one, and suppose that $\cZ \subset M$ is a minimal set for which the local entropy $h(\cGF , \cZ) = 0$. Then the dynamics of $\cGF$ on $X = \cT \cap \cZ$ is semi-conjugate to the pseudogroup of an isometric dense action on $\mS^1$. If $\F$ is $C^2$, and $M$ is connected, then $\cZ = M$ and the action is conjugate to a rotation group.
\end{cor}
\proof 
 Theorem~\ref{thm-GLWmeas0} implies there exists an invariant probability measure for the action of $\cGF$ on $X$, so the conclusions   follow  from  Sacksteder \cite{Sacksteder1965}.
\endproof

 In the remainder of this section, we discuss three   results of the author on geometric entropy. Note that the works \cite{BW2010,Walczak2010}  by Walczak and Bi\'{s} also study the entropy and orbit growth rates of distal groupoids and group actions.

\begin{thm} \cite{Hurder2000b} \label{thm-main1}
Let $M$ be compact with a   $C^{r}$-foliation $\F$ of codimension-$q$. If $q=1$ and $r \geq 1$, or $q \geq 2$ and $r > 1$, then 
 $$\cGF ~ \text{ distal } \quad \Longrightarrow \quad  h(\cGF) = 0$$
 \end{thm}

\begin{thm} \cite{Hurder2000b}  \label{thm-main2}
 Let $M$ be compact with a  codimension one,  $C^{1}$-foliation $\F$. Then 
  $$h(\cGF) > 0 \quad \Longrightarrow \quad \F ~ \text{ has a resilient leaf}$$
    \end{thm}

\begin{thm}  \cite{HLa2000}  \label{thm-main3}
Let $M$ be compact with a  codimension one  $C^{2}$-foliation $\F$. Then 
  $$0 \ne GV(\F) \in H^3(M, \mR)  \quad \Longrightarrow \quad h(\cGF) > 0$$
   where $GV(\F) \in H^3(M, \mR)$ is the Godbillon-Vey class of $\F$.
   \end{thm}
 
 \medskip
 
The proofs of all three results are based on   the existence of \emph{stable transverse manifolds} for  hyperbolic measures for the foliation geodesic flow. The first step is the following:

  \begin{prop}\label{prop-hypmeasF} 
Let $M$ be compact with a   $C^{1}$-foliation $\F$, and suppose that $\cZ \subset M$ is a minimal set for which the relative entropy $h(\cGF , \cZ) > 0$.  Then, there exists a transversally hyperbolic invariant probability measure $\mu_*$  for the foliation geodesic flow, with   the support of $\mu_*$ contained in the unit leafwise tangent bundle to $\cZ$.
 \end{prop}
 \proof
 We give a sketch of the proof. Let $X = \cZ \cap \cT$. The assumption $\lambda = h(\cGF , X) > 0$ implies there exists $\epsilon > 0$ so that $\lambda_{\epsilon} = h(\cGF, X, \epsilon) > \frac{3}{4} \lambda > 0$. Thus, there exists   a sequence of subsets 
   $\{ \cE_{\ell} \subset X \mid \ell \to \infty \}$ such that $\cE_{\ell}$ is $(\epsilon_{\ell}, r_{\ell})$-separated, where $\epsilon_{\ell} \to 0$ and $r_{\ell} \geq \ell$, and $\# \cE_{\ell} \geq \exp\{3  r_{\ell} \lambda /4\}$. 
   
   We can assume without loss that $\cE_{\ell}$ is contained in the transversal for a single coordinate chart, say $\cE_{\ell} \subset \cT_1$. As $\cT_1$ has bounded diameter, this implies there exists pairs $\{x_{\ell} , y_{\ell} \} \subset \cE_{\ell}$ so that 
   $$d_{\cT}(x_{\ell} , y_{\ell}) \lesssim  \exp\{- 3 r_{\ell} \lambda /4\} \cdot  \text{diam} (\cT_1) $$
 and leafwise geodesic segments    $\sigma_{\ell} \colon [0,1] \to L_{x_{\ell}}$ with $\|\sigma_{\ell} \| \leq r_{\ell}$ such that $d_{\cT}(h_{\sigma_{\ell}}(x_{\ell}) , h_{\sigma_{\ell}}(y_{\ell})) \geq \epsilon$.
 
 By the mean value theorem, there exists $z_{\ell} \in B_{\cT}(x_{\ell} ,  \exp\{- 3 r_{\ell} \lambda /4\} \cdot  \text{diam} (\cT_1) )$ such that 
 $\| D_{z_{\ell}} h_{\sigma_{\ell}} \|   \gtrsim  \exp\{ 3 r_{\ell} \lambda /4\}$.
 
 Noting that  $\epsilon_{\ell} \to 0$ and choosing appropriate subsequences, the resulting geodesic segments $\sigma_{\ell}$ define an invariant probability measure $\mu_*$ for the geodesic flow, with support in $\cZ$. Moreover, by the cocycle equation and continuity of the derivative, the measure $\mu_*$ will be hyperbolic. In fact, with careful choices above, the exponent can be made arbitrarily close to $h(\cGF, X)$, modulo the adjustment for the relation between geodesic and word lengths. See \cite{Hurder2000b} for details. 
 \endproof

 The construction sketched in the proof of Proposition~\ref{prop-hypmeasF}  is very ``lossy'' - at each stage, information about the transverse expansion due to the assumption that $h(\cGF, X ) > 0$       gets discarded, especially in that for each $n$ we only consider a pair of points $(x_{\ell} , y_{\ell})$  to obtain a geodesic segment $\sigma_{\ell}$  along which the transverse derivative has exponentially increasing norm. We will return to this point later.

 The next step in the construction of stable manifolds, is to assume we are given    a transversally hyperbolic invariant probability measure $\mu_*$  for the foliation geodesic flow. Then for a typical point $\whx = (x, \vec{v}) \in \whM$ in the support of  $\mu_*$ the geodesic ray at $(x, \vec{v})$ has an exponentially expanding norm of its transverse derivative, and hence  the normal Lyapunov spectrum of the leafwise geodesic flow on $\whM$ contains a non-trivial expanding  eigenspace. 
  By reversing the time flow (via the inversion $\vec{v} \mapsto - \vec{v}$ of $\whM$) we obtain an  invariant probability measure $\mu_*^-$  for the foliation geodesic flow for which the Lyapunov spectrum of the flow contains a non-trivial contracting eigenspace.

  If we assume that the flow is $C^{1+\alpha}$ for some H\"{o}lder exponent $\alpha > 0$, then there   exists non-trivial stable manifolds in $\whM$ for almost every $(x, \vec{v})$ in the support of $\mu_*^-$.  Denote this stable manifold by $\cS(x, \vec{v})$ and note that  its tangent space projects non-trivially onto $\whQ$.  Moreover,  for points   $\why, \whz \in \cS(x, \vec{v})$  the distance $d(\varphi_t(\why), \varphi_t(\whz))$ converges to $0$ exponentially fast, as $t \to \infty$. Thus, the images $y, z \in M$ of these points converge together exponentially fast   under the holonomy of $\F$.
  
 Combining these results we obtain: 
\begin{thm}\label{thm-stable}
 Let $\F$ be $C^{1+\alpha}$ and suppose that $h(\cGF) > 0$. 
 Then there exists   a transversally hyperbolic invariant probability measure $\mu_*$  for the foliation geodesic flow. 
 Moreover, for a typical point $\whx = (x, \vec{v})$ in the support of $\mu_*^{-}$ there is a transverse stable manifold $\cS(x, \vec{v})$ for the geodesic ray starting at $\whx$. 
 \end{thm}

  If the codimension of $\F$ is one, then the differentiability is just required to be $C^1$, as the stable manifold for $\varphi_t$ consists simply of the full transversal to $\whF$. 
  
  Observe that Theorem~\ref{thm-stable} implies  Theorem~\ref{thm-main1}.  
  
   \medskip

   The assumption that $h(\cGF) > 0$ has much stronger consequences than simply implying that the dynamics of $\cGF$ is not distal, but obtaining these results requires much more care. We sketch next some ideas for analyzing these  dynamics in the case of codimension-one foliations.

   In the proof of Proposition~\ref{prop-hypmeasF}, instead of choosing only a single pair of points $(x_{\ell} , y_{\ell})$ at each stage,  one can also use the Pigeon Hole Principle to choose a subset $\cE_{\ell}' \subset \cE$  contained in a fixed ball $B_{\cT}(w, \delta_{\ell})$ where $\# \cE_{\ell}'$ grows exponentially fast as a function of $\ell$, and the diameter $\delta_{\ell}$ of the ball  decreases exponentially fast, although   at a rates   less that $\lambda_{\epsilon}$.
This  leads to the following  notion.
 
 \begin{defn}\label{def-quiver}
 An $(\epsilon_{\ell}, \delta_{\ell},\ell)$-\emph{quiver} is a subset $\cQ_{\ell} =\{(x_i , \vec{v}_i) \mid 1 \leq i \leq k_{\ell}\} \subset \whM$ such that $x_j \in B_{\cT}(x_i , \delta_{\ell})$ for all $1 \leq j \leq k_{\ell}$, and for the unit-speed geodesic segment $\sigma_i \colon [0,s_i] \to L_{x_i}$ of length   $s_i \leq d$, we have 
 $$d_{\cT}(h_{\sigma_i}(x_i), h_{\sigma_i}(x_j)) \geq \epsilon ~, ~ \text{ for all}~ j \ne i$$
 An \emph{exponential quiver} is a collection of quivers $\{\cQ_{\ell} \mid \ell = 1,2, \ldots\}$ such that the function $ \ell \mapsto \# \cQ_{\ell}$ has exponential growth rate.
 \end{defn}
The idea is that one has a collection of points $\{x_i  \mid 1 \leq i \leq k_{\ell}\}$ contained in a ball of radius $\delta_{\ell}$ along with a corresponding geodesic segment based at each point whose transverse holonomy separates points. The term ``quiver'' is based on the intuitive notion that the collection of geodesic segments emanating from the $\delta_{\ell}$-clustered set of basepoints $\{x_i\}$ is like a collection of arrows in a quiver.
 It is     immediate that  $h(\cGF, \epsilon, d) \geq  \# \cQ_{\ell}$. 
\begin{prop}
If $\F$   admits an exponential quiver, then $h(\cGF) > 0$.
\end{prop}

 For codimension-one foliations, the results of \cite{Hurder1988} and   \cite{LW1994a} yields the converse estimate:
 \begin{prop} 
 Let $\F$ be a $C^1$-foliation of codimension-one on a compact manifold $M$. If $h(\cGF) > 0$ then there exists an  exponential quiver.
  \end{prop}
  
   It is an unresolved question whether a similar result holds for higher codimension. The point is that if so, then $h(\cGF)$ is estimated by the entropy of the foliation geodesic flow, and most of the problems we address here can be resolved using a form of  Pesin Theory for flows relative to the foliation $\F$.  (See \cite{Hurder1988} for further discussion of this point.)
 
 The existence of an exponential quiver for a codimension-one foliation of a compact manifold $M$ has strong implications for its dynamics. The basic idea is that   the basepoints of the geodesic rays in the quiver are tightly clustered, and because the ranges of the endpoints of the geodesic rays are lie in a compact set, one can pass to a subsequence for which the endpoints are also tightly clustered.  From this observation, one can show: 
 
 \begin{thm}\cite{Hurder2000b} \label{thm-quiverppg}
 Let $\F$ be a $C^1$-foliation  with codimension-one foliation of a compact manifold $M$. If $h(\cGF) > 0$, then $\cGF$ acting on $\cT$ admits a   ``ping-pong game'',  which implies the existence of a resilient leaf for $\F$.
\end{thm}

This result is a $C^1$-version of one of the main results concerning the dynamical meaning of positive entropy given in  \cite{GLW1988}. In their paper, Ghys, Langevin, and Walczak  require the foliations be $C^2$,  as they invoke the Poincar\'{e}-Bendixson Theory for codimension-one foliations which is only valid for $C^2$-pseudogroups.

\medskip

There is another approach to obtaining exponential quivers for a foliation $\F$, which is based on cohomology assumptions about $\F$.
For a $C^1$-foliation $\F$, there exists a leafwise closed, continuous $1$-form $\eta$ on $M$ whose cohomology class  $[\eta] \in H^1(M, \F)$
in the leafwise foliated cohomology group is well defined. The form $\eta$ has the property that its integral along a leafwise path  gives the logarithmic infinitesimal expansion of the determinant of the linear holonomy defined by the path.  Thus, for codimension-one foliations, this integral is the expansion exponent of the holonomy.  

For a $C^2$-foliation $\F$ of codimension-$q$ the form $\eta$ can be chosen to be $C^1$, and thus the exterior form $\eta \wedge d\eta^q$ is well defined. As observed by Godbillon and Vey \cite{GV1971}, the form $\eta \wedge d\eta^q$  is closed and yields a well defined cohomology class 
$GV(\F)  \in H^{2q+1}(M, \mR)$. One of the basic problems of   foliation theory has been to understand the ``dynamical meaning'' of this class. A fundamental breakthrough was made by Gerard Duminy in the unpublished manuscripts \cite{Duminy1982a, Duminy1982b}, where the study of this problem ``entered its modern phase''. (See also the reformulation of these results by Cantwell and Conlon in \cite{CantwellConlon1984}.)  Based on this breakthrough, the papers  \cite{HH1984, Hurder1986} showed that if $GV(\F) \ne 0$,  then there is a saturated set of positive measure on which $\eta$ is non-zero, and hence the set of hyperbolic leaves $\HF$ has \emph{positive} Lebesgue measure. This study culminated in the following result of the author with Remi Langevin from \cite{HLa2000}:
 
 \begin{thm}\label{thm-HL2000}
 Let $\F$ be a $C^2$-foliation of codimension-one on a compact manifold. If $\HF$ has positive Lebesgue measure, then $\F$ admits exponential quivers, and in particular the dynamics of $\cGF$ admits ping-pong games. Thus, $h(\cGF) > 0$.
 \end{thm}
 
 Combining Theorem~\ref{thm-HL2000} with the previous remarks yields Theorem~\ref{thm-main3}.
Theorem~\ref{thm-HL2000} is the basis for the somewhat-cryptic Problem~\ref{prob-gecoho} given at the start of this section. The assumption that the class  $[\eta] \in H^1(M, \F)$ is non-trivial on a set of positive Lebesgue measure leads  to positive entropy, and raises the question whether there are other leafwise cohomology classes, possibly of higher degrees,   which if non-trivial have implications for the foliation dynamics.

In general, the results of this section are just part of a more general ``Pesin Theory for foliations'' as sketched in the author's overview paper \cite{Hurder1988}, whose study continues to yield new insights into the dynamical properties of foliations for which $\HF$ is non-empty.
There is much work left to do!


\section{Minimal Sets} \label{sec-minsets}

Every foliation of a compact manifold has at least one minimal set, and possibly a continuum of them. Can they be described? What are their topological properties? When does the dynamics restricted to a minimal set have a ``canonical form''?  Is it possible to give an effective classification of the dynamics of foliation minimal sets, at least for some particular classes?

For non-singular flows, this has been a major theme of research beginning with Poincar\'{e}'s work concerning periodic orbits for flows, and continuing with  the work in the 1960's and 70's of Smale \cite{Smale1967} and others,   to more modern questions about which continua arise as invariant sets for flows \cite{KY1995}. Of notable interest  for foliation theory (in higher codimensions) is    Williams' work on the topology of attractors for Axion A systems \cite{Williams1967, Williams1974}, including the introduction of the so-called Williams solenoids.

 In this section, we discuss the differentiable dynamics  properties of minimal sets, applying the concepts    of the last section.  In  sections~\ref{sec-matchbox} and \ref{sec-shape}, we   generalize this discussion to the classification problem for ``matchbox manifolds''  and  their relevance  to the study  of    foliation dynamics.

   Recall that    a  \emph{minimal set} for $\F$ is a closed, saturated subset   $\cZ \subset M$ for which every leaf $L \subset \cZ $ is dense.  A    \emph{transitive set} for $\F$  is a closed saturated subset   $\cZ \subset M$  such that there exists at least one dense leaf $L_0 \subset \cZ$; that is, the transitive sets are the closures of the leaves. Very little is known in general about the transitive sets for foliations; a well-developed theory for transitive sets,   would include   a generalization of the Poincar\'{e}-Bendixson Theory for codimension-one foliations.

Traditionally, the minimal sets are divided into three classes. A compact leaf of $\F$ is a minimal set.  If every leaf of $\F$ is dense, then $M$ itself is a minimal set. 
The third possibility is that  the minimal set $\cZ $ has no interior, but contains more than one leaf, hence the intersection $\cZ \cap \cT$ is always a perfect set.  This third case can be subdivided into further cases:   if the intersection $\cZ \cap \cT$ is a Cantor set, then $\cZ $ is said to be an \emph{exceptional minimal set}, and otherwise if $\cZ \cap \cT$ has no interior but is not totally disconnected, then it is said to be an  \emph{exotic minimal set}. For codimension one foliations, the case of exotic minimal sets cannot occur, but for foliations with codimension greater than one there are various types of constructions of exotic minimal sets \cite{BHS2006,BisHurder2011a}.

\begin{defn}
An invariant set $\cZ$ is said to be \emph{elliptic} if $\cZ \subset  \EF$.
\end{defn}

For example, if $\F$ is a Riemannian foliation, then all holonomy maps are isometries for some smooth transverse metric. Therefore,  the expansion function $ e(\cGF , T, w)$ defined in Definition~\ref{def-expansionfunction}
 is bounded. It follows that    every minimal set for a Riemannian foliation is elliptic.
\begin{prob}
Does there exist an elliptic   minimal set $\cZ$ for a smooth foliation $\F$, such that  $\cZ$ is not a compact leaf and $\F$ is    not Riemannian in some open neighborhood of $\cZ$?
\end{prob}
No such example has been constructed, to the best of the author's knowledge. Note that as remarked previously, the Denjoy minimal sets are parabolic, but not elliptic.

\begin{defn}
A minimal set $\cZ $ is said to be \emph{parabolic} if $\cZ \subset \EF \cup \PF$, but $\cZ \not\subset \EF$. In particular, $\cZ \cap \HF = \emptyset$.
\end{defn}

Various examples of parabolic minimal sets are known, such as the well-known Denjoy minimal sets for $C^1$-diffeomorphisms in codimension-one. The construction    by  Pat McSwiggen in \cite{McSwiggen1993, McSwiggen1995}  of $C^{k+1-\e}$-diffeomorphisms of $\mT^{k+1}$, uses a generalization of Smale's ``DA'' (\emph{derived from Anosov}) construction  to obtain parabolic minimal sets.

Recall that a compact foliation is one for which every leaf is compact \cite{EV1978, Sullivan1975b, Vogt1976, Vogt1977b, Vogt1994}. 
\begin{prop}
Let $\F$ be a   $C^1$-foliation of a compact manifold $M$, with all leaves of $\F$ compact. Then every leaf of $\F$ is a parabolic minimal set.
 \end{prop}
 \proof
A compact foliation is clearly distal, so  by the proof of Theorem~\ref{thm-main1}, we have $\HF = \emptyset$.
\endproof 
  
  The embedding theorems for solenoids in \cite{ClarkHurder2011a} yield another class of parabolic minimal sets for foliations in arbitrary dimension. 
  
This list of examples   exhaust the   constructions of   parabolic minimal sets of $C^1$-foliations, as known to the author. It would be very interesting to have further constructions.  

 Note that  we have seen previously that $h(\cGF, \cZ ) > 0$ implies $\cZ \cap \HF \ne \emptyset$,  so the parabolic minimal sets include the zero entropy case. Also, a minimal set for a foliation for which $\cGF$ acts distally will be parabolic, so this provides a guide for further constructions.

 \begin{defn}
A minimal set $\cZ $ is said to be \emph{hyperbolic} if  $\cZ \cap \HF \ne \emptyset$.
\end{defn}

As remarked above,     $h(\cGF, \cZ) > 0$ implies that $\cZ$ is hyperbolic, and    by Proposition~\ref{prop-hypmeasF}, there exists a transversally hyperbolic invariant probability measure $\mu_*$  for the foliation geodesic flow restricted to $\cZ $. One of the main open problems in foliation dynamics is to obtain a partial converse to this:

\begin{prob}\label{prob-hyperent}
Let $\F$ be a   $C^{r}$-foliation  of codimension $q \geq 1$ on a compact manifold $M$, and let $\cZ $ be a hyperbolic minimal set. Find conditions on $r \geq 1$, the topology of $\cZ $, and/or the Hausdorff dimension of $\cZ \cap \HF \cap \cT$ which are sufficient to imply that $h(\cGF, \cZ ) > 0$. 
\end{prob}

This is easy to show in a very special case:

\begin{thm}\label{thm-hyperent}
Let $\F$ be a   $C^{2}$-foliation  of codimension $q \geq 1$ on a compact manifold $M$, and let $\cZ $ be a hyperbolic minimal set. If the holonomy of $\cGF$ is conformal, then $h(\cGF, \cZ ) > 0$. 
\end{thm}
\proof
The hyperbolic hypothesis implies that the geodesic flow has a stable manifold for some hyperbolic measure. The conformal hypothesis implies that the holonomy is actually transversally contracting for this measure. That is, the stable manifold for this measure has dimension equal to the transverse dimension of $\F$. Minimality of the dynamics then implies there is a ``ping-pong game'' for the action of $\cGF$ restricted to $\cZ \cap \cT$, and thus $h(\cGF, \cZ ) > 0$. 
\endproof

The difficulty with proving results of the kind in Problem~\ref{prob-hyperent} is that in general, the stable manifolds of the hyperbolic measure for the geodesic flow on $\cZ$ will have dimension less than the codimension of $\F$, and hence the ``trapping'' argument employed above requires some additional hypotheses. Exactly what those hypotheses might be, that is the question.

Note that the construction of Bi\'{s}, Nakayama, and Walczak in \cite{BNW2007} give a $C^0$-foliation with an exotic minimal set $\cZ$ that has $h(\cGF, \cZ ) > 0$. Their technique does not extend to smooth foliations, though perhaps some modification of the method may yield $C^1$-foliations. 

There is another construction  of foliations such that the hypotheses of Theorem~\ref{thm-hyperent} are always satisfied. 
 Let $N$ be a Riemannian manifold of dimension $q$ with metric $d_N$. Let $C \subset N$ be a convex subset for the metric. A diffeomorphism  $f \colon N \to N$ is said to be \emph{contracting on $C$} 
if $f(C) \subset C$ and for all $x,y \in C$, we have $d_N(f(x), f(y)) < d_N(x,y)$.  Then define
\begin{defn}
An iterated function system (IFS) on $N$ is a collection of diffeomorphisms $\{f_1, \ldots, f_k\}$ of $N$ and a compact convex subset $C \subset N$ such that each $f_{\ell}$ is contracting on $C$, and for $\ell \ne \ell'$ we have $f_{\ell}(C) \cap f_{\ell'}(C) = \emptyset$. 
\end{defn}
Note that since $C$ is assumed compact, the contracting assumption implies that for each map $f_{\ell}$ the norm of its differential $Df_{\ell}$ is uniformly less than $1$. That is, the maps $f_{\ell}$ are infinitesimal contractions.

The suspension construction \cite{CN1985,CandelConlon2000} yields a foliation $\F$ on a fiber bundle $M$ over a surface of genus $k$ with fiber $N$, for which the maps $\{f_1, \ldots, f_k\}$ define the holonomy of $\F$. If the manifold $N$ is compact then $M$ will also be compact. 

The relevance of this construction is that such a system admits a   minimal set $\cZ \subset C$,  which is necessarily hyperbolic. In fact, $\cZ$ is the unique minimal set for the restriction of the action to $C$  is called the \emph{Markov Minimal Set} associated to the IFS (see \cite{BisHurder2011a}). It is an exercise to show that  $h(\cGF, \cZ) > 0$ for these examples.

The traditional construction of an IFS is for $N = \mR^q$ and the maps $f_{\ell}$ are assumed to be affine contractions. The compact convex set $C$ can then be chosen to be any sufficiently large closed ball about the origin  in $\mR^q$. There is a vast literature on affine IFS's, as well as   beautiful computer-generated illustrations in articles and books of the   invariant sets for various systems. 

Note that every affine map of $\mR^q$ extends to a conformal map of $\mS^q$, so these constructions also provide examples of hyperbolic minimal sets for smooth foliations of compact manifolds. The construction  of affine minimal sets via this method has many generalizations, and leads to a variety of interesting examples, which can be  considered  from the foliation point of view.

 \section{Classification Schemes}\label{sec-classification}

After introducing several dynamical invariants of $C^1$-foliations, it is time to ask: How to ``classify'' all the foliations of fixed codimension-$q$ on a given closed manifold $M$? 

It all depends on the meaning of the word ``classify'' -- modulo homeomorphism? diffeomorphism? concordance? Borel orbit equivalence? measurable orbit equivalence? These are just some of the notions of equivalence that have been used to approach this issue - see the   surveys \cite{Hurder2002,Hurder2009, Lawson1974, Lawson1977}. We discuss    the role of the invariants introduced   in the previous sections    for the study of this   problem. 

Invariants of  foliation dynamics such as orbit growth type,  transverse expansiveness, or local entropy   are constant on leaves, and thus  are associated in some fashion with  the ``leaf space'' $M/\F$. The question is what notions of equivalence preserve the leaf space $M/\F$ and have enough additional restrictions to preserve these invariants, yet are not so restricted as to be effectively uncomputable.

 Given foliated manifolds $(M_1, \F_1)$ and $(M_2 , \F_2)$, and $r \geq 0$,  the most basic equivalence relation is to be  \emph{$C^r$-conjugate}; that is,  there exists a $C^r$-diffeomorphism $f \colon M_1 \to M_2$ such that the leaves of $\F_1$  are mapped to the leaves of $\F_2$. 
 If $r = 0$, then the map $f$ is just a homeomorphism, and we say the foliations are \emph{topologically conjugate}.
 Certainly, two foliations which are $C^r$-conjugate have ``conjugate leaf spaces''. 
 Most invariants in foliation theory are preserved by $C^1$-conjugation, and some such as leaf growth rate  are preserved by topological conjugation. However, conjugation is an extremely   strong equivalence relation. 

 The introduction of secondary classes for $C^2$-foliations in the 1970's suggested classification modulo ``concordance'', a  weaker form of equivalence than $C^2$-conjugation. Two foliations $\F_1$ and $\F_2$ of codimension-$q$ on a manifold $M$ are \emph{concordant} if there exists a foliation $\F$  on $M \times \mR$,  also  of codimension-$q$, so that $\F$ is transverse to the slices $M \times \{t\}$ for $t = 1,2$, and the restrictions $\F | M \times \{t\} = \F_t$ for $t = 1,2$. The   lecture     notes by Milnor \cite{Milnor2009} and  the survey by Lawson  \cite{Lawson1977} discuss this concept further. 
 
 Concordance is the natural notion of equivalence associated to the study  of homotopy classes of maps from $M$ into a foliation classifying space, such as    $B\Gamma_q^r$ introduced in \cite{Haefliger1970}. 
  The celebrated results by Thurston on classification of foliations are formulated   in terms of homotopy classes of maps into  the classifying spaces $B\Gamma_q^r$. (See \cite{Haefliger1970, Haefliger1984, Thurston1974b, Thurston1976, Tsuboi1989a, Tsuboi1989b} and the surveys \cite{Hurder2009, Lawson1977}.)  
  
  On the other hand, it is unknown if any of the invariants of dynamics discussed in these lectures are preserved, in some fashion, by concordance.  For example, given any two linear foliations of $\mT^2$, they are concordant \cite[Lemma 8.5]{Milnor2009}, so that a foliation whose leaves have linear growth rate can be concordant to one with compact leaves.  There appears to be no  relation between concordance of $\F_1$ and $\F_2$  and some form of equivalence of the   leaf spaces $M/\F_1$ and $M/\F_2$. 
  
  \begin{quest} 
  Given concordant foliations $\F_1$ and $\F_2$ of a compact manifold $M$, does this imply any relationship between their dynamically defined  invariants?
  \end{quest}  
  
  At the other extreme from conjugation,  is the notion of \emph{orbit equivalence} [\emph{OE}]. 
  Recall that the equivalence relation on $\cT$ defined by $\F$ is the Borel subset 
  $$\RF \equiv \{ (w,z) \mid w \in \cT, ~ z \in L_w \cap \cT\} \subset \cT \times \cT$$
Two foliations $\F_1, \F_2$ with complete transversals $\cT_1$ and $\cT_2$, respectively, are  \emph{Borel orbit equivalent} (\emph{bOE}) if there exists a Borel map $h \colon \cT_1 \to \cT_2$ which induces a Borel isomorphism  $\mathcal{R}_{\F_1} \cong \mathcal{R}_{\F_2}$. Note that a Borel orbit equivalence $h \colon \cT_1 \to \cT_2$ induces a Borel ``isomorphism'' $h_* \colon  M_1/\F_1 \to M_2/\F_2$.
If two foliations are topologically conjugate, then they are bOE. On the other hand, the assumption that  $\F_1$ and $\F_2$ are \emph{bOE}, does not imply that their  leaves     have the same dimensions, so this is a much weaker equivalence  than   conjugation.

The foliations $\F_1, \F_2$  are said to be \emph{measurably orbit equivalent} (\emph{mOE}) if there exists a Borel measurable map $h \colon \cT_1 \to \cT_2$ which induces a Borel orbit equivalence, up to sets of Lebesgue measure zero. See the works  \cite{Ellis1969,FeldmanMoore1975,Moore1982,Hjorth2000,KechrisMiller2004} for more background on this topic.

For example, a foliation is said to be (measurably) \emph{hyperfinite} if it is mOE to an action of the integers $\mZ$ on the interval $[0,1]$. The celebrated result of Dye \cite{Dye1959, Dye1963, Krieger1976} implies:
\begin{thm}[Dye 1957] 
A $C^1$-foliation defined by a non-singular flow  is always hyperfinite.
\end{thm}

 Caroline Series  generalized this result in \cite{Series1980a} to foliations whose leaves have polynomial growth.
 
\begin{thm}[Series 1980]  Let $\F$ be a $C^1$-foliation of a compact manifold $M$. 
If the growth type of all functions $Gr(w,\ell)$ are uniformly of polynomial type, then the equivalence relation on $\cT$ defined by $\cGF$ is hyperfinite.
\end{thm}
 
 The most general form of such results is due to   Connes,   Feldman, and   Weiss \cite{CFW1981}, and implies:
\begin{thm}[Connes-Feldman-Weiss 1981] \label{thm-CFW} 
Let $\F$ be a $C^1$-foliation of a compact manifold $M$. If the equivalence relation $\cR_{\F}$ is amenable,  then the equivalence relation on $\cT$ it defines is hyperfinite. In particular, if the growth type of all functions $Gr(w,\ell)$ are uniformly of subexponential  type, 
 then the equivalence relation on $\cT$ it defines is hyperfinite.
\end{thm}

One conclusion of these results is that measurable orbit equivalence does   preserves neither the growth rates of leaves, nor many other ``usual'' invariants of smooth foliations. For example, all ergodic actions of $\mZ^n$ which preserve a probability measure are \emph{mOE} for   all $n \geq 1$,     yet have polynomial orbit growth rates of degree $n$. Also, the weak-stable foliation for a geodesic flow of a closed manifold with constant negative curvature has leaves of exponential growth, has an amenable equivalence relation \cite{Bowen1977}, has positive   entropy \cite{GLW1988}, and   has non-trivial Godbillon-Vey class \cite{GV1971}.

Two foliations $\F_1, \F_2$ on manifolds $M_1$ and $M_2$ with complete transversals $\cT_1$ and $\cT_2$, respectively, are said to be \emph{restricted orbit equivalent} (\emph{rOE}) if there exists a Borel isomorphism  $f \colon M_1 \to M_2$ which maps leaves of $\F_1$ homeomorphically to leaves of $\F_2$, and such that its restriction to transversals induces   a Borel map
$h \colon \cT_1 \to \cT_2$ which induces a Borel isomorphism  $\mathcal{R}_{\F_1} \cong \mathcal{R}_{\F_2}$. Thus, a restricted orbit equivalence ``permutes'' the leaves of the foliations.  If the restriction of such a map induces a quasi-isometry   between the leaves, then we say the foliations are \emph{quasi-isometric orbit equivalent} (\emph{qiOE}). It is then obvious, for example, that the growth rate of a leaf is an   invariant of \emph{qiOE}. It is not known if these   refined notions of equivalence preserve the other invariants.

\begin{quest} \label{quest-qioe+entropy}
Suppose that $C^1$-foliations $\F_1, \F_2$ are  {qiOE}. Does $h(\cG_{\F_1}) > 0$ imply $h(\cG_{\F_2}) > 0$?
\end{quest}

\begin{quest} \label{quest-roe+growth}
Suppose that $C^1$-foliations $\F_1, \F_2$ are  {rOE}. If $L_1 \subset M_1$ is a leaf of $\F_1$,  and $L_2 \subset M_2$ is the corresponding leaf for $\F_2$ under a  {rOE}. Must $L_1$ and $L_2$ have the same growth rates?
\end{quest}

 There are many variants of these questions, whose answers are essentially unknown.  
  These sorts of questions seem  of fundamental importance to the study of foliations. While the topological classification of foliations is surely an unsolvable problem, in any sense of the word ``unsolvable'', a variation on the Borel classification problem might be possible when restricted to special subclasses, such as for foliations with uniformly polynomial growth, or amenable foliations. 
 
  In the late 1970's and early 1980's, 
 Cantwell and Conlon, Hector, Nishimori, Tsuchiya in particular \cite{CantwellConlon1978, CantwellConlon1981b, Hector1983, Tsuchiya1979a, Tsuchiya1980b}, 
 developed a  Poincar\'{e}-Bendixson Theory of levels   for  codimension-one $C^2$-foliations.  
 For real analytic foliations with all leaves of polynomial growth type, their results give an algorithmic description of the limit sets of leaves.    
\begin{prob} \label{prob-roe+growth}
Classify the \emph{restricted orbit equivalence} classes of    codimension-one  real analytic   foliations   with all leaves of polynomial growth type.
\end{prob}

For the general case of codimension-one $C^2$-foliations, the      theory of levels    becomes much more complicated, as there are numerous counter-examples which have been constructed to show that the conclusions in the analytic case do not extend so easily.  The theory of levels  is even more problematic   for  $C^1$-foliations of codimension-one, and  non-existent for foliations of codimension $q > 1$. 

The concept of \emph{measurable  amenable} has a generalization to \emph{amenable Borel equivalence relations}, as given for example by  Anantharaman--Delaroche and  Renault     in  \cite{ADR2000}. The class of foliations in Problem~\ref{prob-roe+growth} is amenable in this sense. Of course, every $1$-dimensional foliation also has this property, and the papers \cite{Furstenberg1963, EThomas1973, Williams1970, GPS1999} give classification schemes for special cases of flows (see also \cite{CHL2011b}). 

\begin{prob} \label{prob-roe+amenable}
Find subclasses of  amenable foliations for which    \emph{restricted orbit equivalence} gives a good classification. 
\end{prob}

The conclusion is that the two notions of  equivalence of foliations  discussed above, concordance and orbit equivalence,   yield classification schemes that are at least somewhat effectively computable, but do not preserve the dynamically defined invariants discussed previously. 

There is another invariant  for measurable equivalence relations, their ``cost'',  as  introduced by Gilbert Levitt \cite{Levitt1995}. The ``cost'' is   \emph{mOE}, essentially by definition. All measurably amenable foliations have cost equal zero, so this invariant does not distinguish a large class of foliation dynamics. On the other hand,  Gaboriau's work in  \cite{Gaboriau2000} showed that ``cost'' is  a very effective invariant of \emph{mOE} for   non-amenable foliations, and has led to spectacular results such as that by  Gaboriau and Popa in \cite{GaboriauPopa2005}.  Other applications of the cost of an equivalence relation can be found in the literature, for example in \cite{AM2001,KechrisMiller2004,Paulin1999}, but further discussion takes us too far away from our theme. 

Bounded cohomology invariants can be used to distinguish  measurable orbit equivalence classes, as in \cite{MonodShalom2006}. As with the cost invariant, these classes vanish for measurable amenable group actions and foliations. On the other hand, the bounded cohomology classes are often non-zero for the same classes of foliations which have non-trivial secondary classes (see \cite{HurderKatok1987}), suggesting their study will have further applications to classifying foliations with exponential complexity.

\begin{prob} \label{prob-roe+ec}
Find classes of  foliations with exponential complexity for which their are non-trivial bounded cohomology invariants.  
\end{prob}

\section{Matchbox Manifolds} \label{sec-matchbox}

 Let $M$ be a foliated manifold, with foliation $\F$.  If $\cS \subset M$ is a closed saturated subset,  then it is an example of a foliated space, as discussed for example in \cite{MS2006},  \cite[Chapter 11]{CandelConlon2000}, or  \cite{ClarkHurder2011a, ClarkHurder2011b}.   
\begin{defn}\label{def-foliatedspace}
$\cS$ is a $C^r$-foliated space if it admits a covering by  foliated coordinate charts $\ds \cU = \left\{ \vp_{i} \colon U_i \to [-1,1]^n \times \fT_i \mid 1 \leq i \leq k  \right\}$ where $\fT_i$ are   compact metric spaces.  The transition functions between overlapping charts are assumed to be $C^r$  along leaves,  for $1 \leq r \leq \infty$,  and the derivatives depend uniform--continuously on the transverse parameter.
\end{defn}
 
In particular, the minimal sets of a foliation $\F$  can be studied ``independently'' as foliated spaces. 
 An  exceptional minimal set  for a foliation  can be considered as  a     \emph{transversally} zero-dimensional foliated space. For flows, these spaces have been    called ``matchbox manifolds'' in the topological dynamics literature \cite{AF1991, AM1988, FO2002}. The author, in the works with Alex Clark and Olga Lukina \cite{ClarkHurder2011a, ClarkHurder2011b, CHL2011b}, propose   the term  matchbox manifold for the more general case: 
 
 \begin{defn}
An $n$-dimensional  \emph{matchbox manifold} $\fM$ is a   continuum (i.e., a connected and compact metrizable space)     which is a smooth $n$-dimensional foliated space with codimension zero. 
\end{defn}
 
 For a matchbox manifold $\fM$, the  \emph{transverse model spaces}  $\fT_i$ in Definition~\ref{def-foliatedspace} are totally disconnected. We define 
 their disjoint union $\fT = \fT_1 \cup \cdots \cup \fT_k$ which is a \emph{total transversal} for   $\F_{\fM}$. 
 
The   path connected components of $\fM$ are precisely the leaves of $\F_{\fM}$, and thus the foliation  of $\fM$ is   defined by the topology.  In particular, any homeomorphism $h \colon \fM \to \fM'$ between two such spaces maps leaves to leaves. We often abuse notation, and refer to $\fM$ implying its foliated structure $\F_{\fM}$. 
 
 We say that $\fM$ is \emph{minimal} if every leaf is dense. In this case, the transverse model spaces $\fT_i$ are Cantor sets, and their disjoint union $\fT$ is again a Cantor set.

 Essentially, the concept of a matchbox manifold   is the same   as that of a lamination, except that matchbox manifolds are not regarded as  embedded in any manifold. In fact, whether a given matchbox manifold $\fM$ embeds as a minimal set of a   foliated    $C^r$-manifold is a fundamental  question. 
 
 The holonomy groupoid $\cG_{\fM}$ is generated as in section~\ref{sec-topdyn}, with object space $\fT$, and  the transition functions $\gamma_{i,j}$ between open subsets of transversal spaces $\fT_i$ and $\fT_j$ defined when the open sets $U_i \cap U_j \ne \emptyset$. By a careful choice of the open covering by foliation charts of $\fM$, we can assume the domains and ranges of the generating maps $\gamma_{i,j}$ are clopen subsets. 
  
\eject   
  
A matchbox manifold is said to be a \emph{suspension}, if there exists:
\begin{itemize}
\item a compact manifold $B_0$ with fundamental group $G_0 = \pi_1(B_0, b_0)$ for some basepoint $b_0 \in B_0$
\item a continuous action $\rho_{\fM} \colon G_0 \to {\bf Homeo}(\fT)$ on a totally disconnected space $\fT$
\item a homeomorphism $\fM \cong \widetilde{B_0} \times_{\rho_{\fM}} \fT$. 
\end{itemize}
The holonomy pseudogroup $\cG_{\fM}$ is then equivalent to that generated by the action  of $G_0$ on $\fT$.

If $B_0 = \mT^n$ so $G_0 = \mZ^n$, then the suspension foliation is defined by an action of $\mR^n$ on $\fM$. 

In general, if $G_0$ is generated by $m > 1$ elements, then the fundamental group of a surface $\Sigma_{2m}$ of genus $2m$ maps onto $G_0$, so the representation $\rho_{\fM}$ lifts to an action of the surface group $\pi_1(\Sigma_{2m}, x_0)$ on $\fT$. Then  the resulting suspension foliation has all leaves isometric to some quotient of the hyperbolic disk. This matchbox manifold has   holonomy groupoid   determined by $\rho_{\fM}$, so the ``general suspension case'' is a $2$-dimensional matchbox manifold with hyperbolic leaves, though the leaves certainly     need not be simply connected.

  Analogous to the case for foliated manifolds, for a matchbox manifold $\fM$ one can define the growth rates of leaves, geometric entropy, and also the foliation geodesic flow. The one missing property is the infinitesimal transverse behavior, as the transverse zero-dimension hypotheses implies there are no transverse vectors. 
This issue will be discussed in  section~\ref{sec-shapedynamics}.

Next, we consider a selection of  examples where  matchbox manifolds arise naturally. The reader will note, that whereas section~\ref{sec-basics}  of these notes introduced some of the simplest examples of foliated manifolds which can be visualized, the examples below  are at the opposite extreme, in that they are essentially impossible to visualize. 

If  $\fM \subset M$ is an exceptional  minimal set in a compact foliated manifold $M$, then with the restricted foliation,  $\fM$  is a matchbox manifold.  For codimension-one foliations, the study of exceptional minimal sets was started in 1960's with work of 
Sacksteder \cite{Sacksteder1964, Sacksteder1965, SackstederSchwartz1965}, and Hector's Thesis   \cite{Hector1972} introduced many of the subsequent themes for their study \cite{CandelConlon2000, CandelConlon2003, CantwellConlon1988a,Hurder1991, HLa2000, Hurder2002, Matsumoto1988, Plante1975, Walczak1996}. The dynamical and topological properties of exceptional minimal sets in higher codimensions are not well-understood. The case of exceptional minimal sets will be discussed further in the next section.

Another source of examples of matchbox manifolds is provided by the
space of tilings associated to a given quasi-periodic tiling $\Delta$ of $\mR^n$.
If   $\Delta$ satisfies the conditions: it is repetitive, aperiodic, and has    finite local complexity, then the ``hull closure'' $\Omega_{\Delta}$ of the translates of $\Delta$ by the action of $\mR^n$ defines a matchbox manifold. These assumptions can be relaxed somewhat, as discussed by Franks and Sadun \cite{FrankSadun2009}.  The tiling space $\Omega_{\Delta}$ was introduced by Bellisard in his study of mathematical models of electron transport \cite{Bellissard1986}. This construction is the subject of many papers, as for example   in   \cite{AP1998, BBG2006, BJS2010, Frank2008, GMPS2010, Hurder1990, Senechal1995}. The results have been generalized to   quasi-periodic tilings of $G$-spaces in \cite{BG2003}. Sadun and Williams \cite{SadunWilliams2003} showed that the space $\Omega_{\Delta}$ associated to a tiling of $\mR^n$ is always a Cantor bundle over $\mT^n$, associated to a minimal free action of $\mZ^n$. A striking result of Marcy Barge and Beverly Diamond \cite{BargeDiamond2001} classifies   $1$-dimensional tiling spaces in terms of   cohomology.

For a few classes of quasi-periodic tilings of $\mR^n$, the codimension one canonical cut and project
tiling spaces ~\cite{FHK2002}, it is known that the associated
matchbox manifold $\Omega_{\Delta}$ is a minimal set for a $C^1$-foliation of a torus
$\mT^{n+1}$, where the foliation is a generalized Denjoy example.

The ``Ghys-Kenyon'' construction, introduced by Ghys in    \cite{Ghys1999},   associates a matchbox manifold to translates of subgraphs of a fixed graph $\cG$. This construction has been studied by E.~Blanc in  \cite{Blanc2001, Blanc2003}, and by F.~Alcalde-Cuesta, A.~Lozano--Rojo, and M.~Macho~Stadler in \cite{ALM2009a, LozanoRojo2009}.  This class of examples provides a wide variety of dynamical behavior, related to the properties of the graph $\cG$. For example, in contrast to the tiling spaces, constructions of Lukina     \cite{Lukina2011,Lukina2014} yield   graph matchbox manifolds which are not   minimal, and can have leaves with non-trivial holonomy.  

Next, we discuss a very  general (and very abstract) procedure for obtaining Cantor bundle examples, which has a variety of    important special cases.  Let $\G$ be a countable group, and choose an integer $ M \geq 1$.  Set $\mN_{m} \equiv \{1,2, \ldots, m\}$ with the discrete topology. Then the product space 
$\ds \Omega_{\G, m} \equiv  \prod_{\gamma \in \G} ~ \mN_{m}$
is compact. For a ``word'' $\omega \in \Omega$,  which is  considered as  a function $\omega \colon \G \to \mN_m$, and for $\delta \in \G$, define  $\delta \cdot \omega(\gamma) = \omega(\gamma \cdot \delta)$. This yields a continuous action of $\G$ on $\ds \Omega_{\G, m}$.

For a word $\omega_0$ let $\Omega_{\omega_0} = \overline{\G \cdot \omega_0}$ denote the closure of the translates of the ``basepoint''  $\omega_0$. Then $\Omega_{\omega_0}$ is compact and totally disconnected, and the action of $\G$ restricts to an action on $\Omega_{\omega_0}$, clearly with a dense orbit. If $\Omega_{\omega_0}$ is minimal and not periodic for the $\G$-action, then it is     expansive.  Otherwise,  it is essentially   impossible to predict   the   dynamical properties of the restricted action of $\G$ on $\Omega_{\omega_0}$.

If $\G$ is finitely generated, then the choice of a Riemann surface $\Sigma$ whose fundamental group maps onto $\G$ yields, via the suspension construction, a $2$-dimensional matchbox manifold $\fM$ whose holonomy pseudogroup is defined by the action of $\G$ on $\Omega_{\omega_0}$.  This construction is clearly related to both of the above constructions, using graphs and using translates of tiles. In these cases, the dynamical properties are related to either the structure of the graph, or the geometry  of the tiling.

For the case where $\G = \mZ^n$ there is an alternate approach to choosing invariant closed subsets of $\Omega_{\sigma} \subset \Omega_{\G, m}$, using translation-invariant pattern rules. When $\Omega_{\sigma}$ is non-empty, this yields   generalized subshifts of finite type, which are called algebraic dynamical systems.     There is an extensive literature on these examples, especially relating their dynamical properties  to commutative algebra and number theory.  For example, the textbooks by Graham Everest and Thomas Ward \cite{EverestWard1999} and Klaus Schmidt \cite{Schmidt1995} give introductions to the dynamics of   algebraically defined actions of $\mZ^n$, and the papers \cite{BoyleLind1997,ELMW2001} lead to the more recent works, following the citations to these papers. 

Finally, we discuss a class of examples of matchbox manifolds, the generalized solenoids, which have a more dynamical origin and geometric interpretations. 
The classical ``Vietoris solenoid'',  introduced in \cite{Vietoris1927}, provides examples of $1$-dimensional matchbox manifolds. Given a sequence of smooth covering maps $p_{\ell} \colon \mS^1 \to \mS^1$ of degree $d_{\ell} > 1$,  form the inverse limit space 
$\ds \cS = \varprojlim~ \{ p_{\ell} \colon \mS^1 \to \mS^1\}$. Then $\cS$ has a smooth flow, whose flow boxes give $\cS$ a matchbox manifold structure. An application of Pontryagin duality \cite{Baer1937, Pontrjagin1934} implies that the space $\cS$ is determined up to foliated homeomorphism by the sequence of integers $\{d_{\ell} \mid \ell = 1,2, \ldots\}$, modulo ``tail equivalence''.

   The existence of 1-dimensional Vietoris solenoids as   minimal sets of smooth flows  has an extensive history in topological dynamics. See for example,    \cite{BowenFranks1976, GT1990, GST1994, Kan1986, MM1980, Smale1967, EThomas1973}. The existence is generally shown via an iterated perturbation argument, which is essentially folklore. That is, starting with a closed orbit, $M_0 \cong \mS^1$, it is modified in an open neighborhood of $M_0$ so that the flow now has a nearby closed orbit $M_1 \cong \mS^1$ which covers $M_0$ with degree $d_1 > 1$. This process is inductively repeated for all subsequent closed orbits $M_{\ell}$ with $\ell > 1$. With suitable care in the choices, the resulting flow    will  be   $C^{\infty}$ and has a minimal set homeomorphic to the inverse limit of the   system of closed orbits resulting from the construction. 

A generalization of the Vietoris solenoid construction was introduced by Bob Williams in \cite{Williams1967,Williams1970}  to describe the topology  of $1$-dimensional   attractor of an Axiom A diffeomorphism $f \colon N \to N$, are again matchbox manifolds.    A $1$-dimensional ``Williams solenoid'' is the inverse limit of  the iterations of a single  expanding map $f \colon B \to B$ of a special form,  where $B$ is a branched $1$-manifold. Williams generalized this construction to higher dimensional branched manifolds in \cite{Williams1974}, which again gives rise to matchbox manifolds.  Farrell and Jones showed in \cite{FarrellJones1980, FarrellJones1981} that bizarre topology can arise if higher dimensions, even in this  special case  where the maps $p_{\ell}$ are dynamically defined.

 Finally, we discuss the class of ``weak solenoids'' introduced by McCord   in \cite{McCord1965}. 
 For $\ell \geq 0$,  let $B_{\ell}$ be compact, orientable connected manifolds without boundary of dimension   $n \geq 1$ . Assume there are given      orientation-preserving, smooth, proper covering maps $\cP = \{ p_{\ell} \colon B_{\ell} \to B_{\ell -1} \mid \ell > 0\}$. 
 Then the    inverse limit topological space
 \begin{equation}\label{eq-invlim}
 \cS_{\cP} \equiv \varprojlim~ \{ p_{\ell} \colon B_{\ell} \to B_{\ell -1}\} ~ \subset ~ \prod_{\ell =0}^{\infty} ~ B_{\ell} \quad \stackrel{\pi_0}{\longrightarrow} \quad B_0
\end{equation}
is said to be a \emph{weak solenoid} with base  $B_0$.  The $p_{\ell}$ are the \emph{bonding maps} for the weak solenoid. Let $\cS$ denote the homeomorphism class of $\cS_{\cP}$, then the 
  collection $\cP$ defining the space $\cS_{\cP}$ is said to be a \emph{presentation} for $\cS$.

  \begin{thm}[McCord \cite{McCord1965}]
   $\cS_{\cP}$ has  a natural structure as an orientable, $n$-dimensional smooth  matchbox manifold, with every leaf dense.
  \end{thm}
The foliated homeomorphism types of weak solenoids  are determined by the algebraic structure of the inverse limit of the   maps on fundamental groups \cite{McCord1965,Schori1966, CordierPorter1989, Mardesic2000, CHL2011a}. These maps are induced by the bonding maps in the given presentation $\cP$, which we consider   in more detail. 

Chose a basepoint $b_0 \in B_0$,   inductively chose   $b_{\ell} \in B_{\ell}$ with $p_{\ell}(b_{\ell}) = b_{\ell -1}$. 
Let $G_{\ell} = \pi_1(B_{\ell} , b_{\ell})$ denote the corresponding fundamental groups. We obtain  a descending chain of groups and injective maps
 $$\cP_{\#} \equiv \left\{  \stackrel{p_{\ell +1}}{\longrightarrow} G_{\ell} \stackrel{p_{\ell}}{\longrightarrow} G_{\ell -1} \stackrel{p_{\ell -1}}{\longrightarrow} \cdots \stackrel{p_{2}}{\longrightarrow}  G_1 \stackrel{p_1}{\longrightarrow}  G_0 \right\}$$

Set $q_{\ell, k} = p_{\ell} \circ \cdots \circ p_{k+1} \colon B_{\ell} \longrightarrow  B_k$.
We say that   $\cS_{\cP}$ is a \emph{McCord solenoid}  if for some fixed $\ell_0 \geq 0$,     for all $\ell \geq \ell_0$ the image 
$(q_{\ell, \ell_0})_{\#} \colon G_{\ell} \to H_{\ell} \subset G_{\ell_0}$ is a normal subgroup of $G_{\ell_0}$. 
Replacing $B_0$ with $B_{\ell_0}$, we can reduce to the case where $\ell_0 = 0$. Then define 
$$\G_{\cP} =  \varprojlim ~ \left\{ G_0/G_{\ell} \to G_0/G_{\ell -1} \right\}$$
which is a Cantor group. Then the space $\cS_{\cP}$ is homeomorphic to the principal $\G_{\cP}$-bundle over $B_0$ defined by the canonical representation $G_0 \to  \G_{\cP}$. 
Thus, the  McCord solenoids are the ``natural'' generalizations of the Vietoris solenoids to higher dimensions.

Note that if the base manifold $B_0$ satisfies $G_0  = \pi_1(B_{0} , b_{0})$ is abelian, then  every weak solenoid over $B_0$ is a McCord solenoid. In particular, when $B_0 \cong \mT^n$ this is the case.

Unlike the case of Vietoris solenoids, very little is known about when an $n$-dimensional weak solenoid is homeomorphic to an exceptional  minimal set for a $C^r$-foliation, for $n \geq 2$ and $r \geq 1$. A discussion of   this question, and some  partial realization results for the case $G_0 \cong \mZ^k$, are given in \cite{ClarkHurder2011a}. 

\begin{prob} Let $\cP$ be a presentation of a weak solenoid $\cS_{\cP}$. Find conditions on $\cP$ such that $\cS_{\cP}$ is foliated homeomorphic to an exceptional minimal set of a $C^r$-foliation, for $r > 1$.
\end{prob}

 \medskip
 
 The   varieties of examples of matchbox manifolds described above shows that they form a large class of interesting foliated spaces,  certainly deserving of further study. We can ask the same questions for matchbox manifolds as for foliations, and foliation  minimal sets: Find invariants of their foliated homeomorphism type, and find   classification schemes for their topological dynamics. 

Note that a $1$-dimensional oriented matchbox manifold is defined by a non-singular flow, and all such examples can be obtained by the suspension of a $\mZ$-action on a $0$-dimensional space \cite{AF1991,AM1988, EThomas1973}. The minimal $1$-dimensional matchbox manifolds thus correspond to   suspensions of   minimal Cantor systems, which have been extensively studied, and even classified up to orbit equivalence and homeomorphism  -- 
see for example \cite{BargeDiamond2001,BKM2005,BM2008,GPS1995,HPS1992}. Thus, the questions we pose below can be considered as asking for extensions of these results from $1$-dimensional flows, to higher dimensions. 

Minimal Cantor systems are classified by the ``full groups'' \cite{HPS1992,GPS1995,GPS1999,BM2008}, which suggests the introduction and study  of  an analogous concept for matchbox manifolds. Define the  closed topological  subgroup   of all leaf-preserving homeomorphisms:
 $${\bf Inner}(\fM, \F_{\fM}) = {\bf  Homeo}(\F_{\fM}) \subset  {\bf  Homeo}(\fM, \F_{\fM})$$ 
 That is, $h \in {\bf Inner}(\fM, \F_{\fM})$ maps each leaf of $\F_{\fM}$ to itself. This is a normal subgroup of ${\bf  Homeo}(\fM, \F_{\fM})$. In analogy with the full group concept, and also group theoretic constructions,   we introduce:
 
 \begin{defn}
 The   group  of \emph{outer automorphisms} of a matchbox manifold $\fM$ is the quotient topological group
 \begin{equation}
{\bf Out}(\fM) = {\bf  Homeo}(\fM, \F_{\fM})/{\bf Inner}(\fM, \F_{\fM})
\end{equation}
 \end{defn}

 One can think of ${\bf Out}(\fM)$ as the group of automorphisms of the leaf space $\fM/\F_{\fM}$ and thus  should reflect many aspects of the space $\fM$ -- its topological, dynamical and algebraic properties. Very little is known, in general, concerning some basic questions in higher dimensions:

\begin{prob} 
Let $\fM$ be a matchbox manifold with foliation $\F_{\fM}$. Study ${\bf Out}(\fM)$:
\begin{enumerate}
\item If ${\bf Out}(\fM)$ is not discrete, must it act transitively?  If not, what are the examples?
\item If ${\bf Out}(\fM)$ is   discrete and infinite, what conditions on $\fM$ imply that it is finitely generated?
\item Suppose that $\fM$ is minimal and expansive, must ${\bf Out}(\fM)$ be discrete?
\item For what hypotheses on $\fM$ must  ${\bf Out}(\fM)$ be a finite group?
\end{enumerate}
\end{prob}

A matchbox manifold $\fM$ is said to be \emph{homogeneous} if the group of homeomorphisms ${\bf Homeo}(\fM)$ of $\fM$ acts transitively. For a matchbox manifold, every homeomorphism is a foliated homeomorphism, so we have ${\bf Homeo}(\fM) = {\bf  Homeo}(\fM, \F_{\fM}) $. A result of Bing \cite{Bing1960} showed that if $\fM$ is a homogeneous matchbox manifold of dimension 1, then $\fM$ is   homeomorphic to a Vietoris solenoid. The higher dimensional versions of this result have been an open problem, with one direction proven by McCord:
 
  \begin{thm}[McCord \cite{McCord1965}]
Let $\fM$ be homeomorphic to  a  McCord solenoid $\cS_{\cP}$. Then $\fM$ is   homogeneous, and the pseudogroup associated to it is  equicontinuous.
  \end{thm}

Results of the author with Alex Clark give a converse to this, which generalizes Bing's Theorem.

\begin{thm}[Clark and Hurder 2010 \cite{ClarkHurder2011b}] \label{thm-homog}
Let $\fM$ be a smooth, oriented matchbox manifold. If the pseudogroup associated to $\fM$ is  equicontinuous, then $\fM$ is minimal, and is   homeomorphic to a weak solenoid. If $\fM$ is homogeneous, then $\fM$ is homeomorphic to a McCord solenoid. 
\end{thm}

That is,    if ${\bf Out}(\fM)$ acts transitively on $\fM/\F_{\fM}$, then $\fM$ is homeomorphic to a McCord solenoid.

\begin{prob} 
For what hypotheses on $\fM$ does  the isomorphism class of ${\bf Out}(\fM)$ characterize the homeomorphism class of $\fM$?
\end{prob} 

 There is an analogy between Theorem~\ref{thm-homog}, and the classification theory for Riemannian foliations \cite{Molino1988,MoerdijkMrcunbook2003}. Recall that a Riemannian foliation $\F$ on a compact manifold $M$ is said to be {\it transversally parallelizable} (or TP) if the group of foliation-preserving diffeomorphisms of $M$ acts transitively. In this case, the minimal sets for $\F$ are principle $H$-bundles, where $H$ is the structural Lie group of the foliation. Theorem~\ref{thm-homog} is the analog of this result for matchbox manifolds. It is interesting to compare this result with the theory of equicontinuous foliations on compact manifolds, as in \cite{ALC2009}. 
 
 However, if $\fM$ is equicontinuous, but not homogeneous, then the analogy becomes more tenuous. Clark, Fokkink, and Lukina introduce in \cite{CFL2010} the Schreier continuum for weak solenoids, an invariant of the topology of $\fM$,  which they use to calculate the end structures of leaves. In particular, they show that there exists   weak solenoids for which the number of ends of leaves can be between 2 and infinity, which is impossible for Riemannian foliations (see also \cite{Ghys1995}). 
 
 The classification of equicontinuous matchbox manifolds implies a classification of weak solenoids, and this appears far from being understood, if not simply impossible \cite{Hjorth2000, KechrisMiller2004, Thomas2001, Thomas2003}.

 As in Definitions~\ref{def-distal}, \ref{def-proximal} and \ref{def-expansive}, one can likewise define distal, proximal, and expansive matchbox manifolds.  Here is a basic question:
 
\begin{prob}
Give an algebraic classification for  minimal expansive matchbox manifolds, analogous to the classification of weak solenoids by the   tower of the fundamental groups $\cP_{\#}$.
 \end{prob}

In the case where $\fM$ is a Cantor bundle associated to a free minimal $\mZ^n$-action, all such actions are   \emph{affable Borel equivalence relations}   by work of 
 Giordano, Matui, Putnam, and Skau  \cite{GPS1995, GMPS2010, Putnam2010}. This concept  generalizes to the Borel category the notion of hyperfinite discussed in section~\ref{sec-classification} above. 
The authors   prove that with the above hypotheses, the equivalence relation associated to $\fM$ is  affable. Again, for the case of minimal $\mZ^n$-actions, it then follows that $\fM$ is classified  up to foliated homeomorphism by the directed K-Theory groups associated to the affable structure \cite{GMPS2010,TomsWinter2009}. 

Following along these lines, one approach to a partial  algebraic classification would be to first show:
 
 \begin{prob}\label{prob-affable}
Let $\fM$ be a Cantor bundle associated to a free minimal  action of a countable amenable group $\G$. Show that the equivalence relation associated to $\fM$ is  affable. 
 \end{prob}

 Another approach to classification, in the special case of $2$-dimensional matchbox manifolds and using the leafwise Euler class, was given by  
   Bermudez and   Hector   in \cite{BermudezHector2006}.

The definition of the geometric entropy for a $C^1$-foliation extends to the pseudogroup associated to a matchbox manifold, except that one does not know a priori that the geometric entropy $h(\cGF, \fM)$ is finite.  None the less, the following extension of a result of  Ghys, Langevin, and Walczak holds:

\begin{thm}  [G-L-W \cite{GLW1988}] \label{thm-GLWmeas1}
Let  $\fM$ be a matchbox manifold with   $h(\cGF, \fM) = 0$,  then the  holonomy  pseudogroup associated to $\fM$ admits an invariant probability measure.  Thus, if $\fM$   does not admit a transverse invariant measure, then   $h(\cGF, \fM) > 0$. 
\end{thm}

 It seems that very little is known about the classification of matchbox manifolds with  $h(\cGF, \fM) > 0$. The works on expansive algebraic dynamical systems cited above provide a source of questions and conjectures about this case.

 After discussing the variety of examples and properties of matchbox manifolds, we introduce  the concept of a ``resolution'' of a foliated space by a matchbox manifold.  

\begin{defn}
Let $\cS$ be a foliated space with foliation $\F_{\cS}$. A resolution for $(\cS, \F_{\cS})$ consists of a matchbox manifold $\fM$ with foliation $\F_{\fM}$ and a foliated continuous surjection $\rho \colon \fM \to \cS$ such that the restriction of $\rho$ to a leaf of $\F_{\fM}$ is a covering of a leaf of $\F_{\cS}$.
\end{defn}

Note that we do not assume that $\cS$ is transversally totally disconnected, so the hypothesis that the map is foliated is required.

If $\cS$ is an exceptional minimal set for a foliation $\F$ of a compact manifold $M$, then $\cS$ equipped with the restricted foliation $\F_{\cS} = \F | \cS$  is a resolution of itself.  There are many further examples.

Let $\F_{\alpha}$ be a foliation of $\mT^{n+1}$ by linear hyperplanes of codimension-one, associated to an injective representation  $\alpha \colon \mZ^n \to \mS^1$. Select a leaf $L_0 \subset \mT^{n+1}$, and apply the ``inflation'' technique as in the construction of the Denjoy examples, to obtain a $C^1$-foliation $\F$ on $\mT^{n+1}$ of codimension-one, which then has a unique exceptional minimal set $\cS \subset \mT^{n+1}$. Let $\fM = \cS$ as above. Then using the collapse map, which is the inverse of inflation, we obtain a resolution $\rho \colon \fM \to \cS \to \mT^{n+1}$.  This example is motivated by a standard technique employed in the study of the spectrum of quasi-crystals, and can be generalized to any linear foliation of a torus with contractible leaves. 

Another example is provided by the ``semi-Markov'' examples of foliations constructed in \cite{BHS2006, BisHurder2011a}, for which  there exists a unique \emph{exotic} minimal set $\cS$. The notation \emph{semi-Markov} refers to the property that in these examples, both the resolving  matchbox manifold $\fM$, and the fibers of the resolution map $\rho \colon \fM \to \cS$, admit   descriptions as Markovian dynamics.  

The following problem thus appears quite interesting:

\begin{quest}\label{quest-resolve}
Which minimal sets, or foliated spaces more generally, admit resolutions?
\end{quest}
 
   Note that if $\rho \colon \fM \to \cS$ is a resolution of a minimal set $\cS \subset M$ for the foliation $\F$ of the compact manifold $M$, and   leaf $L \subset \fM$ is a dense leaf, then   $\rho(L) \subset \cS$ is a dense leaf of $\F$.  One version of Question~\ref{quest-resolve} is to ask, given a leaf $L \subset M$ of a $C^r$-foliation $\F$, under what hypotheses on $\F$ does the closure $\cS_L = \overline{L}$ admit a resolution? A solution to this question, along with a better  understanding of how the the topology and dynamics of matchbox manifolds  behave for resolutions, yields a new approach to the study of   the  $C^r$-embedding problem.

\section{Topological Shape} \label{sec-shape}

Next, we discuss the classical notion of shape for topological spaces, and apply these ideas to minimal sets of foliations.

  The concept of shape for a compact metric space was introduced by Borsuk \cite{Borsuk1968} and ``modern shape theory'' develops   algebraic topology of the shape approximations of spaces 
\cite{Mardesic2000,MardesicSegal2001}. The Conley Index of invariant sets for flows is one traditional application of shape theory to the dynamics of flows. The shape of matchbox manifolds is studied in \cite{CHL2014}.

\begin{defn}\label{def-shape}
    Let $\cZ \subset X$ be a compact subset of a complete metric space $X$. 
The \emph{shape} of $\cZ$ is the equivalence class of any descending chain of   open subsets
 $X \supset V_1 \supset \cdots \supset V_k \supset \cdots \supset \cZ $
 with $\cZ = \bigcap_{k =1}^{\infty} ~ V_k$.
 \end{defn}
 The notion of equivalence referred to in the definition is defined by a ``tower of equivalences'' between such approximating neighborhood systems. The reader is referred to \cite{Mardesic2000,MardesicSegal2001} for details and especially  the subtleties of this definition. One property  of shape theory, is that the  shape of $\cZ$ is independent of the space $X$ and the embedding $\cZ \subset X$. We recall an important  notion:
  
 \begin{defn}\label{def-shapestable}
Let   $\cZ \subset X$ be a    compact subset, and $x_0 \in \cZ$ a fixed basepoint. 
Then  $\cZ$  has \emph{stable} shape if the pointed  inclusions 
 $\ds (V_{k+1}, x_0) \subset (V_{k}, x_0) $
 are homotopy equivalences for all $k \gg 0$. 
 \end{defn}
 
  The \emph{shape fundamental group} of $\cZ$   defined by 
  \begin{equation}
\widehat{\pi}_1(\cZ , x_0) = \varprojlim~ \{\pi_1(V_{k+1} , x_0) \to \pi_1(V_k, x_0)  \}
\end{equation}
 is then well defined.
Note that  if $\cZ$ has stable shape, then for $k \gg 0$ we have $\widehat{\pi}_1(\cZ , x_0) \cong \pi_1(V_k, x_0)$. 

The following example from \cite{ClarkSullivan2004} is perhaps the simplest non-trivial  example of stable shape.  Consider a Denjoy flow on the 2-torus $\mT^2$, obtained by applying inflation to an orbit of the flow, as illustrated in Figure~\ref{Figure15}  below.  Let $\cZ$ be the unique minimal set for the flow. Then $\cZ$ is stable, and is  shape equivalent to the pointed wedge of two circles, $\cZ \cong \mS^1 \vee_{x_0} \mS^1$. Consequently,     $\widehat{\pi}_1(\cZ , x_0) \cong \pi_1(\mS^1 \vee_{x_0} \mS^1 , x_0) \cong \mZ * \mZ$.

\begin{figure}[H]
\setlength{\unitlength}{1cm}
\begin{picture}(16,3.2)
\put(2,.2){\includegraphics[width=55mm]{figure2.png}} 
\put(9,.2){\includegraphics[width=28mm]{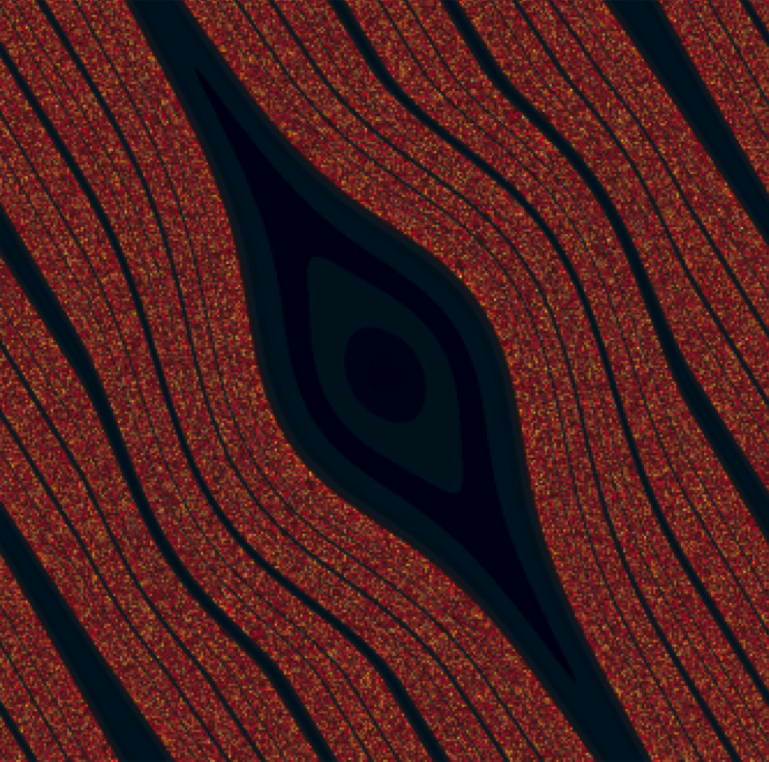}} 
\end{picture}
\caption{Inflating an orbit to obtain a Denjoy flow}
\label{Figure15}
\vspace{-20pt}
\end{figure}

As another example, let $\F$ be a codimension-one foliation with an exceptional minimal set $\fM \subset M$. Then $\fM$ has stable shape if and only if the complement $M - \fM$ consists of a \emph{finite} union of connected open saturated subsets. In the shape framework, one of the long-standing    open problems for codimension-one foliation theory is then:
\begin{prob}\label{prob-stableshape}
 Let $\F$ be a codimension-one, $C^2$-foliation of a compact manifold $M$. Show that     an exceptional minimal set $\fM$ for $\F$ must have stable shape.
 
 More generally, one can ask whether there are other classes of foliations with codimension greater than one, for which the minimal sets are ``expected'' to have stable shape?
 \end{prob}

  Now let $\fM$ be a matchbox manifold, with metric $d_{\fM}$ defining the topology. Choose a basepoint $x_0 \in \fM$ and let $L_0$ be the leaf containing $L_0$.   For $\e > 0$, let    $\tau_{x_0,z} \colon [0,1] \to L_0$ be a leafwise path such that   $d_{\fM}(x_0,z) < \e$.
  Define an equivalence relation on such loops by
 $\tau_0  \stackrel{\e}{\sim}  \tau_1$ if there is a leafwise homotopy $\tau_t$ from $\tau_0$ to $\tau_1$ such that 
 $\tau_t(0) = x_0$ and $d_{\fM}(\tau_t(1) , x_0) < \e$ for all $0 \leq t \leq 1$. 
 The collection of all such approximate loops up to equivalence is denoted by 
\begin{equation}
 \pi_1^{\e}(\fM, x_0) = \{ \widehat{\tau} \mid \widehat{\tau} \stackrel{\e}{\sim} \widehat{\tau'}\}
\end{equation}
 The sets $ \pi_1^{\e}(\fM, x_0)$ do not have a group structure, as concatenation of paths is not necessarily well defined. In any case, there are always maps $\pi_1^{\e'}(\fM, x_0) \subset  \pi_1^{\e}(\fM, x_0)$ for $0 < \e' < \e$.
 
 Note that the sets $\pi_1^{\e}(\fM, x_0)$ may depend strongly on the choice of the basepoint $x_0$.   
  
Next, suppose that $\rho \colon \fM \to \cZ \subset M$ is a resolution of a closed invariant subset $\cZ$ for a foliation $\F$ of a foliated manifold $M$. Let $\delta_0 > 0$ be a Lebesgue number for a covering of $M$ by foliation charts, and let $\e_0 > $ be a modulus of continuity for $\rho$. That is, if $x,y \in \fM$ satisfy $d_{\fM}(x,y) < \e_0$ then $d_M(\rho(x), \rho(y)) < \delta_0$. Set $x_{\rho} = \rho(x_0)$.
 
 \begin{lemma}
 Let $\e < \e_0$ then there is a well defined map $\rho_{\#} \colon \pi_1^{\e}(\fM, x_0)  \to \pi_1(M, x_{\rho})$. 
 \end{lemma}
 \proof
 The assumption  $d_{\fM}(\tau_t(1) , x_0) < \e < \e_0$ for all $0 \leq t \leq 1$ implies that the endpoints satisfy $d_M( \rho(\tau_t(1), x_{\rho}) < \delta_0$, hence are joined by a family of paths contained in a foliation chart. 
 \endproof
 
 There is also a well defined map from the shape fundamental group, $\rho_{\#} \colon  \widehat{\pi}_1(\fM , x_0) \to \pi_1(M, x_{\rho})$. 
 
 For example, suppose the $\fM$ is a McCord solenoid, which is resolution of a minimal set $\cZ \subset M$. Then  the shape of $\fM$ is not stable, and  $\widehat{\pi}_1(\fM , x_0)$ is an non-trivial  inverse limit. Since $\pi_1(M, x_{\rho})$ is always finitely presented, it is a countable group. Thus, the kernel of  $\rho_{\#} \colon  \widehat{\pi}_1(\fM , x_0) \to \pi_1(M, x_{\rho})$  must be non-trivial. 
  We conclude this technical digression with a basic question:
 
  \begin{prob} \label{prob-stablefund}
  Let $\rho \colon \fM \to \cZ$ be a resolution of a closed saturated subset of the foliated manifold $M$.
 How are the   subgroups of $\pi_1(M, x_{\rho})$ given by the images  $\rho_{\#}(\widehat{\pi}_1(\fM , x_0))$ and  $\rho_{\#}(\pi_1^{\e'}(\fM, x_0))$  related to the dynamics of $\F$ and the topology of $M$?
 \end{prob}

It is natural to ask why this problem, and whether these abstract notions have any applications? One point is  that such related ideas have already been introduced in the foliation literature, by Haefliger \cite{Haefliger1988} in his study of Riemannian pseudogroups, and in the study of approximate orbits in foliation dynamics \cite{BW1998,LW1994a, LW1994b, Inaba2000}. We include the above discussion, as the author  believes such considerations are  a fundamental part of the study of shape theory of minimal sets, and these ideas have not been explored. For example, any difference between the subgroups appearing in Problem~\ref{prob-stablefund} will be a measure of how far the set $\cZ$ is from being stable.

  \section{Shape Dynamics}\label{sec-shapedynamics}

 Finally, we introduce the notion of shape dynamics, which is a refinement of the notion of shape for a closed saturated subset $\cZ \subset M$ of the foliated manifold $M$.  The shape dynamics of a foliated space $\cZ$ studies the germinal dynamics of   a sequence of coverings of $\cZ$ which define  its shape, and where the open sets are the union of  foliation charts associated to    a $\G_q^r$-cocycle over the fundamental groupoid of $\F_{\cZ}$. We illustrate this concept with an example.

 Let $\cZ \subset M$ be a closed saturated set. 
 Given $\e > 0$, we can choose a finite covering of $\cZ$ by foliation charts of $M$, whose diameters are bounded above by $\e$. Taking the union of these open sets which intersect $\cZ$, we obtain a shape approximation $\cZ \subset V_{\e} \subset M$. The shape of $\cZ$ can then be defined by the collection of open neighborhoods $\{ V_{\e} \mid \e = 1/\ell ~, ~ \ell = 1, 2, \ldots \}$ for example.

 Associated to each leafwise path $\tau \colon [0,1] \to \cZ$, its    holonomy map $h_{\tau}$ can be defined using a covering of $\cZ$ by foliation charts. In particular, defining the shape approximations of $\cZ$ using foliation charts yields well defined germinal holonomy along all leafwise paths in $\cZ$.   The collection of all such  holonomy maps defines the \emph{shape dynamics} of $\cZ$. 

In terms of sheaf-theoretic approach to foliations of  Haefliger's thesis \cite{Haefliger1958, Haefliger1984}, the foliation $\F$ defines a $\G_q^r$-cocycle over the fundamental groupoid $\G_{\F}$ of $\F$. A closed saturated subset $\cZ \subset M$ induces a subgroupoid  $\G_{\F | \cZ} \subset  \G_{\F}$ given by  the germinal holonomy along leafwise paths of $\F$ which lie in $\cZ$. That is, a shape approximation to a closed saturated subset of $M$ yields more than just the topological shape of $\cZ$, it also yields a $\G_{\F}$-cocycle on the shape approximations. 

Now consider the restriction of the $\G_{\F}$ -cocycle defined by $\F$ to the elements of  $\pi_1^{\e}(\fM, x_0)$. This is well defined, as the germinal holonomy depends only on the leafwise homotopy class of the path. We thus obtain the  holonomy along ``almost closed leafwise paths'', a concept   has a long tradition in foliation folklore. Shape theory simply adds some additional formal structure to their consideration. 
 
  This notion is closely associated to the concept of ``germinal holonomy'' introduced by Timothy Gendron \cite{Gendron2006, Gendron2008}. A related construction has been used by Andr\'{e} Haefliger in his study of the isometry groups associated to the holonomy along a fixed leaf  of a Riemannian foliation  \cite{Haefliger1988}.   
  
  The study of foliation entropy, at its most technical level, often relies on the transformations induced by restriction of the $\G_{\F}$-cocycle   to the elements of  $\pi_1^{\e}(\fM, x_0)$, for $\e> 0$ sufficiently small. This is seen in the works \cite{BW1998,LW1994a, LW1994b, Inaba2000}, and also in the proof of Theorem~\ref{thm-HL2000} in  \cite{HLa2000}. 
  
 Motivated by these examples, we state a very general  problem:
 \begin{prob}
 Given a minimal set $\fM$, what can we say about the ``shape dynamics'' of $\fM$? 
 \end{prob}
   For example, the local entropy $h_{loc}(\cGF, w)$ introduced in Definition~\ref{def-localentropy} is an invariant of the shape dynamics of $\cZ$ with $w \in \cZ$. What other dynamical invariants can be formulated in terms of   shape?

 The $\G_{\F}$-cocycle defined by $\F$ is functorial, so if we are given a resolution $\rho \colon \fM \to \cZ$, then the $\G_{\F}$-cocycle over $\cZ$ lifts to a   $\G_{\F, \rho}$-cocycle  over $\fM$.
 Moreover, the derivative of the   holonomy maps defines a functor $D \colon \G_{\F} \to GL(q, \mR)$, thus a resolution $\rho$ yields a $GL(q, \mR)$-valued cocycle $D \circ \rho$ over the homotopy groupoid of $\fM$. We can then define, exactly as in section~\ref{sec-derivatives}, the normal exponents for the geodesic flow in shape dynamics. 
 
We say that the shape dynamics for $\rho \colon \fM \to \cZ \subset M$ has \emph{hyperbolic type} if  $\rho(\fM) \cap   \HF \ne \emptyset$. The normal cocycle for the leafwise geodesic flow on $\fM$ then has non-zero exponents. What restrictions does this place on the dynamics of $\fM$ and the map $\rho$?
 
  Finally, we reveal the point of our fascination with the formulation of the dynamics of a foliation  in terms of the shape approximations of its closed invariant sets.
   Recall that  the simplicial geometric realization functor (as described for example in \cite{Lawson1977}) yields   classifying map $\nu \colon M \to  BGL(q, \mR) \cong BO(q)$ of the normal bundle to $\F$, and hence induces the universal normal bundle maps $\widehat{\nu} \colon B\G_q^r \to BO(q)$ for all $r \geq 1$.  
  The celebrated Bott Vanishing Theorem \cite{Bott1970} and the very deep works of Tsuboi \cite{Tsuboi1989a, Tsuboi1989b} show that in fact, there is a strong interaction between the degree of differentiability $C^r$, the topology of the classifying map $\widehat{\nu}$, and the dynamics of foliations. One of the deepest open problems of foliation theory is to understand these relationships for $r > 1$. 
 
The functoriality of the    construction of classifying maps implies that if $\rho \colon \fM \to \cZ \subset M$ is a resolution of $\cZ$, then we obtain a universal classifying map $\fH_{\fM} \colon B\G_{\F, \fM} \to B\G_q^r$ which depends only on the shape dynamics of $\fM$.
We can then formulate a very general version of the ``Sullivan Conjecture'' concerning the non-triviality of the Godbillon-Vey classes, extended to the shape dynamics of matchbox manifolds.

\begin{quest}\label{quest-cycles}
How is the homotopy class  of $\fH_{\fM}$ related to the    shape dynamics of $\fM$?
\end{quest}

The point of this problem, is the folklore concept conveyed to the author by  Hans Sah around 1981,  that the topology of the space $B\G_q^r$  is somehow related to algebraic K-Theory invariants of number fields, and the maps $\fH_{\fM}$ represent the sort of generalized cycles for such a theory.  The motivation for this is the celebrated Mather-Thurston Theorem \cite{Mather1975a, Thurston1974a}, which states that the cohomology of the pointed iterated  loop space $\Omega^q B\G_q^r$ is naturally isomorphic to the group cohomology of the group of compactly supported diffeomorphisms of $\mR^q$, so $H^*(\Omega^q B\G_q^r ; \mZ) \cong H^*( {\rm Diff}_c^r(\mR^q); \mZ)$. The point of Question~\ref{quest-cycles}, is to ask whether the ``cycles'' represented by matchbox manifolds resolving a minimal set fit into this scheme, and if so, how the homology classes obtained are related to dynamics in a germinal neighborhood of the minimal set.
(For more on this, see \cite{Hurder2011a, Hurder2011b}.)

 \appendix
 
 \section{Homework}
 
\begin{itemize}
\item[\underline{Monday}:]  Characterize the  transversally hyperbolic invariant probability measures  $\mu_*$ for the foliation geodesic flow of a given foliation. \\
\item[\underline{Tuesday}:] Classify the foliations with subexponential orbit complexity and distal transverse structure.\\
\item[\underline{Wednesday}:]  Find conditions on the geometry of a foliation such that  exponential orbit growth implies positive entropy.\\
\item[\underline{Thursday}:] Find conditions on the Lyapunov spectrum and invariant measures for the geodesic flow which imply positive entropy.\\
\item[\underline{Friday}:]  Characterize the exceptional minimal sets of zero entropy. \\
\item[\underline{Extra Credit}:] Which matchbox manifolds are homeomorphic to an inverse limit of covering maps of branched $n$-manifolds? 
\end{itemize}



\begin{thebibliography}{10}

\bibitem{AF1991}
{J.M.~Aarts and R.~Fokkink},
\newblock {\it The classification of solenoids},
\newblock {\bf Proc. Amer. Math. Soc.}, 111:1161--1163, 1991.

\bibitem{AM1988}
{J.M.~Aarts and M.~Martens},
\newblock {\it Flows on one-dimensional spaces},
\newblock {\bf Fund. Math.}, 131:39--58, 1988.

\bibitem{AA2003}
{E.~Akin and J.~Auslander},
\newblock {\bf Almost periodic sets and subactions in topological dynamics},
\newblock {\bf Proc. Amer. Math. Soc.}, 131:3059--3062 (electronic), 2003.

\bibitem{AAG2008}
{E.~Akin, J.~Auslander and E.~Glasner},
\newblock {\it The topological dynamics of {E}llis actions},
\newblock {\bf Mem. Amer. Math. Soc.}, Vol. 195, 2008.

\bibitem{AM2001}
{F.~Alcalde-Cuesta and M.~Berm{\'u}dez},
\newblock {\it Une remarque sur les relations d'\'equivalence graph\'ees, munies de mesures harmoniques},
\newblock {\bf C. R. Acad. Sci. Paris S\'er. I Math.}, 332:637--640, 2001.

\bibitem{ALM2009a}
{F.~Alcalde-Cuesta, A.~Lozano--Rojo and M.~Macho~Stadler},
\newblock {\it Dynamique transverse de la lamination de {G}hys-{K}enyon},
\newblock {\bf Ast\'erisque}, 323:1--16, 2009.

\bibitem{ALM2009b}
{F.~Alcalde-Cuesta, A.~Lozano--Rojo and M.~Macho~Stadler},
\newblock {\it Affability of {E}uclidean tilings},
\newblock {\bf C. R. Math. Acad. Sci. Paris}, 347:947--952, 2009.

\bibitem{ALM2010}
{F.~Alcalde-Cuesta, A.~Lozano--Rojo and M.~Macho~Stadler},
\newblock {\it Transversely {C}antor laminations as inverse limits},
\newblock {\bf Proc. Amer. Math. Soc.}, 2010.

\bibitem{AR2010}
{F.~Alcalde-Cuesta and A.~Rechtman},
\newblock {\it Minimal {F}{\o}lner foliations are amenable},
\newblock {\bf Discrete and Continuous Dynamical Systems}, to appear. {\bf arXiv:1001.2793}.

\bibitem{ALC2009}
{J.~{\'A}lvarez L{\'o}pez and A.~Candel},
\newblock {\it Equicontinuous foliated spaces},
\newblock {\bf Math. Z.}, 263:725--774, 2009.

\bibitem{ADR2000}
{C.~Anantharaman--Delaroche and J.~Renault},
\newblock {\it Amenable groupoids},
\newblock {\bf Monographies de L'Enseignement Math\'ematique}, Vol. 36, 2000.

\bibitem{AP1998}
{J.~Anderson and I.~Putnam},
\newblock {\it Topological invariants for substitution tilings and their associated {$C\sp *$}-algebras},
\newblock {\bf Ergodic Theory Dyn. Syst.}, 18:509--537, 1998.

\bibitem{Anosov2006}
{D.V.~Anosov},
\newblock {\it Dynamical systems in the 1960s: the hyperbolic revolution},
\newblock In: {\bf Mathematical events of the twentieth century},
\newblock {Springer, Berlin}, 2006: 1--17.

\bibitem{AH1996}
{O.~Attie and S.~Hurder},
\newblock {\it Manifolds which cannot be leaves of foliations},
\newblock {\bf Topology}, 35:335--353, 1996.

\bibitem{Auslander1963}
{L.~Auslander, L.~Green, and F.~Hahn},
\newblock {\it Flows on homogeneous spaces},
\newblock {\bf Annals of Mathematics Studies}, No. 53, 
\newblock {Princeton University Press}, {Princeton, N.J.}, 1963. 

\bibitem{Auslander1988}
{J.~Auslander},
\newblock {\bf Minimal flows and their extensions},
\newblock {North-Holland Mathematics Studies 153}, Mathematical Notes 122, North-Holland, Amsterdam 1988.


\bibitem{Badura2005}
{M.~Badura},
\newblock {\it Realizations of growth types},
\newblock {\bf Ergodic Theory Dynam. Systems}, 25:353--363, 2005.

\bibitem{Baer1937}
{R.~Baer},
\newblock {\it Abelian groups without elements of finite order},
\newblock {\bf Duke Math. Jour.}, 3:68-122, 1937.

\bibitem{BargeDiamond2001}
{M.~Barge and B.~Diamond},
\newblock {\it A complete invariant for the topology of one-dimensional substitution tiling spaces},
\newblock {\bf Ergodic Theory Dynam. Systems}, 21:1333--1358, 2001.
     
\bibitem{BarreiraPesin2002}
{L.~ Barreira and Ya.B.~Pesin}, 
\newblock {\bf Lyapunov exponents and smooth ergodic theory}, 
\newblock {University Lecture Series}, Vol. {23}, {American Mathematical Society, Providence, RI}, {2002}.

\bibitem{Bellissard1986}
{J.~Bellisard},
\newblock  {\it {$K$}-theory of {$C^\ast$}-algebras in solid state physics},
\newblock In: {\bf Statistical mechanics and field theory: mathematical aspects ({G}roningen, 1985)},
\newblock {Lect. Notes in Phys.}, Vol. 257, Springer, Berlin, 1986: 99-156.

\bibitem{BBG2006}
{J.~Bellisard, R.~ Benedetti and J.-M.~Gambaudo},
\newblock {\it Spaces of tilings, finite telescopic approximations and gap-labeling},
\newblock {\bf Comm. Math. Phys.}, 261:1--41, 2006.

\bibitem{BJS2010}
{J.~Bellisard, A.~Julien and J.~Savinien},
\newblock {\it Tiling groupoids and {B}ratteli diagrams},
\newblock {\bf Ann. Henri Poincar\'e}, 11:69--99, 2010.

\bibitem{BG2003}
{R.~Benedetti and J.-M.~Gambaudo},
\newblock {\it On the dynamics of $\mG$-solenoids. {A}pplications to {D}elone sets},
\newblock {\bf Ergodic Theory Dyn. Syst.}, 23:673--691, 2003.

\bibitem{BermudezHector2006}
{M.~Bermudez and G.~Hector},
\newblock {\it Laminations hyperfinies et rev\^etements},
\newblock {\bf Ergodic Theory Dynam. Systems}, 26:305--339, 2006.

\bibitem{BKM2005}
{S.~Bezuglyi, J.~Kwiatkowski and K.~Medynets},
\newblock {\it Approximation in ergodic theory, {B}orel, and {C}antor dynamics},
\newblock In: {\bf Algebraic and topological dynamics}, 
\newblock {Contemp Math. Vol. 385}, American Math. Soc., Providence, RI, 2005: 39--64.

\bibitem{BM2008}
{S.~Bezuglyi  and K.~Medynets},
\newblock {\it Full groups, flip conjugacy, and orbit equivalence of {C}antor minimal systems},
\newblock {\bf Colloq. Math.}, 110:409--429, 2008.
       
\bibitem{Bing1960}
{R.H.~Bing},
\newblock {\it A simple closed curve is the only homogeneous bounded plane continuum that contains an arc},
\newblock {\bf Canad. J. Math.}, 12:209--230, 1960.

\bibitem{BW1998}
{A.~Bi{\'s} and P.~Walczak},
\newblock {\it Pseudo-orbits, pseudoleaves and geometric entropy of foliations},
\newblock {\bf Ergodic Theory Dynam. Systems}, 18:1335--1348, 1998.

\bibitem{BW2010}
{A.~Bi{\'s} and P.~Walczak},
\newblock {\it Entropy of distal groups, pseudogroups, foliations and laminations},
\newblock {\bf Ann. Polon. Math.}, 100:45--54, 2010.

\bibitem{BNW2004}
{A.~Bi{\'s}, H.~ Nakayama and P.~Walczak},
\newblock {\it Locally connected exceptional minimal sets of surface homeomorphisms},
\newblock {\bf Ann. Inst. Fourier (Grenoble)}, 54:711--731, 2004.

\bibitem{BNW2007}
{A.~Bi{\'s}, H.~ Nakayama and P.~ Walczak},
\newblock {\it Modelling minimal foliated spaces with positive entropy},
\newblock {\bf Hokkaido Math. J.}, 36:283--310, 2007.

\bibitem{BHS2006}
{A.~Bi{\'s}, S.~Hurder, and J.~Shive},
\newblock {\it Hirsch foliations in codimension greater than one},
\newblock In: {\bf Foliations 2005},
\newblock {World Scientific Publishing Co. Inc., River Edge, N.J.}, 2006: 71--108.

\bibitem{BisHurder2011a}
{A.~Bi{\'s} and S.~Hurder},
\newblock {\it Markov minimal sets of foliations},
\newblock {in preparation}, 2011.
     
\bibitem{Blanc2001}
{E.~Blanc}, 
\newblock {\it Propri\'{e}t\'{e}es g\'{e}n\'{e}riques des laminations}, 
\newblock {\bf Thesis, Universit'{e} de Claude Bernard-Lyon 1}, Lyon, 2001.

\bibitem{Blanc2003}
{E.~Blanc}, 
\newblock {\it Laminations minimales r\'{e}siduellement \`a 2 bouts}, 
\newblock {\bf Comment. Math. Helv.}, 78:845--864, 2003. 
     
\bibitem{BDV2005}
{C.~Bonatti, L.~D{\'{\i}}az and M.~Viana}
\newblock {\bf Dynamics {B}eyond {U}niform {H}yperbolicity},
\newblock {Encyclopaedia of Mathematical Sciences, Vol. 102},
\newblock {Springer-Verlag, Berlin}, 2005.

\bibitem{Borsuk1968}
{K.~Borsuk},
\newblock {\it Concerning homotopy properties of compacta},
\newblock {\bf Fund. Math.}, 62:223--254, 1968.

\bibitem{Bott1970}
{R.~Bott}, 
\newblock {\it On a topological obstruction to integrability},
\newblock In: {\bf Global Analysis (Proc. Sympos. Pure Math., Vol. XVI, Berkeley, Calif., 1968)},
\newblock {Amer. Math. Soc., Providence, R.I.}, 1970:127--131.

\bibitem{Bowen1971}
{R.~Bowen},
\newblock {\it Entropy for group endomorphisms and homogeneous spaces},
\newblock {\bf Trans. Amer. Math. Soc.}, 153:401--414, 1971.

\bibitem{BowenFranks1976}
{R.~Bowen and J.~Franks},
\newblock {\it The periodic points of maps of the disk and the interval},
\newblock {\bf Topology}, 15:337--342, 1976.

\bibitem{Bowen1977}
{R.~Bowen},
\newblock {\it Anosov foliations are hyperfinite},
\newblock {\bf Ann. of Math. (2)}, 106:549--565, 1977.

\bibitem{BoyleLind1997}
{M.~Boyle and D.~Lind},
\newblock {\it Expansive subdynamics},
\newblock {\bf Trans. Amer. Math. Soc.}, 349:55--102, 1997.

\bibitem{BrinKatok1983}
{M.~Brin and A.~Katok},
\newblock {\it On local entropy},
\newblock In: {\bf Geometric dynamics (Rio de Janeiro, 1981)},
\newblock {Lecture Notes in Math., Vol. 1007}, {Springer, Berlin}, 1983:30--38.

\bibitem{CN1985}
{C.~Camacho and A.~Lins~Neto},
\newblock {\bf Geometric {T}heory of {F}oliations},
\newblock {Translated from the Portuguese by Sue E. Goodman},
\newblock {Progress in Mathematics}, {Birkh\"auser Boston, MA}, 1985.

\bibitem{Candel2003}
{A.~Candel}, 
\newblock {\it  The harmonic measures of Lucy Garnett}, 
\newblock {\bf Advances Math}, 176:187-247, 2003.

\bibitem{CandelConlon2000}
{A.~Candel and L.~Conlon},
\newblock {\bf Foliations {I}},
\newblock Amer. Math. Soc., Providence, RI, 2000.

\bibitem{CandelConlon2003}
{A.~Candel and L.~Conlon},
\newblock {\bf Foliations {II}},
\newblock Amer. Math. Soc., Providence, RI, 2003.

\bibitem{CantwellConlon1978}
{J.~Cantwell and L.~Conlon},
\newblock {\it Growth of leaves},
\newblock {\bf Comment. Math. Helv.}, 53:93--111, 1978.

\bibitem{CantwellConlon1981b}
{J.~Cantwell and L.~Conlon},
\newblock {\it Poincar\'e-Bendixson theory for leaves of codimension one},
\newblock {\bf Transactions Amer. Math. Soc.}, 265:181--209, 1981.

\bibitem{CantwellConlon1983}
{J.~Cantwell and L.~Conlon},
\newblock {\it Analytic foliations and the theory of levels},
\newblock {\bf Math. Annalen}, 265:253--261, 1983.

\bibitem{CantwellConlon1984}
{J.~Cantwell and L.~Conlon},
\newblock {\it The dynamics of open, foliated manifolds and a vanishing theorem for the Godbillon-Vey class},
\newblock {\bf Advances in Math.}, 53:1--27, 1984.

\bibitem{CantwellConlon1988a}
{J.~Cantwell and L.~Conlon},
\newblock {\it Foliations and subshifts},
\newblock {\bf Tohoku Math. J.}, 40:165--187, 1988.

\bibitem{ClarkFokkink2004}
{A.~Clark and R.~Fokkink},
\newblock {\it Embedding solenoids},
\newblock {\bf Fund. Math.}, 181:111--124, 2004.

\bibitem{CFL2010}
{A.~Clark, R.~Fokkink and O.~Lukina},
\newblock {\it The Schreier continuum and ends},
\newblock {\bf Houston J. Math.}   40(2):569-599, 2014; {arXiv:1007.0746v1}.

\bibitem{ClarkSullivan2004}
{A.~Clark and M.~Sullivan},
\newblock {\it The linking homomorphism of one-dimensional minimal sets}, 
\newblock {\bf Topology Appl.}, 141:125--145, 2004.


\bibitem{ClarkHurder2011a}
{A.~Clark and S.~Hurder},
\newblock {\it Embedding solenoids in foliations},
\newblock {\bf Topology Appl.}, 158:1249--1270, 2011.


\bibitem{ClarkHurder2011b}
{A.~Clark and S.~Hurder},
\newblock {\it Homogeneous matchbox manifolds},
\newblock {\bf Trans. Amer. Math. Soc.}, 365:3151--3191, 2013.


\bibitem{CHL2011a}
{A.~Clark, S.~Hurder and O.~ Lukina},
\newblock {\it Shape of matchbox manifolds}, 
\newblock {\bf Indag. Math. (N.S.)}, 25:669--712, 2014.

\bibitem{CHL2011b}
{A.~Clark, S.~Hurder and O.~ Lukina},
\newblock {\it Classifying matchbox manifolds}, 
\newblock {submitted}, 2013; {arXiv:1311.0226}.

\bibitem{CHL2014}
{A.~Clark, S.~Hurder and O.~Lukina},
\newblock {\it {Y}-like matchbox manifolds},
\newblock {in preparation}, 2014.


\bibitem{CFW1981}
{A.~Connes, J.~Feldman and B.~Weiss},
\newblock {\it An amenable equivalence relation is generated by a single transformation},
\newblock {\bf Ergodic Theory Dynamical Systems}, 1:431--450, 1981.

\bibitem{CordierPorter1989}
{J.-M.~Cordier and T.~Porter},
\newblock {\bf Shape theory: Categorical methods of approximation}, 
\newblock {Ellis Horwood Ltd.}, Chichester, 1989. 

\bibitem{delaHarpe1983}
{P.~de~la~Harpe},
\newblock {\it Free groups in linear groups},
\newblock {\bf Enseign. Math. (2)}, 29:129--144, 1983.

\bibitem{Denjoy1932}
{A.~Denjoy}, 
\newblock {\it Sur les courbes d\'efinies par les \'equations diff\'{e}rentielles \`a la surface du tore}, 
\newblock {\bf J. Math. Pures Appl.},11:333--375, 1932.

\bibitem{DK2007}
{B.~Deroin and V.~Klepsyn}, 
\newblock {\it Random conformal dynamical systems}, 
\newblock {\bf Geom. Funct. Anal.}, 17:1043--1105, 2007.

\bibitem{DKN2007}
{B.~Deroin, V.~Klepsyn and A.~Navas}, 
\newblock {\it Sur la dynamique unidimensionnelle en r\'egularit\'e interm\'ediaire}, 
\newblock {\bf Acta Math.}, 199:199--262, 2007.

\bibitem{DKN2009}
{B.~Deroin, V.~Klepsyn and A.~Navas}, 
\newblock {\it On the question of ergodicity for minimal group actions on the circle}, 
\newblock {\bf Mosc. Math. J.}, 9:263--303, 2009.
     
\bibitem{Duminy1982a}
{G.~Duminy},
\newblock {\it L'invariant de Godbillon--Vey d'un feuilletage se localise dans les feuilles ressort},
\newblock {\bf unpublished preprint}, Universit\'e de Lille, I, 1982.

\bibitem{Duminy1982b}
{G.~Duminy},
\newblock {\it Sur les cycles feuillet\'es de codimension un},
\newblock {\bf unpublished manuscript}, Universit\'e de Lille, I, 1982.

\bibitem{Dye1959}
{H.A.~Dye},
\newblock {\it On groups of measure preserving transformation. {I}},
\newblock {\bf Amer. J. Math.}, 81:119--159, 1959.

\bibitem{Dye1963}
{H.A.~Dye},
\newblock {\it On groups of measure preserving transformation. {II}},
\newblock {\bf Amer. J. Math.}, 85:551--576, 1963.

\bibitem{Egashira1993}
{S.~Egashira},
\newblock {\it Expansion growth of foliations},
\newblock {\bf Ann. Fac. Sci. Toulouse Math. (6)}, 2:15--52, 1993.

\bibitem{ELMW2001}
{M.~Einsiedler, D.~Lind, R.~Miles and T.~Ward}, 
\newblock {\it Expansive subdynamics for algebraic {${\mathbb Z}^d$}-actions}, 
\newblock {\bf Ergodic Theory Dynam. Systems}, 21:1695--1729, 2001.

\bibitem{Ellis1969}
{R.~Ellis},
\newblock {\bf Lectures on topological dynamics},
\newblock {W. A. Benjamin, Inc., New York}, 1969.

\bibitem{EEN2001}
{D.~Ellis, R.~Ellis and M.~ Nerurkar},
\newblock {\it The topological dynamics of semigroup actions},
\newblock {\bf Trans. Amer. Math. Soc.}, 353:1279--1320 (electronic), 2001.

\bibitem{EV1978}
{D.B.A.~Epstein and E.~Vogt},  
\newblock {\it A counter-example to the periodic orbit conjecture in codimension 3}
\newblock {\bf   Ann. of Math. (2)},  108:539-552, 1978.

\bibitem{EverestWard1999}
{G.~Everest and T.~Ward},
\newblock {\bf Heights of polynomials and entropy in algebraic dynamics},
\newblock {Universitext}, Springer-Verlag London Ltd., London, 1999.


\bibitem{FarrellJones1980}
{F.T.~Farrell and L.E.~Jones},
\newblock {\it New attractors in hyperbolic dynamics},
\newblock {\bf J. Differential Geom.}, 15:107--133, 1980.
     
\bibitem{FarrellJones1981}
{F.T.~Farrell and L.E.~Jones},
\newblock {\it Expanding immersions on branched manifolds},
\newblock {\bf Amer. J. Math.}, 103:41--101, 1981.
     
\bibitem{FeldmanMoore1975}
{J.~Feldman and {C.C.}~Moore},
\newblock {\it Ergodic equivalence relations, cohomology, and von {N}eumann algebras},
\newblock {\bf Bull. Amer. Math. Soc.}, 81:921--924, 1975.

\bibitem{FO2002}
{F.~Fokkink and L.~Oversteegen},
\newblock {\it Homogeneous weak solenoids},
\newblock {\bf Trans. Amer. Math. Soc.}, 354:3743--3755, 2002.

\bibitem{FHK2002}
{A.~Forrest, J.~Hunton and J.~Kellendonk},
\newblock {\bf Topological invariants for projection method patterns},
\newblock {Mem. Amer. Math. Soc.}, Vol. 159, 2002.

\bibitem{Frank2008}
{N.~Priebe~Frank},
\newblock {\it A primer of substitution tilings of the {E}uclidean plane},
\newblock {\bf Expo. Math.}, 26:295--326, 2008.

\bibitem{FrankSadun2009}
{N.~Priebe~Frank and L.~Sadun},
\newblock {\it Topology of some tiling spaces without finite local complexity},
\newblock {\bf Discrete Contin. Dyn. Syst.}, 23:847--865, 2009.

\bibitem{Furman1999}
{A.~Furman},
\newblock {\it Orbit equivalence rigidity},
\newblock {\bf Ann. of Math. (2)}, 150:1083--1108, 1999.

\bibitem{Furstenberg1963}
{H.~ Furstenberg},
\newblock {\it The structure of distal flows},
\newblock {\bf Amer. J. Math.}, 85:477--515, 1963.

\bibitem{Gaboriau2000}
{D.~Gaboriau}, 
\newblock {\it Co\^ut des relations d'\'equivalence et des groupes}, 
\newblock {\bf Invent. Math.}, 139:41--98, 2000.

\bibitem{GaboriauPopa2005}
{D.~Gaboriau and S.~Popa}, 
\newblock {\it An uncountable family of nonorbit equivalent actions of {$\mathbb F_n$}}, 
\newblock {\bf J. Amer. Math. Soc.}, 18:547--559 (electronic), 2005.

\bibitem{GT1990}
{J.-M.~Gambaudo and C.~Tresser},
\newblock {\it Diffeomorphisms with infinitely many strange attractors},
\newblock {\bf J. Complexity}, 6:409--416, 1990.

\bibitem{GST1994}
{J.-M.~Gambaudo, D.~Sullivan and C.~Tresser},
\newblock {\it Infinite cascades of braids and smooth dynamical systems},
\newblock {\bf Topology}, 33:85--94, 1994.

\bibitem{Garnett1983a}
{L.~Garnett}, 
\newblock {\it Foliations, the ergodic theorem and {B}rownian motion}, 
\newblock {\bf J. Funct. Anal.}, 51(3):285--311, 1983.

\bibitem{Garnett1983b}
{L.~Garnett}, 
\newblock {\it Statistical properties of foliations}, 
\newblock In: {\bf Geometric dynamics (Rio de Janeiro, 1981)}, 
\newblock  {Lecture Notes in Math.}, Vol. 1007,  {Springer--Verlag, Berlin}, 1983, 294--299.

\bibitem{Gendron2006}
{T.M.~Gendron},
\newblock {\it The algebraic theory of the fundamental germ},
\newblock {\bf Bull. Braz. Math. Soc. (N.S.)}, 37:49--87, 2006.

\bibitem{Gendron2008}
{T.M.~Gendron},
\newblock {\it The geometric theory of the fundamental germ},
\newblock {\bf Nagoya Math. J.}, 190:1--34, 2008.

\bibitem{GLW1988}
{\'{E}.~Ghys, R.~Langevin, and P.~Walczak}, 
\newblock {\it Entropie geometrique des feuilletages}, 
\newblock {\bf Acta Math.}, 168:105--142, 1988.

\bibitem{Ghys1995}
{\'E.~Ghys},
\newblock {\it Topologie des feuilles g\'en\'eriques},
\newblock {\bf Ann. of Math. (2)}, 141:387--422, 1995.

\bibitem{Ghys1999}
{\'{E}~Ghys},
\newblock {\it Laminations par surfaces de Riemann},
\newblock In: {\bf Dynamique et G\'{e}om\'{e}trie Complexes},
\newblock {Panoramas \& Synth\`{e}ses}, 8:49--95, 1999.

\bibitem{GhysLey2006}
{\'E.~Ghys and J.~Ley}, 
\newblock {\it Lorenz and Modular Flows: A Visual Introduction}, 
\newblock {\bf Feature Column, Notices AMS}, November 2006.
\newblock {On the web at {\it http://www.ams.org/samplings/feature-column/fcarc-lorenz}}.

\bibitem{GPS1995}
{T.~Giordano, I.~Putnam and C.~Skau},
\newblock {\it Topological orbit equivalence and {$C^*$}-crossed products},
\newblock {\bf J. Reine Angew. Math.}, 469:51--111, 1995.

\bibitem{GPS1999}
{T.~Giordano, I.~Putnam and C.~Skau},
\newblock {\it Full groups of {C}antor minimal systems},
\newblock {\bf Israel J. Math.}, 111:285--320, 1999.
     
\bibitem{GMPS2010}
{T.~Giordano, H.~Matui, I.~Putnam and C.~Skau},
\newblock {\it Orbit equivalence for {C}antor minimal {$\mathbb Z^d$}-systems},
\newblock {\bf Invent. Math.}, 179:119--158, 2010.
     
\bibitem{GlasnerYe2009}
{E.~Glasner and X.~Ye},
\newblock {\it Local entropy theory}, 
\newblock {\bf Ergodic Theory Dynam. Systems}, 29:321--356, 2009.
     
\bibitem{GV1971}
{C.~Godbillon and J.~Vey},
\newblock {\it Un invariant des feuilletages de codimension {$1$}}, 
\newblock {\bf C. R. Acad. Sci. Paris S\'er. A-B}, 273:A92-A95, 1971.

\bibitem{Godbillon1991}
{C.~Godbillon},
\newblock {\bf Feuilletages: Etudes g\'eom\'etriques {I}, {I}{I}},
\newblock {Publ. IRMA Strasbourg (1985-86)};
\newblock {Progress in Math., Vol. 98}, {Birkh\"{a}user, Boston, Mass.}, 1991.

\bibitem{Gromov1981}
{M.~Gromov}, 
\newblock {\it Groups of polynomial growth and expanding maps}, 
\newblock {\bf Inst. Hautes \'Etudes Sci. Publ. Math.}, no. 53:53--73, 1981.

\bibitem{Haefliger1958}
{A.~Haefliger}, 
\newblock {\it Structures feullet{\'{e}}es et cohomologie {\`{a}} valeur dans un faisceau de groupo{\"{i}}des}, 
\newblock {\bf Comment. Math. Helv.}, 32:248--329, 1958.

\bibitem{Haefliger1970}
{A.~Haefliger},
\newblock {\it Feuilletages sur les vari\'et\'es ouvertes},
\newblock {\bf Topology}, 9:183--194, 1970.

\bibitem{Haefliger1984}
{A.~Haefliger},
\newblock {\it Groupo{\"\i}des d'holonomie et classifiants}
\newblock In: {\bf Transversal structure of foliations (Toulouse, 1982)},
\newblock {Asterisque, 177-178, Soci\'et\'e Math\'ematique de France}, 1984:70--97.

\bibitem{Haefliger1988}
{A.~Haefliger},
\newblock {\it Leaf closures in {R}iemannian foliations}
\newblock In: {\bf A f\^ete of topology},
\newblock {Academic Press, Boston, MA}, 1988:3--32.

\bibitem{Haefliger2002a}
{A.~Haefliger},
\newblock {\it Foliations and compactly generated pseudogroups}
\newblock In: {\bf Foliations: geometry and dynamics (Warsaw, 2000)},
\newblock {World Sci. Publ., River Edge, NJ}, 2002:275--295.

\bibitem{Hector1972}
{G.~Hector},
\newblock {\bf Sur un th\'{e}or\`{e}me de structure des feuilletages de codimension un},
\newblock {Thesis, Universit\'{e} Louis Pasteur--Strasbourg}, 1972.

\bibitem{Hector1977}
{G.~Hector},
\newblock {\it Leaves whose growth is neither exponential nor polynomial},
\newblock {\bf Topology}, 16:451--459, 1977.

\bibitem{HecHir1981}
{G.~Hector and U.~Hirsch},
\newblock {\bf Introduction to the Geometry of Foliations, {Parts A,B}},
\newblock Vieweg, Braunschweig, 1981.

\bibitem{Hector1983}
{G.~Hector},
\newblock {\it Architecture des feuilletages de classe {$C\sp{2}$}},
\newblock {\bf Third Schnepfenried geometry conference, Vol. 1 (Schnepfenried, 1982)},
\newblock {Ast\'erisque}, 107, Soci\'et\'e Math\'ematique de France 1983, 243--258.

\bibitem{HH1984}
{J.~Heitsch and S.~Hurder},
\newblock {\it Secondary classes, {Weil} measures and the geometry of foliations},
\newblock {\bf Jour. Differential Geom.}, 20:291--309, 1984.

\bibitem{HPS1992}
{R.~Herman, I.~Putnam and C.~Skau},
\newblock {\it Ordered {B}ratteli diagrams, dimension groups and topological dynamics},
\newblock {\bf Internat. J. Math.}, 3:827--864, 1992.

     
\bibitem{Hirsch1975}
{M.~Hirsch}, 
\newblock {\it A stable analytic foliation with only exceptional minimal sets}, 
\newblock In: {\bf Dynamical Systems, Warwick, 1974}, 
\newblock {Lect. Notes in Math.} vol. 468, , Springer-Verlag, 1975, 9--10.

\bibitem{Hjorth2000}
{G.~Hjorth}, 
\newblock {\it Classification and orbit equivalence relations}, 
\newblock {\bf Mathematical Surveys and Monographs}, Vol. 75, 
\newblock {Amer. Math. Soc.}, Providence, RI, 2000.

\bibitem{HurderKatok1987}
{S.~Hurder and A.~Katok}, 
\newblock {\it Ergodic Theory and {Weil} measures for foliations}, 
\newblock {\bf Annals of Math.}, 126:221--275, 1987.

\bibitem{Hurder1986}
{S.~Hurder},
\newblock {\it The {Godbillon} measure of amenable foliations},
\newblock {\bf Jour. Differential Geom.}, 23:347--365, 1986.

\bibitem{Hurder1988}
{S.~Hurder}, 
\newblock {\it Ergodic theory of foliations and a theorem of {Sacksteder}}, 
\newblock In: {\bf Dynamical Systems: Proceedings, University of Maryland 1986-87}, 
\newblock {Lect. Notes in Math. Vol. 1342}, Springer-Verlag, New York and Berlin, 1988:   291--328. 

\bibitem{Hurder1990}
{S.~Hurder}, 
\newblock {\it Eta invariants and the odd index theorem for coverings}, 
\newblock In: {\bf Geometric and topological invariants of elliptic operators ({B}runswick, {ME}, 1988)}, 
\newblock {Contemp. Math. Vol. 105},  Amer. Math. Soc., Providence, RI,  1990: 47--82. 
      
\bibitem{Hurder1991}
{S.~Hurder},
\newblock {\it Exceptional minimal sets of ${C}^{1+\alpha}$ actions on the circle},
\newblock {\bf Ergodic Theory Dynam. Systems}, 11:455-467, 1991.

\bibitem{Hurder2000b}
{S.~Hurder},
\newblock {\it Entropy and dynamics of $C^1$-foliations},
\newblock {\bf preprint}, August 2000. Available at http://www.math.uic.edu/$\sim$hurder/publications/

 
\bibitem{Hurder2002}
{S.~Hurder},
\newblock {\it Dynamics and the {G}odbillon-{V}ey class: a History and Survey},
\newblock In: {\bf Foliations: Geometry and Dynamics (Warsaw, 2000)},
\newblock {World Scientific Publishing Co. Inc., River Edge, N.J.}, 2002:29--60.

\bibitem{Hurder2009}
{S.~Hurder},
\newblock {\it Classifying foliations},
\newblock In: {\bf Foliations, Geometry and Topology. Paul Schweitzer Festschrift}, 
\newblock {Contemp Math.  Vol. 498}, American Math. Soc., Providence, RI, 2009: 1--65. 




\bibitem{Hurder2010a}
{S.~Hurder},
\newblock {\it Dynamics of Foliations}, {Advanced Course on Foliations}
\newblock {Centre de Recerca Matem\`{a}tica, Barcelona}, May 3--7, 2010.
\newblock {Slides posted at {\it http://www.math.uic.edu/$\sim$hurder/talks/}}

\bibitem{HLa2000}
{S.~Hurder and R.~Langevin},
\newblock {\it Dynamics and the Godbillon-Vey Class of $C^1$ Foliations},
\newblock {\bf submitted}, {arXiv:1403.0494}

\bibitem{HL2014}
{S.~Hurder and O.~Lukina},
\newblock {\it Coarse entropy and transverse dimension},
\newblock {in preparation}, 2014.


\bibitem{Hurder2011a}
{S.~Hurder},
\newblock {\it Shape dynamics and the topology of foliation classifying spaces},
\newblock {in preparation, 2014}.

\bibitem{Hurder2011b}
{S.~Hurder},
\newblock {\it Dynamics and Classification of Foliations},
\newblock {book in preparation, 2014}.

\bibitem{Inaba2000}
{T.~Inaba},
\newblock {\it Expansivity, pseudoleaf tracing property and semistability of foliations},
\newblock {\bf Tokyo J. Math.}, 23:311--323, 2000.

\bibitem{Kaimanovich2001}
{V.~Kaimanovich},
\newblock {\it Equivalence relations with amenable leaves need not be amenable},
\newblock In: {\bf Topology, ergodic theory, real algebraic geometry},
\newblock {Amer. Math. Soc. Transl. Ser. 2} 202:151--166, 2001.

\bibitem{Kan1986}
{I.~Kan},
\newblock {\it Strange attractors of uniform flows},
\newblock {\bf Trans. Amer. Math. Soc.}, 293:135--159, 1986.


\bibitem{KechrisMiller2004}
{A.~Kechris and B.~Miller},
\newblock {\bf Topics in orbit equivalence},
\newblock {Lecture Notes in Math., Vol. 1852}, {Springer, Berlin}, 2004.

\bibitem{KY1995}
{J.~Kennedy and J.~Yorke}, 
\newblock {\it Bizarre topology is natural in dynamical systems}, 
\newblock {\bf Bull. Amer. Math. Soc. (N.S.)}, 32:309--316, 1995.
      
\bibitem{Krieger1976}
{W.~Krieger},
\newblock {\it On ergodic flows and the isomorphism of factors},
\newblock {\bf Math. Ann.} 223:19--70, 1976.


\bibitem{LW1994a}
{R.~Langevin and P.~Walczak}, 
\newblock {\it Entropy, transverse entropy and partitions of unity}, 
\newblock {\bf Ergodic Theory Dynam. Systems}, 14:551--563, 1994.

\bibitem{LW1994b}
{R.~Langevin and P.~Walczak}, 
\newblock {\it Some invariants measuring dynamics of codimension-one foliations}, 
\newblock In: {\bf Geometric study of foliations ({T}okyo, 1993)}, 
\newblock {World Sci. Publ., River Edge, NJ}, 1994: 345--358.
 
 
\bibitem{Lawson1974}
{H.B.~Lawson, Jr.},
\newblock {\it Foliations},
\newblock  {\bf  Bulletin Amer. Math. Soc.}, 8:369--418, 1974.

\bibitem{Lawson1977}
{H.B.~Lawson, Jr.},
\newblock {\bf The Quantitative Theory of Foliations},
\newblock  {NSF Regional Conf. Board   Math. Sci.}, No. 27, 
\newblock {American Math. Society, Providence, RI}, 1977.

     
\bibitem{Levitt1995}
{G.~Levitt}, 
\newblock {\it On the cost of generating an equivalence relation}, 
\newblock {\bf Ergodic Theory Dynam. Systems}, 15:551--563, 1994.
     
\bibitem{Lindenstrauss1999}
{E.~Lindenstrauss},
\newblock {\it Measurable distal and topological distal systems},
\newblock {\bf Ergodic Theory Dynam. Systems}, 19:1173--1181, 1995.

\bibitem{LozanoRojo2009}
{A.~Lozano-Rojo},
\newblock {\it An example of a non-uniquely ergodic lamination}, 
\newblock {preprint}, 2009.

 \bibitem{LRL2013}
{\'A.~Lozano-Rojo and O.~Lukina},
\newblock {\it  Suspensions of {B}ernoulli shifts},
\newblock  {\bf Dyn. Syst.}, 28:551--566, 2013; arXiv:1204.5376.
     

\bibitem{Lukina2011}
{O.~Lukina},
\newblock {\it Hierarchy of graph matchbox manifolds},
\newblock {\bf Topology Appl.}, 159:3461--3485, 2012;  arXiv:1107.5303v3.
     
     
\bibitem{Lukina2014}
{O.~Lukina},
\newblock {\it Hausdorff dimension of graph matchbox manifolds},
\newblock {preprint}, July 2014; {arXiv:1407.0693}.
     
     
     
\bibitem{Mane1987}
{R.~Ma{\~n}{\'e}},
\newblock {\bf Ergodic theory and differentiable dynamics},
\newblock {Ergebnisse der Mathematik und ihrer Grenzgebiete (3)}, Vol. 8, Springer-Verlag, Berlin, 1987.

\bibitem{Manning1979}
{A.~Manning}, 
\newblock {\it Topological entropy for geodesic flows}, 
\newblock {\bf Annals of Math. (2)}, 110:567--573, 1979.

\bibitem{Mardesic2000}
{S.~Marde{\v{s}}i{\'c}},
\newblock {\bf Strong Shape and Homology},
\newblock {Springer Monographs in Mathematics}, Springer-Verlag, Berlin, 2000.

\bibitem{MardesicSegal2001}
{S.~Marde{\v{s}}i{\'c} and J.~Segal}, 
\newblock {\it History of shape theory and its application to general topology}, 
\newblock In: {\bf Handbook of the history of general topology, {V}ol. 3}, 
\newblock {Kluwer Acad. Publ., Dordrecht}, 2001, 1145--1177.

\bibitem{MM1980}
{L.~Markus and K.~Meyer},
\newblock {\it Periodic orbits and solenoids in generic {H}amiltonian dynamical systems},
\newblock {\bf Amer. J. Math.}, 102:25--92, 1980.

\bibitem{Mather1968}
{J.~Mather}, 
\newblock {\it Stability of {$C^{\infty }$} mappings. {I}. {T}he division theorem}, 
\newblock {\bf Annals of Math. (2)}, 87:89--104, 1968.

\bibitem{Mather1975a}
{J.N.~Mather}, 
\newblock {\it Loops and foliations}, 
\newblock In: {\bf Manifolds---Tokyo 1973 (Proc. Internat. Conf., Tokyo, 1973)},  
\newblock {Univ. Tokyo Press, Tokyo}, 1975, 175--180.


\bibitem{Matsumoto1988}
{S.~Matsumoto}, 
\newblock {\it Measure of exceptional minimal sets of codimension one foliations}, 
\newblock In: {\bf A F\^{e}te of Topology}, 
\newblock {Academic Press, Boston}, 1988, 81--94.

\bibitem{Matsumoto2010}
{S.~Matsumoto}, 
\newblock {\it Minimal $C^1$-diffeomorphisms of the circle which admit measurable fundamental domains}, 
\newblock {preprint, arXiv:1005.0585}, May 2010.

\bibitem{McCord1965}
{C.~McCord},
\newblock {\it Inverse limit sequences with covering maps},
\newblock {\bf Trans. Amer. Math. Soc.}, 114:197--209, 1965.

 
\bibitem{McSwiggen1993} 
{P.~McSwiggen}, 
\newblock  {\it Diffeomorphisms of the torus with wandering domains},  
\newblock {\bf Proc. Amer. Math. Soc.}, 117:1175--1186,   1993.

\bibitem{McSwiggen1995}  
{P.~McSwiggen},
\newblock  {\it Diffeomorphisms of the $k$-torus with wandering domains},  
\newblock {\bf Ergodic Theory Dynam. Systems}, 15:1189--1205,  1995.

\bibitem{Milnor1968}
{J.~Milnor}, 
\newblock {\it Curvature and growth of the fundamental group}, 
\newblock {\bf Jour. Differential Geom.}, 2:1--7, 1968.

\bibitem{Milnor2009}
{J.~Milnor},
\newblock {\it Foliations and foliated vector bundles},
\newblock In:{\bf Collected Papers of John Milnor: IV. Homotopy, Homology and Manifolds}, ed. J.~McCleary,
\newblock {American Mathematical Society}, Providence, RI, 2009:  279--320.
\newblock (Updated version of 1969 MIT Lecture Notes.)

\bibitem{MoerdijkMrcunbook2003}
{I.~Moerdijk and J.~Mrcun},
\newblock {\bf Introduction to foliations and Lie groupoids},
\newblock {Cambridge Studies in Advanced Mathematics}, Vol. 91, 2003.

\bibitem{Molino1988}
{P.~Molino},
\newblock {\bf Riemannian foliations},
\newblock {Translated from the French by Grant Cairns, with appendices by Cairns, Y. Carri\`ere, \'E. Ghys, E. Salem and V. Sergiescu},
\newblock {Birkh\"auser Boston Inc., Boston, MA}, 1988.

\bibitem{MonodShalom2006}
{N.~Monod and Y.~Shalom}, 
\newblock {\it Orbit equivalence rigidity and bounded cohomology}, 
\newblock {\bf Ann. of Math. (2)}, 164:825--878, 2006.

\bibitem{Moore1982}
{C.C.~Moore},
\newblock {\it Ergodic theory and von~{N}eumann algebras},
 \newblock In: {\bf Proc. Symp. in Pure Math., Vol 38 Pt. 2} No. 9,
\newblock {American Mathematical Society, Providence, RI}, 1982, 179--226.

\bibitem{MS2006}
{C.C.~Moore and C.~Schochet},
\newblock {\bf Analysis on Foliated Spaces},
\newblock {Math. Sci. Res. Inst. Publ.} Vol. 9, Second Edition,
\newblock Cambridge University Press, New York, 2006.

\bibitem{Oseledets1968}
{V.I.~Oseledets},
\newblock {\it A multiplicative ergodic theorem. Lyapunov characteristic numbers for dynamical systems},
\newblock {\bf Transl. Moscow Math. Soc.}, 19:197--221, 1968.

\bibitem{Ozawa2006}
{N.~Ozawa},
\newblock {\it Amenable actions and applications.},
 \newblock In: {\bf International Congress of Mathematicians. Vol. II}, 1563Ð1580,  
\newblock { Eur. Math. Soc., ZŸrich,}, 2006.

 
\bibitem{Paulin1999}
{F.~Paulin},
\newblock {\it Propri\'et\'es asymptotiques des relations d'\'equivalences mesur\'ees discr\`etes},
\newblock {\bf Markov Process. Related Fields}, 5:163--200, 1999.

\bibitem{Pesin1977}
{Ya.B.~Pesin},
\newblock {\it Characteristic {L}japunov exponents, and smooth ergodic theory},
\newblock {\bf Uspehi Mat. Nauk}, 32:55--112, 1977.
\newblock Transl. {\bf Russian Math. Surveys}, 32(4):55--114, 1977.

\bibitem{PS1981} 
{A.~Phillips and D.~Sullivan}, 
\newblock {\it Geometry of leaves}, 
\newblock {\bf Topology} 20:209--218, 1981.

\bibitem{Plante1975}
{J.~Plante},
\newblock {\it Foliations with measure-preserving holonomy},
\newblock {\bf Ann. of Math.}, 102:327--361, 1975.


\bibitem{Pontrjagin1934}
{L.~Pontryagin},
\newblock {\it The theory of topological commutative groups},
\newblock {\bf Ann. of Math.}, 35:361--388, 1934.

\bibitem{Putnam2010} 
 {I.~Putnam},
\newblock{\it Orbit equivalence of {C}antor minimal systems: a survey and a new proof},
\newblock{\bf Expo. Math.}, 28: 101--131, 2010.

\bibitem{Rechtman2009}
{A.~Rechtman},
\newblock {\it Use and disuse of plugs in foliations},
\newblock {\bf Thesis}, \'{E}cole Normale Sup\'{e}rieure de Lyon, U.M.P.A. , Lyon, France, February 2009.

\bibitem{Richardson2001}
{K.~Richardson},
\newblock {\it The transverse geometry of $G$-manifolds and Riemannian foliations},
\newblock {\bf Illinois J. Math.}, 45:517--535, 2001.

\bibitem{Sacksteder1964}
{R.~Sacksteder},
\newblock {\it On the existence of exceptional leaves in foliations of co-dimension one},
\newblock {\bf Ann. Inst. Fourier (Grenoble)}, 14:221--225, 1964.

\bibitem{Sacksteder1965}
{R.~Sacksteder},
\newblock {\it Foliations and pseudogroups},
\newblock {\bf Amer. J. Math.}, 87:79--102, 1965.

\bibitem{SackstederSchwartz1965}
{R.~Sacksteder and J.~Schwartz},
\newblock {\it Limit sets of foliations},
\newblock {\bf Ann. Inst. Fourier (Grenoble)}, 15:201--213, 1965.

\bibitem{Sadun2003}
{L.~Sadun},
\newblock {\it Tiling spaces are inverse limits},
\newblock {\bf J. Math. Phys.}, 44:5410--5414, 2003.

\bibitem{SadunWilliams2003}
{L.~Sadun and R.F.~Williams},
\newblock {\it Tiling spaces are {C}antor set fiber bundles},
\newblock {\bf Ergodic Theory Dynam. Systems}, 23:307--316, 2003. 

\bibitem{Sadun2008}
{L.~Sadun},
\newblock {\bf Topology of tiling spaces},
\newblock {University Lecture Series}, Vol. 46,  American Mathematical Society, Providence, RI, 2008.

\bibitem{Schmidt1995}
{K.~Schmidt},
\newblock {\bf Dynamical systems of algebraic origin},
\newblock {Progress in Mathematics}, Vol. 128, Birkh\"auser Verlag, Basel, 1995.

\bibitem{Schori1966}
{R.~Schori},
\newblock {\it Inverse limits and homogeneity},
\newblock {\bf Trans. Amer. Math. Soc.}, 124:533--539, 1966.

\bibitem{Schwartzman1957}
{S.~Schwartzman},
\newblock {\it Asymptotic cycles},
\newblock {\bf Ann. of Math. (2)}, 66:270--284, 1957.

\bibitem{Senechal1995}
{M.~Senechal},
\newblock {\bf Quasicrystals and geometry},
\newblock {Cambridge University Press}, Cambridge, 1995.
      
\bibitem{Series1980a}
{C.~Series},
\newblock {\it Foliations of polynomial growth are hyperfinite},
\newblock {\bf Israel J. Math.}, 34:245--258, 1980.

\bibitem{Smale1967}
{S.~Smale}, 
\newblock {\it Differentiable Dynamical Systems}, 
\newblock {\bf Bull. Amer. Math. Soc.}, 73:747-817, 1967.

\bibitem{Sullivan1975b}
{D.~Sullivan},
\newblock {\it A Counterexample to the Periodic Orbit Conjecture},
\newblock {\bf Publ. Math. IHES}, 46:5-14, 1976.

\bibitem{Tamura1992}
{I.~Tamura},
\newblock {\bf Topology of foliations: an introduction}, 
\newblock {Translated from the 1976 Japanese edition and with an afterword by Kiki Hudson, With a foreword by Takashi Tsuboi},
\newblock {American Mathematical Society, Providence, R.I.}, 1992.

\bibitem{EThomas1973}
{E.S.~Thomas, Jr.}
\newblock {\it One-dimensional minimal sets},
\newblock {\bf Topology}, 12:233--242, 1973.

\bibitem{Thomas2001}
{S.~Thomas},
\newblock {\it On the complexity of the classification problem for torsion-free abelian groups of finite rank},
\newblock {\bf Bull. Symbolic Logic}, 7:329--344, 2001.

\bibitem{Thomas2003}
{S.~Thomas},
\newblock {\it The classification problem for torsion-free abelian groups of finite rank},
\newblock {\bf J. Amer. Math. Soc.}, 16:233--258 , 2003.

\bibitem{Thurston1974a}
{W.P.~Thurston},
\newblock {\it Foliations and groups of diffeomorphisms}, 
\newblock {\bf Bull. Amer. Math. Soc.} 80:304--307, 1974. 

\bibitem{Thurston1974b}
{W.P.~Thurston},
\newblock {\it The theory of foliations of codimension greater than one}, 
\newblock {\bf Comment. Math. Helv.} 49:214--231, 1974.

\bibitem{Thurston1976}
{W.P.~Thurston},
\newblock {\it Existence of codimension-one foliations}, 
\newblock {\bf Ann. of Math. (2)} 104: 249--268, 1976.

\bibitem{Tits1972}
{J.~Tits},
\newblock {\it Free subgroups in linear groups},
\newblock {\bf J. of Algebra}, 20:250-270, 1972.

\bibitem{TomsWinter2009}
{A.~Toms and W.~Winter},
\newblock {\it Minimal dynamics and the classification of {$C^*$}-algebras},
\newblock {\bf Proc. Natl. Acad. Sci. USA}, 106:16942--16943, 2009.


\bibitem{Tsuboi1989a}
{T.~Tsuboi},
\newblock {\it On the foliated products of class {$C\sp 1$}}, 
\newblock {\bf Ann. of Math. (2)}, 130:227--271, 1989.

\bibitem{Tsuboi1989b}
{T.~Tsuboi},
\newblock {\it On the connectivity of the classifying spaces for foliations}, 
\newblock In: {\bf Algebraic topology (Evanston, IL, 1988)},
\newblock {Contemp. Math. Vol. 96}, {Amer. Math. Soc., Providence, R.I.}, 1989: 319--331.

\bibitem{Tsuboi2009}
{T.~Tsuboi},
\newblock {\it Classifying spaces for groupoid structures},
\newblock In: {\bf Foliations, Geometry and Topology. Paul Schweitzer Festschrift}, 
\newblock (eds. Nicoalu Saldanha et al),
\newblock {Contemp Math.  Vol. 498}, American Math. Soc., Providence, RI, 2009: 67--81.
      
      
\bibitem{Tsuchiya1979a}
{N.~Tsuchiya},
\newblock {\it Growth and depth of leaves},
\newblock {\bf J. Fac. Sci. Univ. Tokyo Sect. IA Math.}, 26:465--471, 1979.

\bibitem{Tsuchiya1979b}
{N.~Tsuchiya},
\newblock {\it Lower semicontinuity of growth of leaves},
\newblock {\bf J. Fac. Sci. Univ. Tokyo Sect. IA Math.}, 26:473--500, 1979.

\bibitem{Tsuchiya1980b}
{N.~Tsuchiya},
\newblock {\it Leaves of finite depth},
\newblock {\bf Japan. J. Math. (N.S.)}, 6:343--364, 1980.

\bibitem{Veech1970}
{W.~Veech},
\newblock {\it Point-distal flows},
\newblock {\bf Amer. J. Math.}, 92:205--242, 1970.

\bibitem{Vietoris1927}
{L.~ Vietoris},
\newblock {\it  \"{U}ber den h\"oheren {Z}usammenhang kompakter {R}\"aume und eine {K}lasse von zusammenhangstreuen {A}bbildungen},
\newblock {\bf Math. Ann.}, 97:454--472, 1927.

\bibitem{Vogt1976}
{E.~Vogt},
\newblock {\it  Foliations of codimension {$2$} with all leaves compact}, 
\newblock {\bf Manuscripta Math.}, 18:187--212, 1976.

\bibitem{Vogt1977b}
{E.~Vogt},
\newblock {\it  A periodic flow with infinite {Epstein} hierarchy}, 
\newblock {\bf Manuscripta Math.}, 22:403--412, 1977.

\bibitem{Vogt1994}
{E.~Vogt}, 
\newblock {\it  Bad sets of compact foliations of codimension {$2$}},
\newblock In: {\bf Low-dimensional topology (Knoxville, TN, 1992)}, 
\newblock {Conf. Proc. Lecture Notes Geom. Topology, III},
\newblock  {Internat. Press,  Cambridge, MA},1994:187--216.


\bibitem{Walczak1988}
{P.~Walczak},
\newblock {\it Dynamics of the geodesic flow of a foliation},
\newblock {\bf Ergodic Theory Dynam. Systems}, 8:637--650, 1988.

\bibitem{Walczak1996}
{P.~Walczak},
\newblock {\it Hausdorff dimension of Markov invariant sets},
\newblock {\bf Journal Math. Soc. of Japan}, 48:125--133, 1996.

\bibitem{Walczak2004}
{P.~Walczak},
\newblock {\bf Dynamics of foliations, groups and pseudogroups},
\newblock {Instytut Matematyczny Polskiej Akademii Nauk. Monografie Matematyczne (New Series) [Mathematics Institute of the Polish Academy of Sciences. Mathematical Monographs (New Series)]}, Vol. 64.
\newblock {Birkh\"auser Verlag}, Basel, 2004.

\bibitem{Walczak2010}
{P.~Walczak},
\newblock {\it Expansion growth, entropy and invariant measures of distal groups and pseudogroups of homeo- and diffeomorphisms},
\newblock {\bf Discrete Contin. Dyn. Syst.}, 33:4731--4742, 2013.
     
     
\bibitem{Wall1971}
{C.T.C.~Wall}, 
\newblock {\it Lectures on {$C^{\infty}$}-stability and classification}, 
\newblock in {\bf Proceedings of {L}iverpool {S}ingularities--{S}ymposium, {I} (1969/70)}, 
\newblock {Lect. Notes in Math.} vol. 192, , Springer-Verlag, 1971, 178--206.
      
\bibitem{Williams1967}
{R.F.~Williams},
\newblock {\it One-dimensional non-wandering sets},
\newblock {\bf Topology}, 6:473--487, 1967.

\bibitem{Williams1970}
{R.F.~Williams},
\newblock {\it Classification of one dimensional attractors}, 
\newblock In: {\bf Global {A}nalysis ({P}roc. {S}ympos. {P}ure {M}ath., {V}ol.{XIV}, {B}erkeley, {C}alif., 1968)},   
\newblock {American Mathematical Society}, Providence, RI, 1970: 341--361.

\bibitem{Williams1974}
{R.F.~Williams},
\newblock {\it Expanding attractors},
\newblock {\bf Inst. Hautes \'Etudes Sci. Publ. Math.}, 43:169--203, 1974.


\end{thebibliography}
\end{document}